# GAME THEORY, MAXIMUM ENTROPY, MINIMUM DISCREPANCY AND ROBUST BAYESIAN DECISION THEORY[1]

By Peter D. Grünwald and A. Philip Dawid

*CWI Amsterdam and University College London*

We describe and develop a close relationship between two problems that have customarily been regarded as distinct: that of maximizing entropy, and that of minimizing worst-case expected loss. Using a formulation grounded in the equilibrium theory of zero-sum games between Decision Maker and Nature, these two problems are shown to be dual to each other, the solution to each providing that to the other. Although Topsøe described this connection for the Shannon entropy over 20 years ago, it does not appear to be widely known even in that important special case.

We here generalize this theory to apply to arbitrary decision problems and loss functions. We indicate how an appropriate generalized definition of entropy can be associated with such a problem, and we show that, subject to certain regularity conditions, the above-mentioned duality continues to apply in this extended context. This simultaneously provides a possible rationale for maximizing entropy and a tool for finding robust Bayes acts. We also describe the essential identity between the problem of maximizing entropy and that of minimizing a related discrepancy or divergence between distributions. This leads to an extension, to arbitrary discrepancies, of a well-known minimax theorem for the case of Kullback–Leibler divergence (the "redundancy-capacity theorem" of information theory).

Received February 2002; revised May 2003.

[1]Supported in part by the EU Fourth Framework BRA NeuroCOLT II Working Group EP 27150, the European Science Foundation Programme on Highly Structured Stochastic Systems, Eurandom and the Gatsby Charitable Foundation. A four-page abstract containing an overview of part of this paper appeared in the *Proceedings of the 2002 IEEE Information Theory Workshop* [see Grünwald and Dawid (2002)].

*AMS 2000 subject classifications.* Primary 62C20; secondary 94A17.

*Key words and phrases.* Additive model, Bayes act, Bregman divergence, Brier score, convexity, duality, equalizer rule, exponential family, Gamma-minimax, generalized exponential family, Kullback–Leibler divergence, logarithmic score, maximin, mean-value constraints, minimax, mutual information, Pythagorean property, redundancy-capacity theorem, relative entropy, saddle-point, scoring rule, specific entropy, uncertainty function, zero–one loss.

For the important case of families of distributions having certain mean values specified, we develop simple sufficient conditions and methods for identifying the desired solutions. We use this theory to introduce a new concept of "generalized exponential family" linked to the specific decision problem under consideration, and we demonstrate that this shares many of the properties of standard exponential families.

Finally, we show that the existence of an equilibrium in our game can be rephrased in terms of a "Pythagorean property" of the related divergence, thus generalizing previously announced results for Kullback–Leibler and Bregman divergences.

**1. Introduction.** Suppose that, for purposes of inductive inference or choosing an optimal decision, we wish to select a single distribution $P^*$ to act as representative of a class $\Gamma$ of such distributions. The maximum entropy principle [Jaynes (1989), Csiszár (1991) and Kapur and Kesavan (1992)] is widely applied for this purpose, but its rationale has often been controversial [see, e.g., van Fraassen (1981), Shimony (1985), Skyrms (1985), Jaynes (1985), Seidenfeld (1986) and Uffink (1995, 1996)]. Here we emphasize and generalize a reinterpretation of the maximum entropy principle [Topsøe (1979), Walley (1991), Chapter 5, Section 12, and Grünwald (1998)]: that the distribution $P^*$ that maximizes the entropy over $\Gamma$ also minimizes the worst-case expected logarithmic score (log loss). In the terminology of decision theory [Berger (1985)], $P^*$ is a *robust Bayes*, or $\Gamma$-*minimax*, act, when loss is measured by the logarithmic score. This gives a decision-theoretic interpretation of maximum entropy.

In this paper we extend this result to apply to a generalized concept of entropy, tailored to whatever loss function $L$ is regarded as appropriate, not just logarithmic score. We show that, under regularity conditions, maximizing this generalized entropy constitutes the major step toward finding the robust Bayes ("$\Gamma$-minimax") act against $\Gamma$ with respect to $L$. For the important special case that $\Gamma$ is described by mean-value constraints, we give theorems that in many cases allow us to find the maximum generalized entropy distribution explicitly. We further define *generalized exponential families* of distributions, which, for the case of the logarithmic score, reduce to the usual exponential families. We extend generalized entropy to *generalized relative entropy* and show how this is essentially the same as a general decision-theoretic definition of *discrepancy*. We show that the family of divergences between probability measures known as *Bregman divergences* constitutes a special case of such discrepancies. A discrepancy can also be used as a loss function in its own right: we show that a minimax result for relative entropy [Haussler (1997)] can be extended to this more general case. We further show that a "Pythagorean property" [Csiszár (1991)] known to hold for relative entropy and for Bregman divergences in fact applies much more generally; and we give a precise characterization of those discrepancies for which it holds.

Our analysis is game-theoretic, a crucial concern being the existence and properties of a *saddle-point*, and its associated minimax and maximin acts, in a suitable zero-sum game between Decision Maker and Nature.

1.1. *A word of caution.* It is not our purpose either to advocate or to criticize the maximum entropy or robust Bayes approach: we adopt a philosophically neutral stance. Rather, our aim is mathematical unification. By generalizing the concept of entropy beyond the standard Shannon framework, we obtain a variety of interesting characterizations of maximum generalized entropy and display its connections with other known concepts and results.

The connection with $\Gamma$-minimax might be viewed, by those who already regard robust Bayes as a well-founded principle, as a justification for maximizing entropy—but it should be noted that $\Gamma$-minimax, like all minimax approaches, is not without problems of its own [Berger (1985)]. We must also point out that some of the more problematic aspects of maximum entropy inference, such as the incompatibility of maximum entropy with Bayesian updating [Seidenfeld (1986) and Uffink (1996)], carry over to our generalized setting: in the words of one referee, rather than resolving this problem, we "spread it to a new level of abstraction and generality." Although these dangers must be firmly held in mind when considering the implications of this work for inductive inference, they do not undermine the mathematical connections established.

**2. Overview.** We start with an overview of our results. For ease of exposition, we make several simplifying assumptions, such as a finite sample space, in this section. These assumptions will later be relaxed.

2.1. *Maximum entropy and game theory.* Let $\mathcal{X}$ be a finite sample space, and let $\Gamma$ be a family of distributions over $\mathcal{X}$. Consider a Decision Maker (DM) who has to make a decision whose consequences will depend on the outcome of a random variable $X$ defined on $\mathcal{X}$. DM is willing to assume that $X$ is distributed according to some $P \in \Gamma$, a known family of distributions over $\mathcal{X}$, but he or she does not know which such distribution applies. DM would like to pick a single $P^* \in \Gamma$ to base decisions on. One way of selecting such a $P^*$ is to apply the *maximum entropy principle* [Jaynes (1989)], which advises DM to pick that distribution $P^* \in \Gamma$ maximizing $H(P)$ over all $P \in \Gamma$. Here $H(P)$ denotes the *Shannon entropy* of $P$, $H(P) := -\sum_{x \in \mathcal{X}} p(x) \log p(x) = \mathrm{E}_P\{-\log p(X)\}$, where $p$ is the probability mass function of $P$. However, the various rationales offered in support of this advice have often been unclear or disputed. Here we shall present a game-theoretic rationale, which some may find attractive.

Let $\mathcal{A}$ be the set of all probability mass functions defined over $\mathcal{X}$. By the information inequality [Cover and Thomas (1991)], we have that, for any distribution $P$, $\inf_{q \in \mathcal{A}} \mathrm{E}_P\{-\log q(X)\}$ is achieved uniquely at $q = p$, where it takes the value $H(P)$. That is, $H(P) = \inf_{q \in \mathcal{A}} \mathrm{E}_P\{-\log q(X)\}$, and so the maximum entropy can be written as

$$(1) \qquad \sup_{P \in \Gamma} H(P) = \sup_{P \in \Gamma} \inf_{q \in \mathcal{A}} \mathrm{E}_P\{-\log q(X)\}.$$

Now consider the "log loss game" [Good (1952)], in which DM has to specify some $q \in \mathcal{A}$, and DM's ensuing loss if Nature then reveals $X = x$ is measured by $-\log q(x)$. Alternatively, we can consider the "code-length game" [Topsøe (1979) and Harremoës and Topsøe (2001)], wherein we require DM to specify a prefix-free code $\sigma$, mapping $\mathcal{X}$ into a suitable set of finite binary strings, and to measure his or her loss when $X = x$ by the length $\kappa(x)$ of the codeword $\sigma(x)$. Thus DM's objective is to minimize expected code-length. Basic results of coding theory [see, e.g., Dawid (1992)] imply that we can associate with $\sigma$ a probability mass function $q$ having $q(x) = 2^{-\kappa(x)}$. Then, up to a constant, $-\log q(x)$ becomes identical with the code-length $\kappa(x)$, so that the log loss game is essentially equivalent to the code-length game.

By analogy with minimax results of game theory, one might conjecture that

$$(2) \qquad \sup_{P \in \Gamma} \inf_{q \in \mathcal{A}} \mathrm{E}_P\{-\log q(X)\} = \inf_{q \in \mathcal{A}} \sup_{P \in \Gamma} \mathrm{E}_P\{-\log q(X)\}.$$

As we have seen, $P$ achieving the supremum on the left-hand side of (2) is a maximum entropy distribution in $\Gamma$. However, just as important, $q$ achieving the infimum on the right-hand side of (2) is a *robust Bayes* act against $\Gamma$, or a $\Gamma$-*minimax* act [Berger (1985)], for the log loss decision problem.

Now it turns out that, when $\Gamma$ is closed and convex, (2) does indeed hold under very general conditions. Moreover the infimum on the right-hand side is achieved uniquely for $q = p^*$, the probability mass function of the maximum entropy distribution $P^*$. Thus, in this game between DM and Nature, the maximum entropy distribution $P^*$ may be viewed, simultaneously, as defining both Nature's maximin and—in our view more interesting—DM's minimax strategy. In other words, *maximum entropy is robust Bayes*. This decision-theoretic reinterpretation might now be regarded as a plausible justification for selecting the maximum entropy distribution. Note particularly that we do *not* restrict the acts $q$ available to DM to those corresponding to a distribution in the restricted set $\Gamma$: that the optimal act $p^*$ does indeed turn out to have this property is a consequence of, not a restriction on, the analysis.

The maximum entropy method has been most commonly applied in the setting where $\Gamma$ is described by *mean-value constraints* [Jaynes (1989) and Csiszár (1991)]: $\Gamma = \{P : \mathrm{E}_P(T) = \tau\}$, where $T = t(X) \in \mathcal{R}^k$ is some given real- or vector-valued statistic. As pointed out by Grünwald (1998), for such constraints the property (2) is particularly easy to show. By the general theory of exponential families [Barndorff-Nielsen (1978)], under some mild conditions on $\tau$ there will exist a distribution $P^*$ satisfying the constraint $\mathrm{E}_{P^*}(T) = \tau$ and having probability mass function of the form $p^*(x) = \exp\{\alpha_0 + \alpha^\mathrm{T} t(x)\}$ for some $\alpha \in \mathcal{R}^k$, $\alpha_0 \in \mathcal{R}$. Then, for any $P \in \Gamma$,

$$(3) \qquad \mathrm{E}_P\{-\log p^*(X)\} = -\alpha_0 - \alpha^\mathrm{T} \mathrm{E}_P(T) = -\alpha_0 - \alpha^\mathrm{T} \tau = H(P^*).$$

We thus see that $p^*$ is an "equalizer rule" against $\Gamma$, having the same expected loss under any $P \in \Gamma$.

To see that $P^*$ maximizes entropy, observe that, for any $P \in \Gamma$,

$$(4) \qquad H(P) = \inf_{q \in \mathcal{A}} \mathrm{E}_P\{-\log q(X)\} \leq \mathrm{E}_P\{-\log p^*(X)\} = H(P^*),$$

by (3).

To see that $p^*$ is robust Bayes and that (2) holds, note that, for any $q \in \mathcal{A}$,

$$(5) \quad \sup_{P \in \Gamma} \mathrm{E}_P\{-\log q(X)\} \geq \mathrm{E}_{P^*}\{-\log q(X)\} \geq \mathrm{E}_{P^*}\{-\log p^*(X)\} = H(P^*),$$

where the second inequality is the information inequality [Cover and Thomas (1991)]. Hence

$$(6) \qquad H(P^*) \leq \inf_{q \in \mathcal{A}} \sup_{P \in \Gamma} \mathrm{E}_P\{-\log q(X)\}.$$

However, it follows trivially from the "equalizer" property (3) of $p^*$ that

$$(7) \qquad \sup_{P \in \Gamma} \mathrm{E}_P\{-\log p^*(X)\} = H(P^*).$$

From (6) and (7), we see that the choice $q = p^*$ achieves the infimum on the right-hand side of (2) and is thus robust Bayes. Moreover, (2) holds, with both sides equal to $H(P^*)$.

The above argument can be extended to much more general sample spaces (see Section 7). Although this game-theoretic approach and result date back at least to Topsøe (1979), they seem to have attracted little attention so far.

2.2. *This work*: *generalized entropy*. The above robust Bayes view of maximum entropy might be regarded as justifying its use in those decision problems, such as *discrete coding* and *Kelly gambling* [Cover and Thomas (1991)], where the log loss is clearly an appropriate loss function to use. But what if we are interested in other loss functions? This is the principal question we address in this paper.

2.2.1. *Generalized entropy and robust Bayes acts.* We first recall, in Section 3, a natural generalization of the concept of "entropy" (or "uncertainty inherent in a distribution"), related to a specific decision problem and loss function facing DM. The generalized entropy thus associated with the log loss problem is just the Shannon entropy. More generally, let $\mathcal{A}$ be some space of actions or decisions and let $\mathcal{X}$ be the (not necessarily finite) space of possible outcomes to be observed. Let the loss function be given by $L : \mathcal{X} \times \mathcal{A} \to (-\infty, \infty]$, and let $\Gamma$ be a convex set of distributions over $\mathcal{X}$. In Sections 4–6 we set up a statistical game $\mathcal{G}^\Gamma$ based on these ingredients and use this to show that, under a variety of broad regularity conditions, the distribution $P^*$ maximizing, over $\Gamma$, the generalized entropy associated with the loss function $L$ has a Bayes act $a^* \in \mathcal{A}$ [achieving $\inf_{a \in \mathcal{A}} L(P^*, a)$] that is a robust Bayes ($\Gamma$-minimax) decision relative to $L$—thus generalizing the result for the log loss described in Section 2.1. Some variations on this result are also given.

2.2.2. *Generalized exponential families.* In Section 7 we consider in detail the case of *mean-value constraints*, of the form $\Gamma = \{P : \mathrm{E}_P(T) = \tau\}$. For fixed loss function $L$ and statistic $T$, as $\tau$ varies we obtain a family of maximum generalized entropy distributions, one for each value of $\tau$. For Shannon entropy, this turns out to coincide with the *exponential family* having natural sufficient statistic $T$ [Csiszár (1975)]. In close analogy we define the collection of maximum generalized entropy distributions, as we vary $\tau$, to be the *generalized exponential family* determined by $L$ and $T$, and we give several examples of such generalized exponential families. In particular, Lafferty's "additive models based on Bregman divergences" [Lafferty (1999)] are special cases of our generalized exponential families (Section 8.4.2).

2.2.3. *Generalized relative entropy and discrepancy.* In Section 8 we describe how generalized entropy extends to *generalized relative entropy* and show how this in turn is intimately related to a *discrepancy* or *divergence function*. Maximum generalized relative entropy then becomes a special case of the minimum discrepancy method. For the log loss, the associated discrepancy function is just the familiar Kullback–Leibler divergence, and the method then coincides with the "classical" minimum relative entropy method [Jaynes (1989); note that, for Jaynes, "relative entropy" is the same as Kullback–Leibler divergence; for us it is the negative of this].

2.2.4. *A generalized redundancy-capacity theorem.* In many statistical decision problems it is more natural to seek minimax decisions with respect to the discrepancy associated with a loss, rather than with respect to the loss directly. With any game we thus associate a new "derived game," in which the discrepancy constructed from the loss function of the original game now serves as a new loss function. In Section 9 we show that our minimax theorems apply to games of this form too: broadly, whenever the conditions for such a theorem hold for the original game, they also hold for the derived game. As a special case, we reprove a minimax theorem for the Kullback–Leibler divergence [Haussler (1997)], known in information theory as the redundancy-capacity theorem [Merhav and Feder (1995)].

2.2.5. *The Pythagorean property.* The Kullback–Leibler divergence has a celebrated property reminiscent of squared Euclidean distance: it satisfies an analogue of the Pythagorean theorem [Csiszár (1975)]. It has been noted [Csiszár (1991), Jones and Byrne (1990) and Lafferty (1999)] that a version of this property is shared by the broader class of Bregman divergences. In Section 10 we show that a "Pythagorean inequality" in fact holds for the discrepancy based on an arbitrary loss function $L$, so long as the game $\mathcal{G}^\Gamma$ has a value; that is, an analogue of (2) holds. Such decision-based discrepancies include Bregman divergences as special cases. We demonstrate that, even for the case of mean-value constraints, the Pythagorean inequality for a Bregman divergence may be strict.

2.2.6. Finally, Section 11 takes stock of what has been achieved and presents some suggestions for further development.

**3. Decision problems.** In this section we set out some general definitions and properties we shall require. For more background on the concepts discussed here, see Dawid (1998).

A DM has to take some action $a$ selected from a given *action space* $\mathcal{A}$, after which Nature will reveal the value $x \in \mathcal{X}$ of a quantity $X$, and DM will then suffer a loss $L(x,a)$ in $(-\infty, \infty]$. We suppose that Nature takes no account of the action chosen by DM. Then this can be considered as a zero-sum game between Nature and DM, with both players moving simultaneously, and DM paying Nature $L(x,a)$ after both moves are revealed. We call such a combination $\mathcal{G} := (\mathcal{X}, \mathcal{A}, L)$ a *basic game*.

Both DM and Nature are also allowed to make randomized moves, such a move being described by a probability distribution $P$ over $\mathcal{X}$ (for Nature) or $\zeta$ over $\mathcal{A}$ (for DM). We assume that suitable $\sigma$-fields, containing all singleton sets, have been specified in $\mathcal{X}$ and $\mathcal{A}$, and that any probability distributions considered are defined over the relevant $\sigma$-field; we denote the family of all such probability distributions on $\mathcal{X}$ by $\mathcal{P}_0$. We further suppose that the loss function $L$ is jointly measurable.

3.1. *Expected loss.* We shall permit algebraic operations on the extended real line $[-\infty, \infty]$, with definitions and exceptions as in Rockafellar (1970), Section 4.

For a function $f:\mathcal{X} \to [-\infty, \infty]$, and $P \in \mathcal{P}_0$, we may denote $\mathrm{E}_P\{f(X)\}$ [i.e., $\mathrm{E}_{X \sim P}\{f(X)\}$] by $f(P)$. When $f$ is bounded below, $f(P)$ is construed as $\infty$ if $P\{f(X) = \infty\} > 0$. When $f$ is unbounded, we interpret $f(P)$ as $f^+(P) - f^-(P) \in [-\infty, +\infty]$, where $f^+(x) := \max\{f(x), 0\}$ and $f^-(x) := \max\{-f(x), 0\}$, allowing either $f^+(P)$ or $f^-(P)$ to take the value $\infty$, but not both. In this last case $f(P)$ is undefined, else it is *defined* (either as a finite number or as $\pm\infty$).

If DM knows that Nature is generating $X$ from $P$ or, in the absence of such knowledge, DM is using $P$ to represent his or her own uncertainty about $X$, then the undesirability to DM of any act $a \in \mathcal{A}$ will be assessed by means of its *expected loss*,

$$(8) \qquad L(P, a) := \mathrm{E}_P\{L(X, a)\}.$$

We can similarly extend $L$ to randomized acts: $L(x, \zeta) := \mathrm{E}_{A \sim \zeta}\{L(x, A)\}$, $L(P, \zeta) = \mathrm{E}_{(X,A) \sim P \times \zeta}\{L(X, A)\}$.

Throughout this paper we shall mostly confine attention to probability measures $P \in \mathcal{P}_0$ such that $L(P, a)$ is defined for all $a \in \mathcal{A}$, and we shall denote the family of all such $P$ by $\mathcal{P}$. We further confine attention to randomized acts $\zeta$ such that $L(P, \zeta)$ is defined for all $P \in \mathcal{P}$, denoting the set of all such $\zeta$ by $\mathcal{Z}$. Note that any distribution degenerate at a point $x \in \mathcal{X}$ is in $\mathcal{P}$, and so $L(x, \zeta)$ is defined for all $x \in \mathcal{X}$, $\zeta \in \mathcal{Z}$.

LEMMA 3.1.  *For all $P \in \mathcal{P}$, $\zeta \in \mathcal{Z}$,*

$$(9) \qquad L(P, \zeta) = \mathrm{E}_{X \sim P}\{L(X, \zeta)\} = \mathrm{E}_{A \sim \zeta}\{L(P, A)\}.$$

PROOF. When $L(P, \zeta)$ is finite this is just Fubini's theorem.

Now consider the case $L(P, \zeta) = \infty$. First suppose $L \geq 0$ everywhere. If $L(x, \zeta) = \infty$ for $x$ in a subset of $\mathcal{X}$ having positive $P$-measure, then (9) holds, both sides being $+\infty$. Otherwise, $L(x, \zeta)$ is finite almost surely $[P]$. If $\mathrm{E}_P\{L(X, \zeta)\}$ were finite, then by Fubini it would be the same as $L(P, \zeta)$. So once again $\mathrm{E}_P\{L(X, \zeta)\} = L(P, \zeta) = +\infty$.

This result now extends easily to possibly negative $L$, on noting that $L^-(P, \zeta)$ must be finite; a parallel result holds when $L(P, \zeta) = -\infty$.

Finally the whole argument can be repeated after interchanging the roles of $x$ and $a$ and of $P$ and $\zeta$.  □

COROLLARY 3.1.  *For any $P \in \mathcal{P}$,*

$$(10) \qquad \inf_{\zeta \in \mathcal{Z}} L(P, \zeta) = \inf_{a \in \mathcal{A}} L(P, a).$$

PROOF. Clearly $\inf_{\zeta \in \mathcal{Z}} L(P, \zeta) \leq \inf_{a \in \mathcal{A}} L(P, a)$. If $\inf_{a \in \mathcal{A}} L(P, a) = -\infty$ we are done. Otherwise, for any $\zeta \in \mathcal{Z}$, $L(P, \zeta) = \mathrm{E}_{A \sim \zeta} L(P, A) \geq \inf_{a \in \mathcal{A}} L(P, a)$.  □

We shall need the fact that, for any $\zeta \in \mathcal{Z}$, $L(P,\zeta)$ is linear in $P$ in the following sense.

LEMMA 3.2. *Let $P_0, P_1 \in \mathcal{P}$, and let $P_\lambda := (1-\lambda)P_0 + \lambda P_1$. Fix $\zeta \in \mathcal{Z}$, such that the pair $\{L(P_0,\zeta), L(P_1,\zeta)\}$ does not contain both the values $-\infty$ and $+\infty$. Then, for any $\lambda \in (0,1)$, $L(P_\lambda,\zeta)$ is finite if and only if both $L(P_1,\zeta)$ and $L(P_0,\zeta)$ are. In this case $L(P_\lambda,\zeta) = (1-\lambda)\,L(P_0,\zeta) + \lambda\,L(P_1,\zeta)$.*

PROOF. Consider a bivariate random variable $(I, X)$ with joint distribution $P^*$ over $\{0,1\} \times \mathcal{X}$ specified by the following: $I = 1, 0$ with respective probabilities $\lambda$, $1-\lambda$; and, given $I = i$, $X$ has distribution $P_i$. By Fubini we have

$$\mathrm{E}_{P^*}\{L(X,\zeta)\} = \mathrm{E}_{P^*}[\mathrm{E}_{P^*}\{L(X,\zeta)|I\}],$$

in the sense that, whenever one side of this equation is defined and finite, the same holds for the other, and they are equal. Noting that, under $P^*$, the distribution of $X$ is $P_\lambda$ marginally, and $P_i$ conditional on $I = i$ ($i = 0, 1$), the result follows. □

3.2. *Bayes act.* Intuitively, when $X \sim P$ an act $a_P \in \mathcal{A}$ will be optimal if it minimizes $L(P, a)$ over all $a \in \mathcal{A}$. Any such act $a_P$ is a *Bayes act* against $P$. More generally, to allow for the possibility that $L(P, a)$ may be infinite as well as to take into account randomization, we call $\zeta_P \in \mathcal{Z}$ a (*randomized*) *Bayes act*, or simply *Bayes*, against $P$ (not necessarily in $\mathcal{P}$) if

(11) $$\mathrm{E}_P\{L(X,\zeta) - L(X,\zeta_P)\} \in [0,\infty]$$

for all $\zeta \in \mathcal{Z}$. We denote by $\mathcal{A}_P$ (resp. $\mathcal{Z}_P$) the set of all nonrandomized (resp. randomized) Bayes acts against $P$. Clearly $\mathcal{A}_P \subseteq \mathcal{Z}_P$, and $L(P, \zeta_P)$ is the same for all $\zeta_P \in \mathcal{Z}_P$.

The loss function $L$ will be called $\Gamma$-*strict* if, for each $P \in \Gamma$, there exists $a_P \in \mathcal{A}$ that is the unique Bayes act against $P$; $L$ is $\Gamma$-*semistrict* if, for each $P \in \Gamma$, $\mathcal{A}_P$ is nonempty, and $a, a' \in \mathcal{A}_P \Rightarrow L(\cdot, a) \equiv L(\cdot, a')$. When $L$ is $\Gamma$-strict, and $P \in \Gamma$, it can never be optimal for DM to choose a randomized act; when $L$ is $\Gamma$-semistrict, even though a randomized act can be optimal there is never any point in choosing one, since its loss function will be identical with that of any nonrandomized optimal act.

Semistrictness is clearly weaker than strictness. For our purposes we can replace it by the still weaker concept of *relative strictness*: $L$ is $\Gamma$-*relatively strict* if for all $P \in \Gamma$ the set of Bayes acts $\mathcal{A}_P$ is nonempty and, for all $a, a' \in \mathcal{A}_P$, $L(P', a) = L(P', a')$ for all $P' \in \Gamma$.

3.3. *Bayes loss and entropy.* Whether or not a Bayes act exists, the *Bayes loss* $H(P) \in [-\infty, \infty]$ of a distribution $P \in \mathcal{P}$ is defined by

$$H(P) := \inf_{a \in \mathcal{A}} L(P, a). \tag{12}$$

It follows from Corollary 3.1 that it would make no difference if the infimum in (12) were extended to be over $\zeta \in \mathcal{Z}$. We shall mostly be interested in Bayes acts of distributions $P$ with finite $H(P)$. In the context of Section 2.1, with $L(x, q)$ the log loss $-\log q(x)$, $H(P)$ is just the Shannon entropy of $P$.

PROPOSITION 3.1. *Let $P \in \mathcal{P}$ and suppose $H(P)$ is finite. Then the following hold:*

(i) *$\zeta_P \in \mathcal{Z}$ is Bayes against $P$ if and only if*

$$\mathrm{E}_P\{L(X, a) - L(X, \zeta_P)\} \in [0, \infty] \tag{13}$$

*for all $a \in \mathcal{A}$.*

(ii) *$\zeta_P$ is Bayes against $P$ if and only if $L(P, \zeta_P) = H(P)$.*

(iii) *If $P$ admits some randomized Bayes act, then $P$ also admits some nonrandomized Bayes act; that is, $\mathcal{A}_P$ is not empty.*

PROOF. Items (i) and (ii) follow easily from (10) and finiteness. To prove (iii), let $f(P, a) := L(P, a) - H(P)$. Then $f(P, a) \geq 0$ for all $a$, while $\mathrm{E}_{A \sim \zeta_P} f(P, A) = L(P, \zeta_P) - H(P) = 0$. We deduce that $\{a \in \mathcal{A}: f(P, a) = 0\}$ has probability 1 under $\zeta_P$ and so, in particular, must be nonempty. □

We express the well-known concavity property of the Bayes loss [DeGroot (1970), Section 8.4] as follows.

PROPOSITION 3.2. *Let $P_0, P_1 \in \mathcal{P}$, and let $P_\lambda := (1-\lambda)P_0 + \lambda P_1$. Suppose that $H(P_i) < \infty$ for $i = 0, 1$. Then $H(P_\lambda)$ is a concave function of $\lambda$ on $[0, 1]$ (and thus, in particular, continuous on $(0, 1)$ and lower semicontinuous on $[0, 1]$). It is either bounded above on $[0, 1]$ or infinite everywhere on $(0, 1)$.*

PROOF. Let $\mathcal{B}$ be the set of all $a \in \mathcal{A}$ such that $L(P_\lambda, a) < \infty$ for some $\lambda \in (0, 1)$—and thus, by Lemma 3.2, for all $\lambda \in [0, 1]$. If $\mathcal{B}$ is empty, then $H(P_\lambda) = \infty$ for all $\lambda \in (0, 1)$; in particular, $H(P_\lambda)$ is then concave on $[0, 1]$. Otherwise, taking any fixed $a \in \mathcal{B}$ we have $H(P_\lambda) \leq L(P_\lambda, a) \leq \max_i L(P_i, a)$, so $H(P_\lambda)$ is bounded above on $[0, 1]$. Moreover, as the pointwise infimum of the nonempty family of concave functions $\{L(P_\lambda, a): a \in \mathcal{A}\}$, $H(P_\lambda)$ is itself a concave function of $\lambda$ on $[0, 1]$. □

COROLLARY 3.2. *If for all $a \in \mathcal{A}$, $L(P_\lambda, a) < \infty$ for some $\lambda \in (0, 1)$, then for all $\lambda \in [0, 1]$, $H(P_\lambda) = \lim\{H(P_\mu): \mu \in [0, 1], \mu \to \lambda\}$ [it being allowed that $H(P_\lambda)$ is not finite].*

PROOF. In this case $\mathcal{B} = \mathcal{A}$, so that $H(P_\lambda) = \inf_{a \in \mathcal{B}} L(P_\lambda, a)$. Each function $L(P_\lambda, a)$ is finite and linear, hence a closed concave function of $\lambda$ on $[0, 1]$. This last property is then preserved on taking the infimum. The result now follows from Theorem 7.5 of Rockafellar (1970). □

COROLLARY 3.3. *If in addition $H(P_i)$ is finite for $i = 0, 1$, then $H(P_\lambda)$ is a bounded continuous function of $\lambda$ on $[0, 1]$.*

Note that Corollary 3.3 will always apply when the loss function is bounded.

Under some further regularity conditions [see Dawid (1998, 2003) and Section 3.5.4 below], a general concave function over $\mathcal{P}$ can be regarded as generated from some decision problem by means of (12). Concave functions have been previously proposed as general measures of the uncertainty or diversity in a distribution [DeGroot (1962) and Rao (1982)], generalizing the Shannon entropy. We shall thus call the Bayes loss $H$, as given by (12), the (*generalized*) *entropy function* or *uncertainty function* associated with the loss function $L$.

3.4. *Scoring rule.* Suppose the action space $\mathcal{A}$ is itself a set $\mathcal{Q}$ of distributions for $X$. Note we are not here considering $Q \in \mathcal{Q}$ as a randomized act over $\mathcal{X}$, but rather as a simple act in its own right (e.g., a decision to quote $Q$ as a description of uncertainty about $X$). We typically write the loss as $S(x, Q)$ in this case and refer to $S$ as a *scoring rule* or *score*. Such scoring rules are used to assess the performance of probability forecasters [Dawid (1986)]. We say $S$ is $\Gamma$-*proper* if $\Gamma \subseteq \mathcal{Q} \subseteq \mathcal{P}$ and, for all $P \in \Gamma$, the choice $Q = P$ is Bayes against $X \sim P$. Then for $P \in \Gamma$,

$$H(P) = S(P, P). \tag{14}$$

Suppose now we start from a general decision problem, with loss function $L$ such that $\mathcal{Z}_Q$ is nonempty for all $Q \in \mathcal{Q}$. Then we can define a scoring rule by

$$S(x, Q) := L(x, \zeta_Q), \tag{15}$$

where for each $Q \in \mathcal{Q}$ we suppose we have selected some specific Bayes act $\zeta_Q \in \mathcal{Z}_Q$. Then for $P \in \mathcal{Q}$, $S(P, Q) = L(P, \zeta_Q)$ is clearly minimized when $Q = P$, so that this scoring rule is $\mathcal{Q}$-proper. If $L$ is $\mathcal{Q}$-semistrict, then (15) does not depend on the choice of Bayes act $\zeta_Q$. More generally, if $L$ is $\mathcal{Q}$-relatively strict, then $S(P, Q)$ does not depend on such a choice, for all $P, Q \in \mathcal{Q}$.

We see that, for $P \in \mathcal{Q}$, $\inf_{Q \in \mathcal{Q}} S(P, Q) = S(P, P) = L(P, \zeta_P) = H(P)$. In particular, the generalized entropy associated with the constructed scoring rule (15) is identical with that determined by the original loss function $L$. In this way, almost any decision problem can be reformulated in terms of a proper scoring rule.

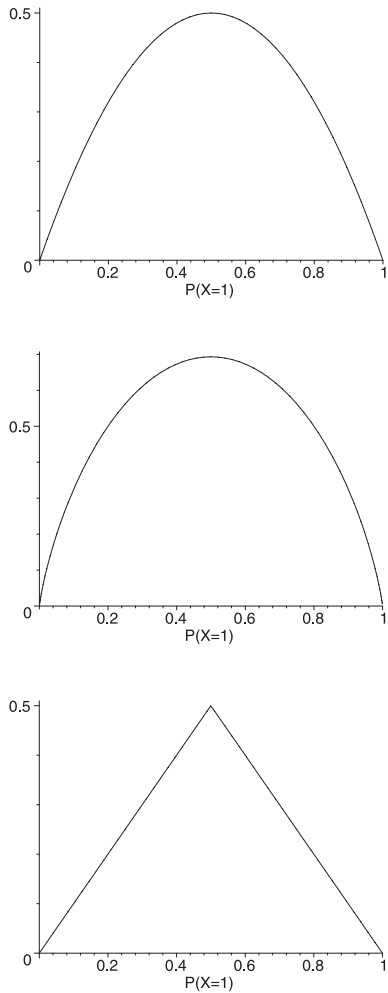

Fig. 1. *Brier, log and zero–one entropies for the case $\mathcal{X} = \{0,1\}$.*

3.5. *Some examples.* We now give some simple examples, both to illustrate the above concepts and to provide a concrete focus for later development. Further examples may be found in Dawid (1998) and Dawid and Sebastiani (1999).

3.5.1. *Brier score.* Although it can be generalized, we restrict our treatment of the *Brier score* [Brier (1950)] to the case of a finite sample space $\mathcal{X} = \{x_1, \ldots, x_N\}$. A distribution $P$ over $\mathcal{X}$ can be represented by its probability vector $p = (p(1), \ldots, p(N))$, where $p(x) := P(X = x)$. A point $x \in \mathcal{X}$ may also be represented by the $N$-vector $\delta^x$ corresponding to the point-mass distribution on $\{x\}$ having entries $\delta^x(j) = 1$ if $j = x$, 0 otherwise. The Brier scoring rule is then defined by

$$S(x, Q) := \|\delta^x - q\|^2 \tag{16}$$

$$= \sum_{j=1}^{N} \{\delta^x(j) - q(j)\}^2$$

$$= \sum_{j} q(j)^2 - 2q(x) + 1. \tag{17}$$

Then

$$S(P, Q) = \sum_{j} q(j)^2 - 2\sum_{j} p(j)\, q(j) + 1, \tag{18}$$

which is uniquely minimized for $Q = P$, so that this is a $\mathcal{P}$-strict proper scoring rule. The corresponding entropy function is (see Figure 1)

$$H(P) = 1 - \sum_{j} p(j)^2. \tag{19}$$

3.5.2. *Logarithmic score.* An important scoring rule is the *logarithmic score*, generalizing the discrete-case log loss as already considered in Section 2. For a general sample space $\mathcal{X}$, let $\mu$ be a fixed $\sigma$-finite measure (the *base measure*) on a suitable $\sigma$-algebra in $\mathcal{X}$, and take $\mathcal{A}$ to be the set of all finite nonnegative measurable real functions $q$ on $\mathcal{X}$ for which $\int q(x)\, d\mu(x) = 1$. Any $q \in \mathcal{A}$ can be regarded as the density of a distribution $Q$ over $\mathcal{X}$ which is absolutely continuous with respect to $\mu$. We denote the set of such distributions by $\mathcal{M}$. However, because densities are only defined up to a set of measure 0, different $q$'s in $\mathcal{A}$ can correspond to the same $Q \in \mathcal{M}$. Note moreover that the many–one correspondence between $q$ and $Q$ depends on the specific choice of base measure $\mu$ and will change if we change $\mu$.

We define a loss function by

$$S(x, q) = -\log q(x). \tag{20}$$

If (and only if) $P \in \mathcal{M}$, then $S(P, q)$ will be the same for all versions $q$ of the density of the same distribution $Q \in \mathcal{M}$. Hence for $P, Q \in \mathcal{M}$ we can write $S(P, Q)$ instead of $S(P, q)$, and we can consider $S$ to be a scoring rule. It is well known that, for $P, Q, Q^* \in \mathcal{M}$, $\mathrm{E}_P\{S(X, Q) - S(X, Q^*)\} = -\int p(x) \log\{q(x)/q^*(x)\}\, d\mu$ is nonnegative for all $Q$ if and only if $Q^* = P$. That is, $Q^*$ is Bayes against $P$ if and only if $Q^* = P$, so that this scoring rule is $\mathcal{M}$-strictly proper.

We have, for $P \in \mathcal{M}$,

$$H(P) = -\int p(x) \log p(x)\, d\mu, \tag{21}$$

the usual definition of the *entropy* of $P$ with respect to $\mu$. When $\mathcal{X}$ is discrete and $\mu$ is counting measure, we recover the Shannon entropy. For the simple case $\mathcal{X} = \{0,1\}$ this is depicted in Figure 1. Note that the whole decision problem, and in particular the value of $H(P)$ as given by (21), will be altered if we change (even in a mutually absolutely continuous way) the base measure $\mu$.

Things simplify when $\mu$ is itself a probability measure. In this case $\mathcal{A}$ contains the constant function 1. For any distribution $P$ whatsoever, whether or not $P \in \mathcal{M}$, we have $L(P, 1) = 0$, whence we deduce $H(P) \leq 0$ (with equality if and only if $P = \mu$). When $P \in \mathcal{M}$, (21) asserts $H(P) = -\mathrm{KL}(P, \mu)$, where KL is the Kullback–Leibler divergence [Kullback (1959)]. [Note that it is possible to have $\mathrm{KL}(P, \mu) = \infty$, and thus $H(P) = -\infty$, even for $P \in \mathcal{M}$.] If $P \notin \mathcal{M}$, there exist a measurable set $N$ and $\alpha > 0$ such that $\mu(N) = 0$ but $P(N) = \alpha$. Define $q_n(x) = 1$ $(x \notin N)$, $q_n(x) = n$ $(x \in N)$. Then $q_n \in \mathcal{A}$ and $L(P, q_n) = -\alpha \log n$. It follows that $H(P) = -\infty$. Since the usual definition [Csiszár (1975) and Posner (1975)] has $\mathrm{KL}(P, \mu) = \infty$ when $P \not\ll \mu$, we thus have $H(P) = -\mathrm{KL}(P, \mu)$ in all cases. This formula exhibits clearly the dependence of the entropy on the choice of $\mu$.

3.5.3. *Zero–one loss.* Let $\mathcal{X}$ be finite or countable, take $\mathcal{A} = \mathcal{X}$ and consider the loss function

$$L(x, a) = \begin{cases} 0, & \text{if } a = x, \\ 1, & \text{otherwise.} \end{cases} \tag{22}$$

Then $L(P, a) = 1 - P(X = a)$, and a nonrandomized Bayes act under $P$ is any mode of $P$. When $P$ has (at least) two modes, say $a_P$ and $a'_P$, then $L(x, a_P)$ and $L(x, a'_P)$ are not identical, so that this loss function is not $\mathcal{P}$-semistrict. This means that we may have to take account of randomized strategies $\zeta$ for DM. Then, writing $\zeta(x) := \zeta(A = x)$, we have

$$L(x, \zeta) = 1 - \zeta(x) \tag{23}$$

and

$$L(P, \zeta) = 1 - \sum_{x \in \mathcal{X}} p(x)\zeta(x). \tag{24}$$

A randomized act $\zeta$ is Bayes against $P$ if and only if it puts all its mass on the set of modes of $P$.

We have generalized entropy function

$$H(P) = 1 - p_{\max}, \tag{25}$$

with $p_{\max} := \sup_{x \in \mathcal{X}} p(x)$. For the simple case $\mathcal{X} = \{0,1\}$, this is depicted in Figure 1.

3.5.4. *Bregman score.* Suppose that $\#(\mathcal{X}) = N < \infty$ and that we represent a distribution $P \in \mathcal{P}$ over $\mathcal{X}$ by its probability mass function $p \in \Delta$, the unit simplex in $\mathcal{R}^N$, which can in turn be considered as a subset of $(N-1)$-dimensional Euclidean space. The interior $\Delta^\circ$ of $\Delta$ then corresponds to the subset $\mathcal{Q} \subset \mathcal{P}$ of distributions giving positive probability to each point of $\mathcal{X}$.

Let $H$ be a finite concave real function on $\Delta$. For any $q \in \Delta^\circ$, the set $\nabla H(q)$ of supporting hyperplanes to $H$ at $q$ is nonempty [Rockafellar (1970), Theorem 27.3]—having a unique member when $H$ is differentiable at $q$. Select for each $q \in \Delta^\circ$ some specific member of $\nabla H(q)$, and let the height of this hyperplane at arbitrary $p \in \Delta$ be denoted by $l_q(p)$: this affine function must then have equation of the form

$$l_q(p) = H(q) + \alpha_q^{\mathrm{T}}(p - q). \tag{26}$$

Although the coefficient vector $\alpha_q \in \mathcal{R}^{\mathcal{X}}$ in (26) is only defined up to addition of a multiple of the unit vector, this arbitrariness will be of no consequence. We shall henceforth reuse the notation $\nabla H(q)$ in place of $\alpha_q$.

By the supporting hyperplane property,

$$l_q(p) \geq H(p), \tag{27}$$

$$l_q(q) = H(q). \tag{28}$$

Now consider the function $S: \mathcal{X} \times \mathcal{Q}$ defined by

$$S(x, Q) = H(q) + \nabla H(q)^{\mathrm{T}}(\delta^x - q), \tag{29}$$

where $\delta^x$ is the vector having $\delta^x(j) = 1$ if $j = x$, 0 otherwise.

Then we easily see that $S(P, Q) = l_q(p)$, so that, by (27) and (28), $S(P, Q)$ is minimized in $Q$ when $Q = P$. Thus $S$ is a $\mathcal{Q}$-proper scoring rule.

We note that

$$\begin{aligned}0 \leq d(P, Q) &:= S(P, Q) - S(P, P) \\ &= H(q) + \nabla H(q)^{\mathrm{T}}(p - q) - H(p).\end{aligned} \tag{30}$$

With further regularity conditions (including in particular differentiability), (30) becomes the *Bregman divergence* [Brègman (1967), Csiszár (1991) and Censor and Zenios (1997)] associated with the convex function $-H$. We therefore call $S$, defined as in (29), a *Bregman score* associated with $H$. This will be unique when $H$ is differentiable on $\Delta^\circ$. In Section 8 we introduce a more general decision-theoretic notion of divergence.

We note by (28) that the generalized entropy function associated with this score is $H^*(P) = S(P, P) = l_p(p) = H(p)$ (at any rate inside $\Delta^\circ$). That is to say, we have exhibited a decision problem for which a prespecified concave function $H$ is the entropy. This construction can be extended to the whole of $\Delta$ and to certain concave functions $H$ that are not necessarily finite [Dawid (2003)]. Extensions can also be made to more general sample spaces.

3.5.5. *Separable Bregman score.* A special case of the construction of Section 3.5.4 arises when we take $H(q)$ to have the form $-\sum_{x \in \mathcal{X}} \psi\{q(x)\}$, with $\psi$ a real-valued differentiable convex function of a nonnegative argument. In this case we can take $(\nabla H(q))(x) = -\psi'\{q(x)\}$, and the associated proper scoring rule has

$$(31) \qquad S(x, Q) = -\psi'\{q(x)\} - \sum_{t \in \mathcal{X}} [\psi\{q(t)\} - q(t)\psi'\{q(t)\}].$$

We term this the *separable Bregman scoring rule* associated with $\psi$. The corresponding *separable Bregman divergence* [confusingly, this special case of (30) is sometimes also referred to simply as a Bregman divergence] is

$$(32) \qquad d_\psi(P, Q) = \sum_{x \in \mathcal{X}} \Delta_\psi\{p(x), q(x)\},$$

where we have introduced

$$(33) \qquad \Delta_\psi(a, b) := \psi(a) - \psi(b) - \psi'(b)(a - b).$$

The nonnegative function $\Delta_\psi$ measures how much the convex function $\psi$ deviates at $a$ from its tangent at $b$; this can be considered as a measure of "how convex" $\psi$ is.

We can easily extend the above definition to more general sample spaces. Thus let $\mathcal{X}$, $\mu$, $\mathcal{A}$ and $\mathcal{M}$ be as in Section 3.5.2, and, in analogy with (31), consider the following loss function:

$$(34) \qquad S(x, q) := -\psi'\{q(x)\} - \int [\psi\{q(t)\} - q(t)\psi'\{q(t)\}] \, d\mu(t).$$

Clearly if $q$, $q'$ are both $\mu$-densities of the same $Q \in \mathcal{M}$, then $S(x, q) = S(x, q')$ a.e. $[\mu]$, and so, for any $P \in \mathcal{M}$, $S(P, q) = S(P, q')$. Thus once again, for $P, Q \in \mathcal{M}$, we can simply write $S(P, Q)$. We then have

$$(35) \qquad S(P, Q) = \int [\{q(t) - p(t)\} \psi'\{q(t)\} - \psi\{q(t)\}] \, d\mu(t),$$

whence

$$(36) \qquad S(P, P) = -\int \psi\{p(t)\} \, d\mu(t),$$

and so, if $S(P, P)$ is finite,

$$(37) \qquad d_\psi(P, Q) := S(P, Q) - S(P, P) = \int \Delta_\psi\{p(t), q(t)\} \, d\mu(t).$$

Thus, for $P, Q \in \mathcal{M}$, if $S(P, P)$ is finite, $S(P, P) \leq S(P, Q)$. Using the extended definition (11) of Bayes acts, we can show that $P$ is Bayes against $P$ even when $S(P, P)$ is infinite. That is, $S$ is an $\mathcal{M}$-proper scoring rule. If $\psi$ is strictly convex, $S$ is $\mathcal{M}$-strict.

The quantity $d_\psi(P, Q)$ defined by (37) is identical with the (*separable*) *Bregman divergence* [Brègman (1967) and Csiszár (1991)] $B_\psi(p, q)$, based on $\psi$ (and $\mu$), between the densities $p$ and $q$ of $P$ and $Q$. Consequently, we shall term $S(x, q)$ given by (34) a *separable Bregman score*. For $P \in \mathcal{M}$ the associated *separable Bregman entropy* is then, by (36),

$$H_\psi(P) = -\int \psi\{p(t)\}\,d\mu(t). \tag{38}$$

The logarithmic score arises as a special case of the separable Bregman score on taking $\psi(s) \equiv s \log s$; and the Brier score arises on taking $\mu$ to be counting measure and $\psi(s) \equiv s^2 - 1/N$.

3.5.6. *More examples.* Since every decision problem generates a generalized entropy function, an enormous range of such functions can be constructed. As a very simple case, consider the *quadratic loss* problem, with $\mathcal{X} = \mathcal{A} = \mathcal{R}$, $L(x, a) = (x - a)^2$. Then $a_P = \mathrm{E}_P(X)$ is Bayes against $P$, and the associated proper scoring rule and entropy are $S(x, P) = \{x - \mathrm{E}_P(X)\}^2$ and $H(P) = \mathrm{var}_P(X)$ — a very natural measure of uncertainty. This cannot be expressed in the form (38), so it is not associated with a separable Bregman divergence. Dawid and Sebastiani (1999) characterize all those generalized entropy functions that depend only on the variance of a (possibly multivariate) distribution.

**4. Maximum entropy and robust Bayes.** Suppose that Nature may be regarded as generating $X$ from a distribution $P$, but DM does not know $P$. All that is known is that $P \in \Gamma$, a specified family of distributions over $\mathcal{X}$. The consequence DM faces if he or she takes act $a \in \mathcal{A}$ when Nature chooses $X = x$ is measured by the loss $L(x, a)$. How should DM act?

4.1. *Maximum entropy.* One way of proceeding is to replace the family $\Gamma$ by some "representative" member $P^* \in \Gamma$, and then choose an act that is Bayes against $P^*$. A possible criterion for choosing $P^*$, generalizing the standard maximum Shannon entropy procedure, might be:

*Maximize, over $P \in \Gamma$, the generalized entropy $H(P)$.*

4.2. *Robust Bayes rules.* Another approach is to conduct a form of "robust Bayes analysis" [Berger (1985)]. In particular we investigate the $\Gamma$-*minimax criterion*, a compromise between Bayesian and frequentist decision theory. For a recent tutorial overview of this criterion, see Vidakovic (2000).

When $X \sim P \in \Gamma$, the loss of an act $a$ is evaluated by $L(P, a)$. We can form a new *restricted game*, $\mathcal{G}^\Gamma = (\Gamma, \mathcal{A}, L)$, where Nature selects a distribution $P$ from $\Gamma$, DM an act $a$ from $\mathcal{A}$, and the ensuing loss to DM is taken to be $L(P, a)$. Again, we allow DM to take randomized acts $\zeta \in \mathcal{Z}$, yielding loss $L(P, \zeta)$ when Nature generates $X$ from $P$. In principle we could also let Nature choose her distribution $P$ in some random fashion, described by means of a law (distribution) for a random distribution $\widetilde{P}$ over $\mathcal{X}$. However, with the exception of Section 10, where randomization is in any case excluded, in all the cases we shall consider $\Gamma$ will be convex, and then every randomized act for Nature can be replaced by a nonrandomized act (the mean of the law of $\widetilde{P}$) having the identical loss function. Consequently we shall not consider randomized acts for Nature.

In the absence of knowledge of Nature's choice of $P$, we might apply the minimax criterion to this restricted game. This leads to the prescription for DM:

*Choose* $\zeta = \zeta^* \in \mathcal{Z}$, *to achieve*

$$\inf_{\zeta \in \mathcal{Z}} \sup_{P \in \Gamma} L(P, \zeta). \tag{39}$$

We shall term any act $\zeta^*$ achieving (39) *robust Bayes* against $\Gamma$, or $\Gamma$-*minimax*.

When the basic game is defined in terms of a $\mathcal{Q}$-proper scoring rule $S(x, Q)$, and $\Gamma \subseteq \mathcal{Q}$, this robust Bayes criterion becomes:

*Choose* $Q = Q^*$, *to achieve*

$$\inf_{Q \in \mathcal{Q}} \sup_{P \in \Gamma} S(P, Q). \tag{40}$$

Note particularly that in this case there is no reason to require $\mathcal{Q} = \Gamma$; we might want to take $\mathcal{Q}$ larger than $\Gamma$ (typically, $\mathcal{Q} = \mathcal{P}$). Also, we have not considered randomized acts in (40)—we shall see later that, for the problems we consider, this has no effect.

Below we explore the relationship between the above two methods. In particular, we shall show that, in very general circumstances, they produce identical results. That is, maximum generalized entropy is robust Bayes. This will be the cornerstone of all our results to come.

First note that from (12) the maximum entropy criterion can be expressed as:

*Choose* $P = P^*$, *to achieve*

$$\sup_{P \in \Gamma} \inf_{\zeta \in \mathcal{Z}} L(P, \zeta). \tag{41}$$

There is a striking duality with the criterion (39).

In the general terminology of game theory, (41) defines the extended real *lower value*,

$$\underline{V} := \sup_{P \in \Gamma} \inf_{\zeta \in \mathcal{Z}} L(P, \zeta), \tag{42}$$

and (39) the *upper value*,

(43)
$$\overline{V} := \inf_{\zeta \in \mathcal{Z}} \sup_{P \in \Gamma} L(P, \zeta),$$

of the restricted game $\mathcal{G}^\Gamma$. In particular, the maximum achievable entropy is exactly the lower value. We always have $\underline{V} \leq \overline{V}$. When these two are equal and finite, we say the game $\mathcal{G}^\Gamma$ has a *value*, $V := \underline{V} = \overline{V}$.

DEFINITION 4.1. The pair $(P^*, \zeta^*) \in \Gamma \times \mathcal{Z}$ is a *saddle-point* (or *equilibrium*) in the game $\mathcal{G}^\Gamma$ if $H^* := L(P^*, \zeta^*)$ is finite, and the following hold:

(44)
    (a) $L(P^*, \zeta^*) \leq L(P^*, \zeta)$     for all $\zeta \in \mathcal{Z}$;
    (b) $L(P^*, \zeta^*) \geq L(P, \zeta^*)$     for all $P \in \Gamma$.

In Sections 5 and 6 we show for convex $\Gamma$ the existence of a saddle-point in $\mathcal{G}^\Gamma$ under a variety of broadly applicable conditions.

In certain important special cases [see, e.g., Section 2.1, (3)], we may be able to demonstrate (b) above by showing that $\zeta^*$ is an equalizer rule:

DEFINITION 4.2. $\zeta \in \mathcal{Z}$ is an *equalizer rule* in $\mathcal{G}^\Gamma$ if $L(P, \zeta)$ is the same finite constant for all $P \in \Gamma$.

LEMMA 4.1. *Suppose that there exist both a maximum entropy distribution $P^* \in \Gamma$ achieving (42), and a robust Bayes act $\zeta^* \in \mathcal{Z}$ achieving (43). Then $\underline{V} \leq L(P^*, \zeta^*) \leq \overline{V}$. If, further, the game has a value, $V$ say, then $V = H^* := L(P^*, \zeta^*)$, and $(P^*, \zeta^*)$ is a saddle-point in the game $\mathcal{G}^\Gamma$.*

PROOF. $\underline{V} = \inf_\zeta L(P^*, \zeta) \leq L(P^*, \zeta^*)$, and similarly $L(P^*, \zeta^*) \leq \overline{V}$. If the game has a value $V$, then $L(P^*, \zeta^*) = V = \inf_{\zeta \in \mathcal{Z}} L(P^*, \zeta)$, and $L(P^*, \zeta^*) = V = \sup_{P \in \Gamma} L(P, \zeta^*)$. □

Note that, even when the game has a value, either or both of $P^*$ and $\zeta^*$ may fail to exist.

Conversely, we have the following theorem.

THEOREM 4.1. *Suppose that a saddle-point $(P^*, \zeta^*)$ exists in the game $\mathcal{G}^\Gamma$. Then*:

(i) *The game has value $H^* = L(P^*, \zeta^*)$.*
(ii) *$\zeta^*$ is Bayes against $P^*$.*
(iii) *$H(P^*) = H^*$.*
(iv) *$P^*$ maximizes the entropy $H(P)$ over $\Gamma$.*
(v) *$\zeta^*$ is robust Bayes against $\Gamma$.*

PROOF. Part (i) follows directly from (44) and the definitions of $\underline{V}$, $\overline{V}$. Part (ii) is immediate from (44)(a) and finiteness, and in turn implies (iii). For any $P \in \Gamma$, $H(P) \leq L(P, \zeta^*) \leq H^*$ by (44)(b). Then (iv) follows from (iii). For any $\zeta \in \mathcal{Z}$, $\sup_P L(P, \zeta) \geq L(P^*, \zeta)$, so that, by (44)(a),

$$\sup_P L(P, \zeta) \geq H^*. \tag{45}$$

Also, by (44)(b),

$$\sup_P L(P, \zeta^*) = H^*. \tag{46}$$

Comparing (45) and (46), we see that $\zeta^*$ achieves (39); that is, (v) holds. □

COROLLARY 4.1. *Suppose that $L$ is $\Gamma$-relatively strict, that there is a unique $P^* \in \Gamma$ maximizing the generalized entropy $H$ and that $\zeta^* \in \mathcal{Z}$ is a Bayes act against $P^*$. Then, if $\mathcal{G}^\Gamma$ has a saddle-point, $\zeta^*$ is robust Bayes against $\Gamma$.*

COROLLARY 4.2. *Let the basic game $\mathcal{G}$ be defined in terms of a $\mathcal{Q}$-strictly proper scoring rule $S(x, Q)$, and let $\Gamma \subseteq \mathcal{Q}$. If a saddle-point in the restricted game $\mathcal{G}^\Gamma$ exists, it will have the form $(P^*, P^*)$. The distribution $P^*$ will then solve each of the following problems*:

(i) *Maximize over $P \in \Gamma$ the generalized entropy $H(P) \equiv S(P, P)$.*
(ii) *Minimize over $Q \in \mathcal{Q}$ the worst-case expected score, $\sup_{P \in \Gamma} S(P, Q)$.*

It is notable that, when Corollary 4.2 applies, the robust Bayes distribution solving problem (ii) turns out to belong to $\Gamma$, even though this constraint was not imposed.

We see from Theorem 4.1 that, when a saddle-point exists, the robust Bayes problem reduces to a maximum entropy problem. This property can thus be regarded as an indirect justification for applying the maximum entropy procedure. In the light of Theorem 4.1, we shall be particularly interested in the sequel in characterizing those decision problems for which a saddle-point exists in the game $\mathcal{G}^\Gamma$.

4.3. *A special case.* A partial characterization of a saddle-point can be given in the special case that the family $\Gamma$ is *closed under conditioning*, in the sense that, for all $P \in \Gamma$ and $B \subseteq \mathcal{X}$ a measurable set such that $P(B) > 0$, $P_B$, the conditional distribution under $P$ for $X$ given $X \in B$, is also in $\Gamma$. This will hold, most importantly, when $\Gamma$ is the set of all distributions supported on $\mathcal{X}$ or on some measurable subset of $\mathcal{X}$.

For the following lemma, we suppose that there exists a saddle-point $(P^*, \zeta^*)$ in the game $\mathcal{G}^\Gamma$, and write $H^* = L(P^*, \zeta^*)$. In particular, we have $L(P, \zeta^*) \leq H^*$ for all $P \in \Gamma$. We introduce $U := \{x \in \mathcal{X} : L(x, \zeta^*) = H^*\}$.

LEMMA 4.2. *Suppose that $\Gamma$ is closed under conditioning and that $P \in \Gamma$ is such that $L(P, \zeta^*) = H^*$. Then $P$ is supported on $U$.*

PROOF. Take $h < H^*$, and define $B := \{x \in \mathcal{X} : L(x, \zeta^*) \leq h\}$, $\pi := P(B)$. By linearity, we have $H^* = L(P, \zeta^*) = \pi L(P_B, \zeta^*) + (1 - \pi) L(P_{B^c}, \zeta^*)$ (where $B^c$ denotes the complement of $B$). However, by the definition of $B$, $L(P_B, \zeta^*) \leq h$, while (if $\pi \neq 1$) $L(P_{B^c}, \zeta^*) \leq H^*$, by Definition 4.1(b) and the fact that $P_{B^c} \in \Gamma$. It readily follows that $\pi = 0$. Since this holds for any $h < H^*$, we must have $P\{L(X, \zeta^*) \geq H^*\} = 1$. However, $\mathrm{E}_P\{L(X, \zeta^*)\} = L(P, \zeta^*) = H^*$, and the result follows. □

COROLLARY 4.3. $L(X, \zeta^*) = H^*$ *almost surely under $P^*$.*

COROLLARY 4.4. *If there exists $P \in \Gamma$ that is not supported on $U$, then $\zeta^*$ is not an equalizer rule in $\mathcal{G}^\Gamma$.*

Corollary 4.4 will apply, in particular, when $\Gamma$ is the family of all distributions supported on a subset $A$ of $\mathcal{X}$ and (as will generally be the case) $A$ is not a subset of $U$. Furthermore, since $\Gamma$ then contains the point mass at $x \in A$, we must have $L(x, \zeta^*) \leq H^*$, all $x \in A$, so that $U$ is the subset of $A$ on which the function $L(\cdot, \zeta^*)$ attains its maximum. In a typical such problem having a continuous sample space, the maxima of this function will be isolated points, and then we deduce that the maximum entropy distribution $P^*$ will be discrete (and the robust Bayes act $\zeta^*$ will *not* be an equalizer rule).

**5. An elementary minimax theorem.** Throughout this section we suppose that $\mathcal{X} = \{x_1, \ldots, x_N\}$ is finite and that $L$ is bounded. In particular, $L(P, a)$ and $H(P)$ are finite for all distributions $P$ over $\mathcal{X}$, and the set $\mathcal{P}$ of these distributions can be identified with the unit simplex in $\mathcal{R}^N$. We endow $\mathcal{P}$ with the topology inherited from this identification.

In this case we can show the existence of a saddle-point under some simple conditions. The following result is a variant of von Neumann's original minimax theorem [von Neumann (1928)]. It follows immediately from the general minimax theorem of Corollary A.1, whose conditions are here readily verified.

THEOREM 5.1. *Let $\Gamma$ be a closed convex subset of $\mathcal{P}$. Then the restricted game $\mathcal{G}^\Gamma$ has a finite value $H^*$, and the entropy $H(P)$ achieves its maximum $H^*$ over $\Gamma$ at some distribution $P^* \in \Gamma$.*

Theorem 5.1 does not automatically ensure the existence of a robust Bayes act. For this we impose a further condition on the action space. This involves the *risk-set* $S$ of the unrestricted game $\mathcal{G}$, that is, the convex subset of $\mathcal{R}^N$ consisting of all points $l(\zeta) := (L(x_1, \zeta), \ldots, L(x_N, \zeta))$ arising as the risk function of some possibly randomized act $\zeta \in \mathcal{Z}$.

THEOREM 5.2. *Suppose that $\Gamma$ is convex, and that the unrestricted risk-set $S$ is closed. Then there exists a robust Bayes act $\zeta^* \in \mathcal{Z}$. Moreover, there exists $P^*$ in the closure $\overline{\Gamma}$ of $\Gamma$ such that $\zeta^*$ is Bayes against $P^*$ and $(P^*, \zeta^*)$ is a saddle-point in the game $\mathcal{G}^{\overline{\Gamma}}$.*

PROOF. First assume $\Gamma$ closed. By Theorem 5.1 the game $\mathcal{G}^\Gamma$ has a finite value $H^*$. Then there exists a sequence $(\zeta_n)$ in $\mathcal{Z}$ such that $\lim_{n\to\infty} \sup_{P\in\Gamma} L(P, \zeta_n) = \inf_{\zeta\in\mathcal{Z}} \sup_{P\in\Gamma} L(P, \zeta) = H^*$. Since $S$ is compact, on taking a subsequence if necessary we can find $\zeta^* \in \mathcal{Z}$ such that $l(\zeta_n) \to l(\zeta^*)$. Then, for all $Q \in \Gamma$,

$$(47) \qquad L(Q, \zeta^*) = \lim_{n\to\infty} L(Q, \zeta_n) \leq \lim_{n\to\infty} \sup_{P\in\Gamma} L(P, \zeta_n) = H^*,$$

whence

$$(48) \qquad \sup_{P\in\Gamma} L(P, \zeta^*) \leq H^*.$$

However, for $P = P^*$, as given by Theorem 5.1, we have $L(P^*, \zeta^*) \geq H(P^*) = H^*$, so that $L(P^*, \zeta^*) = H^*$. The result now follows.

If $\Gamma$ is not closed, we can apply the above argument with $\Gamma$ replaced by $\overline{\Gamma}$ to obtain $\zeta^* \in \mathcal{Z}$ and $P^* \in \overline{\Gamma}$. Then $\sup_{\overline{\Gamma}} L(P, \zeta^*) \leq \sup_{\overline{\Gamma}} L(P, \zeta)$, all $\zeta \in \mathcal{Z}$. Since $L(P, \zeta)$ is linear, hence continuous, in $P$ for all $\zeta$, $\sup_\Gamma L(P, \zeta) = \sup_{\overline{\Gamma}} L(P, \zeta)$, and the general result follows. □

Note that $S$ is the convex hull of $S_0$, the set of risk functions of nonrandomized acts. A sufficient condition for $S$ to be closed is that $S_0$ be closed. In particular this will always hold if $\mathcal{A}$ is finite.

The above theorem gives a way of restricting the search for a robust Bayes act $\zeta^*$: first find a distribution $P^*$ maximizing the entropy over $\overline{\Gamma}$, then look for acts that are Bayes against $P^*$. In some cases this will yield a unique solution, and we are done. However, as will be seen below, this need not always be the case, and then further principles may be required.

5.1. *Examples.*

5.1.1. *Brier score.* Consider the Brier score (16) for $\mathcal{X} = \{0, 1\}$ and $\Gamma = \mathcal{P}$. Let $H$ be the corresponding entropy as in (19). From Figure 1, or directly, we see that the entropy is maximized for $P^*$ having $p^*(0) = p^*(1) = 1/2$. Since the Brier score is $\mathcal{P}$-strictly proper, the unique Bayes act against $P^*$ is $P^*$ itself. It follows that $P^*$ is the robust Bayes act against $\Gamma$. Hence in this case we can find the robust Bayes act simply by maximizing the entropy.

5.1.2. *Zero–one loss.* Now consider the zero–one loss (22) for $\mathcal{X} = \{0, 1\}$ and $\Gamma = \mathcal{P}$. Let $H$ be the corresponding entropy as in (25). From Figure 1, or directly, we see that the entropy is again maximized for $P^*$ with $p^*(0) = p^*(1) = 1/2$. However, in contrast to the case of the Brier score, $P^*$ now has several Bayes acts. In fact, *every* distribution $\zeta$ over $\mathcal{A} = \{0, 1\}$ is Bayes against $P^*$—yet only one of them (namely, $\zeta^* = P^*$) is robust Bayes. Therefore finding the maximum entropy $P^*$ is of no help whatsoever in finding the robust Bayes act $\zeta^*$ here. As we shall see in Section 7.6.3, however, this does not mean that the procedure described here (find a robust Bayes act by first finding the maximum entropy $P^*$ and then determine the Bayes acts of $P^*$) is never useful for zero–one loss: if $\Gamma \neq \mathcal{P}$, it may help in finding $\zeta^*$ after all.

**6. More general minimax theorems.** We are now ready to formulate more general minimax theorems. The proofs are given in the Appendix.

Let $(\mathcal{X}, \mathcal{B})$ be a metric space together with its Borel $\sigma$-algebra. Recall [Billingsley (1999), Section 5] that a family $\Gamma$ of distributions on $(\mathcal{X}, \mathcal{B})$ is called (*uniformly*) *tight* if, for all $\varepsilon > 0$, there exists a compact set $C \in \mathcal{B}$ such that $P(C) > 1 - \varepsilon$ for all $P \in \Gamma$.

THEOREM 6.1. *Let $\Gamma \subseteq \mathcal{P}$ be a convex, weakly closed and tight set of distributions. Suppose that for each $a \in \mathcal{A}$ the loss function $L(x, a)$ is bounded above and upper semicontinuous in $x$. Then the restricted game $\mathcal{G}^\Gamma = (\Gamma, \mathcal{A}, L)$ has a value. Moreover, a maximum entropy distribution $P^*$, attaining*

$$\sup_{P \in \Gamma} \inf_{a \in \mathcal{A}} L(P, a),$$

*exists.*

We note that if $\mathcal{X}$ is finite or countable and endowed with the discrete topology, then $L(x, a)$ is automatically a continuous, hence upper semicontinuous, function of $x$.

Theorem 6.1 cannot be applied to the logarithmic score, which is not bounded above in general. In such cases we may be able to use the theorems below. Note that these all refer to possibly randomized Bayes acts $\zeta^*$, but by Proposition 3.1 it will always be possible to choose such acts to be nonrandomized.

THEOREM 6.2. *Let $\Gamma \subseteq \mathcal{P}$ be convex, and let $P^* \in \Gamma$, with Bayes act $\zeta^*$, be such that $-\infty < H(P^*) = H^* := \sup_{P \in \Gamma} H(P) < \infty$. Suppose that for all $P \in \Gamma$ there exists $P_0 \in \mathcal{P}$ such that, on defining $Q_\lambda := (1 - \lambda)P_0 + \lambda P$, the following hold:*

  (i) $P^* = Q_{\lambda^*}$ *for some* $\lambda^* \in (0, 1)$.

(ii) *The function $H(Q_\lambda)$ is differentiable at $\lambda = \lambda^*$.*

*Then $(P^*, \zeta^*)$ is a saddle-point in $\mathcal{G}^\Gamma$.*

Theorem 6.2 essentially gives differentiability of the entropy as a condition for the existence of a saddle-point. This condition is strong but often easy to check. We now introduce a typically weaker condition, which may, however, be harder to check.

CONDITION 6.1. Let $(Q_n)$ be a sequence of distributions in $\Gamma$, with respective Bayes acts $(\zeta_n)$, such that the sequence $(H(Q_n))$ is bounded below and $(Q_n)$ converges weakly to some distribution $Q_0 \in \mathcal{P}_0$. Then we require that $Q_0 \in \mathcal{P}$, $Q_0$ has a Bayes act $\zeta_0$ and, for some choice of the Bayes acts $(\zeta_n)$ and $\zeta_0$, $L(P, \zeta_0) \leq \liminf_{n \to \infty} L(P, \zeta_n)$ for all $P \in \Gamma$.

One would typically aim to demonstrate Condition 6.1 in its stronger "$\Gamma$-free" form, wherein all mentions of $\Gamma$ are replaced by $\mathcal{P}$, or both $\Gamma$ and $\mathcal{P}$ are replaced by some family $\mathcal{Q}$ with $\Gamma \subseteq \mathcal{Q} \subseteq \mathcal{P}$. In particular, in the case of a $\mathcal{Q}$-proper scoring rule $S$, Condition 6.1 is implied by the following.

CONDITION 6.2. Let $(Q_n)$ be a sequence of distributions in $\mathcal{Q}$ such that the sequence $(H(Q_n))$ is bounded below and $(Q_n)$ converges weakly to $Q_0$. Then we require $Q_0 \in \mathcal{Q}$ and $S(P, Q_0) \leq \liminf_{n \to \infty} S(P, Q_n)$ for all $P \in \mathcal{Q}$.

This displays the condition as one of weak lower semicontinuity of the score in its second argument.

We shall further consider the following possible conditions on $\Gamma$:

CONDITION 6.3. $\Gamma$ is convex; every $P \in \Gamma$ has a Bayes act $\zeta_P$ and finite entropy $H(P)$; and $H^* := \sup_{P \in \Gamma} H(P) < \infty$.

CONDITION 6.4. Furthermore, there exists $P^* \in \Gamma$ with $H(P^*) = H^*$.

THEOREM 6.3. *Suppose Conditions 6.1, 6.3 and 6.4 hold. Then there exists $\zeta^* \in \mathcal{Z}$ such that $(P^*, \zeta^*)$ is a saddle-point in the game $\mathcal{G}^\Gamma$.*

If $H(P)$ is not upper-semicontinuous or if $\Gamma$ is not closed in the weak topology, then $\sup_{P \in \Gamma} H(P)$ may not be achieved. As explained in the Appendix, for a general sample space these are both strong requirements. If they do not hold, then Theorem 6.3 will not be applicable. In that case we may instead be able to apply Theorem 6.4:

THEOREM 6.4. *Suppose Conditions* 6.1 *and* 6.3 *hold and, in addition,* $\Gamma$ *is tight. Then there exists* $\zeta^* \in \mathcal{Z}$ *such that*

$$(49) \qquad \sup_{P \in \Gamma} L(P, \zeta^*) = \inf_{\zeta \in \mathcal{Z}} \sup_{P \in \Gamma} L(P, \zeta) = \sup_{P \in \Gamma} \inf_{a \in \mathcal{A}} L(P, a) = H^*.$$

*In particular, the game* $\mathcal{G}^\Gamma$ *has value* $H^*$, *and* $\zeta^*$ *is robust Bayes against* $\Gamma$.

In the Appendix we prove the more general Theorem A.2, which implies Theorem 6.4. We also prove Proposition A.1, which shows that (under some restrictions) the conditions of Theorem A.2 are satisfied when $L$ is the logarithmic score.

The theorems above supply sufficient conditions for the existence of a robust Bayes act, but do not give any further characterization of it, nor do they assist in finding it. In the next sections we shall consider the important special case of $\Gamma$ defined by linear constraints, for which we can develop explicit characterizations.

**7. Mean-value constraints.** Let $T \equiv t(X)$, with $t : \mathcal{X} \to \mathcal{R}^k$, be a fixed real- or vector-valued statistic. An important class of problems arises on imposing *mean-value constraints*, where we take

$$(50) \qquad \Gamma = \Gamma_\tau := \{P \in \mathcal{P} : \mathrm{E}_P(T) = \tau\},$$

for some $\tau \in \mathcal{R}^k$. This is the type of constraint for which the maximum entropy and minimum relative entropy principles have been most studied [Jaynes (1957a, b) and Csiszár (1975)].

We denote the associated restricted game $(\Gamma_\tau, \mathcal{A}, L)$ by $\mathcal{G}^\tau$. We call $T$ the *generating statistic*.

In some problems of this type (e.g., with logarithmic score on a continuous sample space), the family $\Gamma_\tau$ will be so large that the conditions of the theorems of Section 6 will not hold. Nevertheless, the special linear structure will often allow other arguments for showing the existence of a saddle-point.

7.1. *Duality.* Before continuing our study of saddle-points, we note some simple duality properties of such mean-value problems.

DEFINITION 7.1. The *specific entropy function* $h : \mathcal{R}^k \to [-\infty, \infty]$ (associated with the loss function $L$ and generating statistic $T$) is defined by

$$(51) \qquad h(\tau) := \sup_{P \in \Gamma_\tau} H(P).$$

In particular, if $\Gamma_\tau = \varnothing$, then $h(\tau) = -\infty$.

Now define $\mathcal{T} := \{\tau \in \mathcal{R}^k : h(\tau) > -\infty\}$ and $\mathcal{P}^* := \{P \in \mathcal{P} : \mathrm{E}_P(T) \in \mathcal{T}\}$.

LEMMA 7.1. *The set $\mathcal{T} \subseteq \mathcal{R}^k$ is convex, and the function $h$ is concave on $\mathcal{T}$.*

PROOF. Take $\tau_0, \tau_1 \in \mathcal{T}$ and $\lambda \in (0,1)$, and let $\tau_\lambda := (1-\lambda)\tau_0 + \lambda \tau_1$. There exist $P_0, P_1 \in \mathcal{P}$ with $P_i \in \Gamma_{\tau_i}$ and $H(P_i) > -\infty$, $i = 0, 1$. Let $P_\lambda := (1-\lambda)P_0 + \lambda P_1$. Then, for any $a \in \mathcal{A}$, $L(P_i, a) \geq H(P_i) > -\infty$, so that $L(P_\lambda, a) = (1-\lambda)L(P_0, a) + \lambda L(P_1, a)$ is defined, that is, $P_\lambda \in \mathcal{P}$. Moreover, clearly $P_\lambda \in \Gamma_{\tau_\lambda}$. We thus have $h(\tau_\lambda) \geq H(P_\lambda) \geq (1-\lambda)H(P_0) + \lambda H(P_1) > -\infty$. Thus $\tau_\lambda \in \mathcal{T}$; that is, $\mathcal{T}$ is convex. Now letting $P_0$ and $P_1$ vary independently, we obtain $h(\tau_\lambda) \geq (1-\lambda)h(\tau_0) + \lambda h(\tau_1)$; that is, $h$ is concave. □

For $\tau \in \mathcal{T}$, define

$$P_\tau := \arg\sup_{P \in \Gamma_\tau} H(P) \tag{52}$$

whenever this supremum is finite and is attained. It is allowed that $P_\tau$ is not unique, in which case we consider an arbitrary such maximizer. Then $H(P_\tau) = h(\tau)$. By Theorem 4.1(iv), (52) will hold if $(P_\tau, \zeta_\tau)$ is a saddle-point in $\mathcal{G}^\tau$.

Dually, for $\beta \in \mathcal{R}^k$, we introduce

$$Q_\beta := \arg\sup_{P \in \mathcal{P}^*}\{H(P) - \beta^{\mathrm{T}} \mathrm{E}_P(T)\}, \tag{53}$$

whenever this supremum is finite and is attained. Again, $Q_\beta$ is not necessarily unique. For any such $Q_\beta$ we can define a corresponding value of $\tau$ by

$$\tau = \mathrm{E}_{Q_\beta}(T). \tag{54}$$

Then $Q_\beta \in \Gamma_\tau$, and on restricting the supremum in (53) to $P \in \Gamma_\tau$, we see that we can take $Q_\beta$ for $P_\tau$ in (52). More generally, we write $\tau \leftrightarrow \beta$ whenever there is a common distribution that can serve as both $P_\tau$ in (52) and $Q_\beta$ in (53) (in cases of nonuniqueness this correspondence may not define a function in either direction).

It follows easily from (53) that, when $\tau \leftrightarrow \beta$,

$$h(\sigma) - \beta^{\mathrm{T}}\sigma \leq h(\tau) - \beta^{\mathrm{T}}\tau, \tag{55}$$

or equivalently

$$h(\sigma) \leq h(\tau) + \beta^{\mathrm{T}}(\sigma - \tau) \tag{56}$$

for all $\sigma \in \mathcal{T}$. Equation (56) expresses the fact that the hyperplane through the point $(\tau, h(\tau))$ with slope coefficients $\beta$ is a supporting hyperplane to the concave function $h: \mathcal{T} \to \mathcal{R}$. Thus $\tau$ and $\beta$ can be regarded as dual coordinates for the specific entropy function. In particular, if $\tau \leftrightarrow \beta$ and $h$ is differentiable at $\tau$, we must have

$$\beta = h'(\tau). \tag{57}$$

More generally, if $\tau_1 \leftrightarrow \beta_1$ and $\tau_2 \leftrightarrow \beta_2$, then on combining two applications of (55) we readily obtain

$$(\tau_2 - \tau_1)^{\mathrm{T}}(\beta_2 - \beta_1) \leq 0. \tag{58}$$

In particular, when $k = 1$ the correspondence $\tau \leftrightarrow \beta$ is nonincreasing in the sense that $\tau_2 > \tau_1 \Rightarrow \beta_2 \leq \beta_1$.

7.2. *Linear loss condition.* Theorem 7.1 gives a simple sufficient condition for an act to be robust Bayes against $\Gamma_\tau$ of the form (50). We first introduce the following definition.

DEFINITION 7.2. An act $\zeta \in \mathcal{Z}$ is *linear* (with respect to loss function $L$ and statistic $T$) if, for some $\beta_0 \in \mathcal{R}$ and $\beta = (\beta_1, \ldots, \beta_k)^{\mathrm{T}} \in \mathcal{R}^k$ and all $x \in \mathcal{X}$,

$$L(x, \zeta) = \beta_0 + \beta^{\mathrm{T}} t(x). \tag{59}$$

A distribution $P \in \mathcal{P}$ is *linear* if it has a Bayes act $\zeta$ that is linear. In this case we call $(P, \zeta)$ a *linear pair*. If $\mathrm{E}_P(T) = \tau$ is finite, we then call $\tau$ a *linear point* of $\mathcal{T}$. In all cases we call $(\beta_0, \beta)$ the associated *linear coefficients*.

Note that, if the problem is formulated in terms of a $\mathcal{Q}$-strictly proper scoring rule $S$, and $P \in \mathcal{Q}$, the conditions "$P$ is a linear distribution," "$P$ is a linear act" and "$(P, P)$ is a linear pair" are all equivalent, holding when we have

$$S(x, P) = \beta_0 + \sum_{j=1}^{k} \beta_j \, t_j(x) \tag{60}$$

for all $x \in \mathcal{X}$.

THEOREM 7.1. *Let $\tau \in \mathcal{T}$ be linear, with associated linear pair $(P_\tau, \zeta_\tau)$ and linear coefficients $(\beta_0, \beta)$. Let $\Gamma_\tau$ be given by (50). Then the following hold:*

(i) $\zeta_\tau$ *is an equalizer rule against* $\Gamma_\tau$.
(ii) $(P_\tau, \zeta_\tau)$ *is a saddle-point in* $\mathcal{G}^\tau$.
(iii) $\zeta_\tau$ *is robust Bayes against* $\Gamma_\tau$.
(iv) $h(\tau) = H(P_\tau) = \beta_0 + \beta^{\mathrm{T}} \tau$.
(v) $\tau \leftrightarrow \beta$.

PROOF. For any $P \in \mathcal{P}^*$ we have

$$L(P, \zeta_\tau) = \beta_0 + \beta^{\mathrm{T}} \mathrm{E}_P(T). \tag{61}$$

By (61) $L(P, \zeta_\tau) = \beta_0 + \beta^{\mathrm{T}} \tau = L(P_\tau, \zeta_\tau)$ for all $P \in \Gamma$. Thus (44)(b) holds with equality, showing (i). Since $L(P_\tau, \zeta_\tau)$ is finite and $\zeta_\tau$ is Bayes against $P_\tau$, (44)(a) holds. We have thus shown (ii). Then (iii) follows from Theorem 4.1(v), and (iv) follows from Theorem 4.1(i), (iii) and (iv). For (v), we have from (61) that, for $P \in \mathcal{P}^*$,

$$H(P) - \beta^{\mathrm{T}} \mathrm{E}_P(T) \leq L(P, \zeta_\tau) - \beta^{\mathrm{T}} \mathrm{E}_P(T) \tag{62}$$

$$= \beta_0 \tag{63}$$

$$= H(P_\tau) - \beta^{\mathrm{T}} \mathrm{E}_{P_\tau}(T) \tag{64}$$

from (iv). Thus we can take $Q_\beta$ in (53) to be $P_\tau$. □

COROLLARY 7.1. *The same result holds if* (59) *is only required to hold with probability* 1 *under every* $P \in \Gamma_\tau$.

We now develop a partial converse to Theorem 7.1, giving a necessary condition for a saddle-point. This will be given in Theorem 7.2.

DEFINITION 7.3. A point $\tau \in \mathcal{T}$ is *regular* if there exists a saddle-point $(P_\tau, \zeta_\tau)$ in $\mathcal{G}^\tau$, and there exists $\beta = (\beta_1, \ldots, \beta_k)^{\mathrm{T}} \in \mathcal{R}^k$ such that:

(i) $P_\tau$ can serve as $Q_\beta$ in (53) (so that $\tau \leftrightarrow \beta$).
(ii) With $\zeta = \zeta_\tau$ and (necessarily)

$$\beta_0 := h(\tau) - \beta^{\mathrm{T}} \tau, \tag{65}$$

the linear loss property (59) holds with $P_\tau$-probability 1.

If $\tau$ satisfies the conditions of Theorem 7.1 or of Corollary 7.1 it will be regular, but in general the force of the "almost sure" linearity requirement in (ii) above is weaker than needed for Corollary 7.1.

We shall denote the set of regular points of $\mathcal{T}$ by $\mathcal{T}^r$, and its subset of linear points by $\mathcal{T}^l$. For discrete $\mathcal{X}$, $\tau \in \mathcal{T}^r$ will by (ii) be linear whenever $P_\tau$ gives positive probability to every $x \in \mathcal{X}$. More generally, as soon as we know $\tau \in \mathcal{T}^r$, the following property, which follows trivially from (ii), can be used to simplify the search for a saddle-point:

LEMMA 7.2. *If $\tau$ is regular, the support $\mathcal{X}_\tau$ of $P_\tau$ is such that, for some $\zeta \in \mathcal{Z}$, $L(x, \zeta)$ is a linear function of $t(x)$ on $\mathcal{X}_\tau$.*

The following lemma and corollary are equally trivial.

LEMMA 7.3. *Suppose $\tau \in \mathcal{T}^r$. If $P \in \Gamma_\tau$ and $P \ll P_\tau$, then $L(P, \zeta_\tau) = h(\tau)$.*

COROLLARY 7.2. *If $\tau \in \mathcal{T}^r$ and $P \ll P_\tau$ for all $P \in \Gamma_\tau$, then $\zeta_\tau$ is an equalizer rule in $\mathcal{G}^\tau$.*

We now show that, under mild conditions, a point $\tau$ in the relative interior [Rockafellar (1970), page 44] $\mathcal{T}^0$ of $\mathcal{T}$ will be regular. Fix $\tau \in \mathcal{T}^0$ and consider $\Gamma_\tau$, given by (50). We shall suppose that there exists a saddle-point $(P_\tau, \zeta_\tau)$ for the game $\mathcal{G}^\tau$—this could be established by the theory of Section 5 or 6, for example. The value $L(P_\tau, \zeta_\tau)$ of the game will then be $h(\tau)$, which will be finite.

Consider the function $\psi_\tau$ on $\mathcal{T}$ defined by

(66) $$\psi_\tau(\sigma) := \sup_{P \in \Gamma_\sigma} L(P, \zeta_\tau).$$

In particular, $\psi_\tau(\tau) = h(\tau)$.

PROPOSITION 7.1. *$\psi_\tau$ is finite and concave on $\mathcal{T}$.*

PROOF. For $\sigma \in \mathcal{T}$ there exists $P \in \Gamma_\sigma$ with $H(P) > -\infty$; so $\psi_\tau(\sigma) \geq L(P, \zeta_\tau) \geq H(P) > -\infty$.

Now take $\sigma_0, \sigma_1 \in \mathcal{T}$ and $\lambda \in (0,1)$, and consider $\sigma_\lambda := (1-\lambda)\sigma_0 + \lambda\sigma_1$. Then $\Gamma_{\sigma_\lambda} \supseteq \{(1-\lambda)P_0 + \lambda P_1 : P_0 \in \Gamma_{\sigma_0}, P_1 \in \Gamma_{\sigma_1}\}$, so that $\psi_\tau(\sigma_\lambda) \geq (1-\lambda) \times \psi_\tau(\sigma_0) + \lambda \psi_\tau(\sigma_1)$. Thus $\psi_\tau$ is concave on $\mathcal{T}$.

Finally, if $\psi_\tau$ were to take the value $+\infty$ anywhere on $\mathcal{T}$, then by Lemma 4.2.6 of Stoer and Witzgall (1970) it would do so at $\tau \in \mathcal{T}^0$, which is impossible since $\psi_\tau(\tau) = h(\tau)$ has been assumed finite. □

For the proof of Theorem 7.2 we need to impose a condition allowing the passage from (70) to (71). For the examples considered in this paper, we can use the simplest such condition:

CONDITION 7.1. *For all $x \in \mathcal{X}$, $t(x) \in \mathcal{T}$.*

This is equivalent to $t(\mathcal{X}) \subseteq \mathcal{T}$, or in turn to $\mathcal{T}$ being the convex hull of $t(\mathcal{X})$. For other applications (e.g., involving unbounded loss functions on continuous sample spaces) this may not hold, and then alternative conditions may be more appropriate.

THEOREM 7.2. *Suppose that $\tau \in \mathcal{T}^0$ and $(P_\tau, \zeta_\tau)$ is a saddle-point for the game $\mathcal{G}^\tau$. If Condition 7.1 holds, then $\tau$ is regular.*

PROOF. $\mathcal{T}$ is convex, $\psi_\tau : \mathcal{T} \to \mathcal{R}$ is concave, and $\tau \in \mathcal{T}^0$. The supporting hyperplane theorem [Stoer and Witzgall (1970), Corollary 4.2.9] then implies that there exists $\beta \in \mathcal{R}^k$ such that, for all $\sigma \in \mathcal{T}$,

$$\psi_\tau(\tau) + \beta^{\mathrm{T}}(\sigma - \tau) \geq \psi_\tau(\sigma). \tag{67}$$

That is, for any $P \in \mathcal{P}^*$,

$$h(\tau) + \beta^{\mathrm{T}}\{\mathrm{E}_P(T) - \tau\} \geq \psi_\tau\{\mathrm{E}_P(T)\}. \tag{68}$$

However, for $P \in \mathcal{P}^*$,

$$\psi_\tau\{\mathrm{E}_P(T)\} \geq L(P, \zeta_\tau) \geq \inf_\zeta L(P, \zeta) = H(P). \tag{69}$$

Thus, for all $P \in \mathcal{P}^*$,

$$h(\tau) + \beta^{\mathrm{T}}\{\mathrm{E}_P(T) - \tau\} \geq H(P),$$

with equality when $P = P_\tau$. This yields Definition 7.3(i).

For (ii), (68) and (69) imply that

$$h(\tau) - L(P, \zeta_\tau) + \beta^{\mathrm{T}}\{\mathrm{E}_P(T) - \tau\} \geq 0 \qquad \text{for all } P \in \mathcal{P}^*. \tag{70}$$

Take $x \in \mathcal{X}$, and let $P_x$ be the point mass on $x$. By Condition 7.1, $P_x \in \mathcal{P}^*$, and so

$$h(\tau) - L(x, \zeta_\tau) + \beta^{\mathrm{T}}\{t(x) - \tau\} \geq 0 \qquad \text{for all } x \in \mathcal{X}. \tag{71}$$

On the other hand,

$$\mathrm{E}_{P_\tau}[h(\tau) - L(X, \zeta_\tau) + \beta^{\mathrm{T}}\{t(X) - \tau\}] = 0. \tag{72}$$

Together (71) and (72) imply that

$$P_\tau[h(\tau) - L(X, \zeta_\tau) + \beta^{\mathrm{T}}\{t(X) - \tau\} = 0] = 1. \tag{73}$$

The result follows. □

7.3. *Exponential families.* Here we relate the above theory to familiar properties of exponential families [Barndorff-Nielsen (1978)].

Let $\mu$ be a fixed $\sigma$-finite measure on a suitable $\sigma$-algebra in $\mathcal{X}$. The set of all distributions $P \ll \mu$ having a $\mu$-density $p$ that can be expressed in the form

$$p(x) = \exp\left\{\alpha_0 + \sum_{j=1}^{k} \alpha_j\, t_j(x)\right\} \tag{74}$$

for all $x \in \mathcal{X}$ is the *exponential family* $\mathcal{E}$ generated by the base measure $\mu$ and the statistic $T$.

We remark that (74) is trivially equivalent to

(75) $$S(x,p) = \beta_0 + \sum_{j=1}^{k} \beta_j\, t_j(x),$$

for all $x \in \mathcal{X}$, where $S$ is the logarithmic score (20), and $\beta_j = -\alpha_j$. In particular, $(P, p)$ is a linear pair.

Now under regularity conditions on $\mu$ and $T$ [Barndorff-Nielsen (1978), Chapter 9; see also Section 7.4.1 below], for all $\tau \in \mathcal{T}^0$ there will exist a unique $P_\tau \in \Gamma_\tau \cap \mathcal{E}$; that is, $P_\tau$ has a density $p_\tau$ of the form (74), and $\mathrm{E}_{P_\tau}(T) = \tau$. Comparing (75) with (59), it follows from Theorem 7.1 that (as already demonstrated in detail in Section 2.1) $(P_\tau, p_\tau)$ is a saddle-point in $\mathcal{G}_\tau$. In particular, as is well known [Jaynes (1989)], the distribution $P_\tau$ will maximize the entropy (21), subject to the mean-value constraints (50). However, we regard this property as less fundamental than the concomitant dual property: that $p_\tau$ is the robust Bayes act under the logarithmic score when all that we know of Nature's distribution $P$ is that it satisfies the mean-value constraint $P \in \Gamma_\tau$. Furthermore, by Theorem 7.1(i), in this case $p_\tau$ will be an equalizer strategy against $\Gamma_\tau$ [cf. (3)].

We remark that $p_\tau$ of the form (74) is only one version of the density for $P_\tau$ with respect to $\mu$; any other such density can differ from $p_\tau$ on a set of $\mu$-measure 0. However, our game requires DM to specify a density, rather than a distribution, and from this point of view certain other versions of the density of $P_\tau$ (which are of course still Bayes against $P_\tau$) will not do: they are not robust Bayes. For example, let $\mathcal{X} = \mathcal{R}$, let $\mu =$ Lebesgue measure and consider the constraints $\mathrm{E}_P(X) = 0$, $\mathrm{E}_P(X^2) = 1$. Let $P_0$ be the standard Normal distribution $N(0,1)$, and let $p_0$ be its usual density formula, $p_0(x) = (2\pi)^{-1/2} \exp -\frac{1}{2}x^2$. Then the conditions of Theorem 7.1 hold, $P_0$ is maximum entropy (as is well known) and the choice $p_0$ for its density is robust Bayes against the set $\Gamma_0$ of *all* distributions $P$—including, importantly, discrete distributions—that satisfy the constraints. This would not have been true if instead of $p_0$ we had taken $p'_0$, identical with $p_0$ except for $p'_0(x) = p_0(x)/2$ at $x = \pm 1$. While $p'_0$ is still Bayes against $P_0$, its Bayes loss against the distribution in $\Gamma_0$ that puts equal probability $1/2$ at $-1$ and $+1$ exceeds the (constant) Bayes loss of $p_0$ by $\log 2$. Consequently, $p'_0$ is *not* a robust Bayes act. It is in fact easy to see that a density $p$ will be robust Bayes in this problem if and only if $p(x) \geq p_0(x)$ everywhere (the set on which strict inequality holds necessarily having Lebesgue measure 0).

We further remark that none of the theorems of Section 6 applies to the above problem. The boundedness and weak closure requirements of Theorem 6.1 both fail; condition (ii) of Theorem 6.2 fails; and although Condition 6.2 holds, the existence of a Bayes act and finite entropy required for Condition 6.3 fail for those distributions in $\Gamma_\tau$ having a discrete component.

7.4. *Generalized exponential families.* We now show how our game-theoretic approach supports the extension of many of the concepts and properties of standard exponential family theory to apply to what we shall term a *generalized exponential family*, specifically tailored to the relevant decision problem. Although the link to exponentiation has now vanished, analogues of familiar duality properties of exponential families [Barndorff-Nielsen (1978), Chapter 9] can be based on the theory of Section 7.1.

Consider the following condition.

CONDITION 7.2. For all $\tau \in \mathcal{T}$, $h(\tau) = \sup_{P \in \Gamma_\tau} H(P)$ is finite and is achieved for a unique $P_\tau \in \Gamma_\tau$.

In particular, this will hold if (i) $\mathcal{X}$ is finite, (ii) $L$ is bounded and (iii) $H$ is strictly convex. For under (i) and (ii) Theorem 5.1 guarantees that a maximum generalized entropy distribution $P_\tau$ exists, which must then be unique by (iii).

Under Condition 7.2 we can introduce the following parametric family of distributions over $\mathcal{X}$:

$$\mathcal{E}^m := \{P_\tau : \tau \in \mathcal{T}\}. \tag{76}$$

We call $\mathcal{E}^m$ the *full generalized exponential family* generated by $L$ and $T$; and we call $\tau$ its *mean-value parameter*. Condition 7.2 ensures that the map $\tau \mapsto P_\tau$ is one-to-one.

Alternatively, consider the following condition:

CONDITION 7.3. For all $\beta \in \mathcal{R}^k$, $\sup_{P \in \mathcal{P}^*}\{H(P) - \beta^\mathrm{T} \mathrm{E}_P(T)\}$ is finite and is achieved for a unique distribution $Q_\beta \in \mathcal{P}^*$.

Again, this will hold if, in particular, (i)–(iii) below Condition 7.2 are satisfied.

Under Condition 7.3 we can introduce the parametric family

$$\mathcal{E}^n := \{Q_\beta : \beta \in \mathcal{R}^k\}. \tag{77}$$

We call this family the *natural generalized exponential family* generated by the loss function $L$ and statistic $T$; we call $\beta$ its *natural parameter*. This definition extends a construction of Lafferty (1999) based on Bregman divergence: see Section 8.4.2. Note that in general the natural parameter $\beta$ in $\mathcal{E}^n$ need not be identified; that is, the map $\beta \mapsto Q_\beta$ may not be one-to-one. See, however, Proposition 7.2, which sets limits to this nonidentifiability.

From this point on, we suppose that both Conditions 7.2 and 7.3 are satisfied. For any $\beta \in \mathcal{R}^k$, (54) yields $\tau \in \mathcal{T}$ with $\tau \leftrightarrow \beta$, that is, $P_\tau = Q_\beta$. It follows that $\mathcal{E}^n \subseteq \mathcal{E}^m$.

We further define $\mathcal{E}^r := \{P_\tau : \tau \in \mathcal{T}^r\}$, the *regular generalized exponential family*, and $\mathcal{E}^l := \{P_\tau : \tau \in \mathcal{T}^l\}$, the *linear generalized exponential family*, generated by $L$ and $T$. Then $\mathcal{E}^l \subseteq \mathcal{E}^r \subseteq \mathcal{E}^m$. In general, $\mathcal{E}^l$ may be a proper subset of $\mathcal{E}^r$: then for $P_\tau \in \mathcal{E}^r \setminus \mathcal{E}^l$ we can only assert the "almost sure linear loss" property of Lemma 7.2.

The following result follows immediately from Definition 7.3(ii).

PROPOSITION 7.2. *If* $Q_{\beta_1} = Q_{\beta_2} = Q \in \mathcal{E}^r$, *then* $(\beta_1 - \beta_2)^\mathrm{T} T = 0$ *almost surely under* $Q$.

For $\tau \in \mathcal{T}^r$ choose $\beta$ as in Definition 7.3. Then $\tau \leftrightarrow \beta$, and it follows that $\mathcal{E}^r \subseteq \mathcal{E}^n$. We have thus demonstrated the following.

PROPOSITION 7.3. *When Conditions* 7.2 *and* 7.3 *both apply,*
$$\mathcal{E}^r \subseteq \mathcal{E}^n \subseteq \mathcal{E}^m.$$

Now consider $\mathcal{E}^0 := \{P_\tau : \tau \in \mathcal{T}^0\}$, the *open generalized exponential family* generated by $L$ and $T$. From Theorem 7.2 we have the following:

PROPOSITION 7.4. *Suppose Conditions* 7.1–7.3 *all apply and a saddle-point exists in* $\mathcal{G}^\tau$ *for all* $\tau \in \mathcal{T}^0$. *Then*
$$\mathcal{E}^0 \subseteq \mathcal{E}^r \subseteq \mathcal{E}^n \subseteq \mathcal{E}^m. \tag{78}$$

7.4.1. *Application to standard exponential families.* We now consider more closely the relationship between the above theory and standard exponential family theory.

Let $\mathcal{E}^*$ be the standard exponential family (74) generated by some base measure $\mu$ and statistic $T$. Taking as our loss function the logarithmic score $S$, (75) shows that $\mathcal{E}^l \subseteq \mathcal{E}^*$ (distributions in $\mathcal{E}^* \setminus \mathcal{E}^l$ being those for which the expectation of $T$ does not exist). We can further ask: What is the relationship between $\mathcal{E}^*$ and $\mathcal{E}^n$? As a partial answer to this, we give sufficient conditions for $\mathcal{E}^*$, $\mathcal{E}^l$ and $\mathcal{E}^n$ to coincide.

For $\beta = (\beta_1, \ldots, \beta_k) \in \mathcal{R}^k$, define

$$\kappa(\beta) := \log \int e^{-\beta^\mathrm{T} t(x)} \, d\mu, \tag{79}$$

$$\chi(\beta) := \sup_{P \in \mathcal{P}^*} \{H(P) - \beta^\mathrm{T} \mathrm{E}_P(T)\}. \tag{80}$$

Let $\mathcal{B}$ denote the convex set $\{\beta \in \mathcal{R}^k : \kappa(\beta) < \infty\}$, and let $\mathcal{B}^0$ denote its relative interior. For $\beta \in \mathcal{B}$, let $Q^*_\beta$ be the distribution in $\mathcal{E}^*$ with $\mu$-density $q^*_\beta(x) := \exp\{-\kappa(\beta) - \beta^\mathrm{T} t(x)\}$, and let $Q_\beta$, if it exists, achieve the supremum in (80).

PROPOSITION 7.5. (i) *For all $\beta \in \mathcal{B}^0$, the act $q_\beta^*$ is linear, and $Q_\beta = Q_\beta^*$ uniquely. Moreover, $\chi(\beta) = \kappa(\beta)$.*

(ii) *If $\mathcal{B} = \mathcal{R}^k$, then* Condition 7.3 *holds and $\mathcal{E}^* = \mathcal{E}^l = \mathcal{E}^n$.*

(iii) *If Condition 7.3 holds, $\mathcal{B}$ is nonempty and $\mathcal{E}^*$ is minimal and steep, then $\mathcal{B} = \mathcal{R}^k$ and $\mathcal{E}^* = \mathcal{E}^l = \mathcal{E}^n$.*

[Note that the condition for (ii) will apply whenever the sample space $\mathcal{X}$ is finite.]

PROOF OF PROPOSITION 7.5. Linearity of the act $q_\beta^*$ ($\beta \in \mathcal{B}$) is immediate, the associated linear coefficients being $(\beta_0, \beta)$ with $\beta_0 = \kappa(\beta)$. Suppose $\beta \in \mathcal{B}^0$. Then $\tau := E_{Q_\beta^*}(T)$ exists [Barndorff-Nielsen (1978), Theorem 8.1]. We may also write $P_\tau$ for $Q_\beta^*$. Then $\tau$ is a linear point, with $(P_\tau, p_\tau)$ the associated linear pair. By Theorem 7.1(iv) $\kappa(\beta) = H(P_\tau) - \beta^\mathrm{T} \tau$. Also, by Theorem 7.1(v) we can take $P_\tau = Q_\beta^*$ as $Q_\beta$. The supremum in (80) thus being achieved by $P_\tau$, we have $\chi(\beta) = H(P_\tau) - \beta^\mathrm{T} \tau = \kappa(\beta)$.

To show that the supremum in (80) is achieved uniquely at $Q_\beta^*$, note that any $P$ achieving this supremum must satisfy

$$
\begin{aligned}
H(P) - \beta^\mathrm{T} \mathrm{E}_P(T) &= H(Q_\beta^*) - \beta^\mathrm{T} \mathrm{E}_{Q_\beta^*}(T) \\
&= \kappa(\beta) = S(P, q_\beta^*) - \beta^\mathrm{T} \mathrm{E}_P(T),
\end{aligned}
\tag{81}
$$

the last equality deriving from the definition of $q_\beta^*$. It follows that $S(P, q_\beta^*) = H(P) = S(P, p)$, whence $\int \log\{p(x)/q_\beta^*(x)\} p(x) \, d\mu = 0$. However, this can only hold if $P = Q_\beta^*$.

Part (ii) follows immediately.

For part (iii), assume Condition 7.3 holds. Then, for all $\beta \in \mathcal{R}^k$,

$$
\chi(\beta) = \sup_{\tau \in \mathcal{T}} \sup_{P \in \Gamma_\tau} \{H(P) - \beta^\mathrm{T} \tau\} = \sup_{\tau \in \mathcal{T}} \{h(\tau) - \beta^\mathrm{T} \tau\},
\tag{82}
$$

with $h(\tau)$ as in (51). By Lemma 7.1 $\mathcal{T}$ is convex. It follows that $\chi$ is a closed convex function on $\mathcal{R}^k$.

Steepness of $\mathcal{E}^*$ means that $|\kappa(\beta_n)| \to \infty$ whenever $(\beta_n)$ is a sequence in $\mathcal{B}^0$ converging to a relative boundary point $\beta^*$ of $\mathcal{B}$. Since $\kappa$ is convex [Barndorff-Nielsen (1978), Chapter 8] and $\chi$ coincides with $\kappa$ on $\mathcal{B}^0$, we must thus have $|\chi(\beta_n)| \to \infty$ as $(\beta_n) \to \beta^*$. Since by Condition 7.3 the closed convex function $\chi$ is finite on $\mathcal{R}^k$, $\mathcal{B}$ cannot have any relative boundary points—hence, under minimality, any boundary points—in $\mathcal{R}^k$. Since $\mathcal{B}$ is nonempty, it must thus coincide with $\mathcal{R}^k$. Then, by (ii) $\mathcal{E}^* = \mathcal{E}^l = \mathcal{E}^n$. □

To see that even under the above conditions we need not have $\mathcal{E}^* = \mathcal{E}^m$, consider the case $\mathcal{X} = \{0, 1\}$, $T = X$. Then $\mathcal{E}^m$ consists of all distributions on $\mathcal{X}$, whereas $\mathcal{E}^* = \mathcal{E}^l = \mathcal{E}^n$ excludes the one-point distributions at 0 and 1.

7.4.2. *Characterization of specific entropy.* We now generalize a result of Kivinen and Warmuth (1999). For the case of finite $\mathcal{X}$, they attack the problem of minimizing the Kullback–Leibler discrepancy $\mathrm{KL}(P, P_0)$ over all $P$ such that $\mathrm{E}_P(T) = 0$. Equivalently (see Section 3.5.2), they are maximizing the entropy $H(P) = -\mathrm{KL}(P, P_0)$, associated with the logarithmic score relative to base measure $P_0$, subject to $P \in \Gamma_0$.

Let $\mathcal{E}^*$ be the standard exponential family (74) generated by base measure $P_0$ and statistic $T$, with typical member $Q_\beta^*$ ($\beta \in \mathcal{R}^k$) having probability mass function of the form

$$q_\beta^*(x) = p_0(x) e^{-\kappa(\beta) - \beta^\mathrm{T} t(x)} \tag{83}$$

and entropy $h(\tau) = \kappa(\beta) + \beta^\mathrm{T} \tau$, where $\tau = \mathrm{E}_{Q_\beta}(T)$.

Suppose $0 \in \mathcal{T}^0$. By Chapter 9 of Barndorff-Nielsen (1978), there then exists within $\Gamma_0$ a unique member $Q_{\beta^*}^*$ of $\mathcal{E}^*$. By Theorem 7.1 the maximum of the entropy $-\mathrm{KL}(P, P_0)$ is achieved for $P = Q_{\beta^*}^*$; its maximized value is thus $h(0) = \kappa(\beta^*)$, where

$$\kappa(\beta) = \log \sum_x p_0(x) e^{-\beta^\mathrm{T} t(x)}. \tag{84}$$

Equation (1.5) of Kivinen and Warmuth (1999) essentially states that the maximized entropy $h(0)$ over $\Gamma_0$ can equivalently be obtained as

$$h(0) = \min_{\beta \in \mathcal{R}^k} \kappa(\beta). \tag{85}$$

By Proposition 7.5(i) this can also be written as

$$h(0) = \min_{\beta \in \mathcal{R}^k} \chi(\beta). \tag{86}$$

We now extend the above property to a more general decision problem, satisfying Conditions 7.2 and 7.3. Let $\tau \leftrightarrow \beta$, $\sigma \leftrightarrow \gamma$ ($\tau, \sigma \in \mathcal{T}$). Then $\chi(\beta) = \beta_0 = h(\tau) - \beta^\mathrm{T} \tau$, $\chi(\gamma) = \gamma_0 = h(\sigma) - \gamma^\mathrm{T} \sigma$, with $\beta_0$, and correspondingly $\gamma_0$, as in (65). From (56) we have

$$h(\sigma) \leq \beta_0 + \beta^\mathrm{T} \sigma. \tag{87}$$

Moreover, we have equality in (87) when $\beta = \gamma$. It follows that for $\sigma \in \mathcal{T}$

$$h(\sigma) = \inf_{\beta \in \mathcal{R}^k} \{\chi(\beta) + \beta^\mathrm{T} \sigma\}, \tag{88}$$

the infimum being attained when $\beta \leftrightarrow \sigma$. In particular, when $0 \in \mathcal{T}$ we recover (86) in this more general context. Equations (82) and (88) express a conjugacy relation between the convex function $\chi$ and the concave function $h$.

7.5. *Support.* Fix $x \in \mathcal{X}$. For any act $\zeta \in \mathcal{Z}$ we term the negative loss $s_x(\zeta) := -L(x, \zeta)$ the *support* for act $\zeta$ based on data $x$. Likewise, $s_P(\zeta) := -L(P, \zeta)$ is the support for $\zeta$ based on a (theoretical or empirical) distribution $P$ for $X$. If $\mathcal{F} \subseteq \mathcal{Z}$ is a family of contemplated acts, then the function $\zeta \mapsto s_P(\zeta)$ on $\mathcal{F}$ is the *support function* over $\mathcal{F}$ based on "data" $P$. When the maximum of $s_P(\zeta)$ over $\zeta \in \mathcal{F}$ is achieved at $\hat{\zeta} \in \mathcal{F}$, we may term $\hat{\zeta}$ the *maximum support act* (in $\mathcal{F}$, based on $P$). Then $\hat{\zeta}$ is just the Bayes act against $P$ in the game with loss function $L(x, \zeta)$, when $\zeta$ is restricted to the set $\mathcal{F}$.

For the special case of the logarithmic score (20), $s_x(q) = \log q(x)$ is the log-likelihood of a tentative explanation $q(\cdot)$, on the basis of data $x$; if $P$ is the empirical distribution formed from a sample of $n$ observations, $s_P(q)$ is ($n^{-1}$ times) the log-likelihood for the explanation whereby these were independently and identically generated from density $q(\cdot)$. Thus our definition of the support function generalizes that used in likelihood theory [Edwards (1992)], while our definition of maximum support act generalizes that of maximum likelihood estimate. In particular, maximum likelihood is Bayes in the sense of the previous paragraph.

Typically we are only interested in differences of support (between different acts, for fixed data $x$ or distribution $P$), so that we can regard this function as defined only up to an additive constant; this is exactly analogous to regarding a likelihood function as defined only up to a positive multiplicative constant.

7.5.1. *Maximum support in generalized exponential families.* Let $T \equiv t(X)$ be a statistic, and let $\mathcal{E}^r$ be the regular generalized exponential family generated by $L$ and $T$. Fix a distribution $P^*$ over $\mathcal{X}$, and consider the associated support function $s^*(\cdot) := s_{P^*}(\cdot)$ over the family $\mathcal{F}^r := \{\zeta_\tau : \tau \in \mathcal{T}^r\}$. It is well known [Barndorff-Nielsen (1978), Section 9.3] that, in the case of an ordinary exponential family (when $L$ is logarithmic score and $\mathcal{F}^r = \{p_\tau(\cdot) : \tau \in \mathcal{T}^r\}$ is the set of densities of distributions in $\mathcal{E}^r$), the likelihood over $\mathcal{F}^r$ based on data $x^*$ (or more generally on a distribution $P^*$) is under regularity conditions maximized at $p_{\tau^*}$, where $\tau^* = t(x^*)$ [or $\tau^* = \mathrm{E}_{P^*}(T)$]. The following result gives a partial generalization of this property.

THEOREM 7.3. *Suppose $\tau^* := \mathrm{E}_{P^*}(T) \in \mathcal{T}^r$. Let $\tau \in \mathcal{T}^r$ be such that either of the following holds:*

(i) *$\zeta_\tau$ is linear;*
(ii) *$P^* \ll P_\tau$.*

*Then*

$$s^*(\zeta_{\tau^*}) \geq s^*(\zeta_\tau). \tag{89}$$

PROOF. Since $P^* \in \Gamma_{\tau^*}$ and $(P_{\tau^*}, \zeta_{\tau^*})$ is a saddle-point in $\mathcal{G}^{\tau^*}$, we have

(90) $$s^*(\zeta_{\tau^*}) \geq -h(\tau^*).$$

Under (i), (59) holds everywhere; under (ii), by Definition 7.3(ii) it holds with $P^*$-probability 1. In either case we obtain

(91) $$L(P^*, \zeta_\tau) = h(\tau) + \beta^{\mathrm{T}}(\tau^* - \tau).$$

By (56), the right-hand side is at least as large as $h(\tau^*)$, whence $s^*(\zeta_\tau) \leq -h(\tau^*)$. Combining this with (90), the result follows. $\square$

COROLLARY 7.3. *If for all $\tau \in \mathcal{E}^r$ either $\zeta_\tau$ is linear or $P^* \ll P_\tau$, then $\zeta_{\tau^*}$ is the maximum support act in $\mathcal{F}^r$.*

For the case of the logarithmic score (20) over a continuous sample space, with $P^*$ a discrete distribution (e.g., the empirical distribution based on a sample), Theorem 7.3(ii) may fail, and we need to apply (i). For this we must be sure to take as the Bayes act $p(\cdot)$ against $P \in \mathcal{E}$ the specific choice where (74) holds everywhere (rather than almost everywhere). Then Corollary 7.3 holds.

See Section 7.6.1 for a case where neither (i) nor (ii) of Theorem 7.3 applies, leading to failure of Corollary 7.3.

7.6. *Examples.* We shall now illustrate the above theory for the Brier score, the logarithmic score and the zero–one loss. In particular we analyze in detail the simple case having $\mathcal{X} = \{-1, 0, 1\}$ and $T \equiv X$. For each decision problem we (i) show how Theorems 7.1 and 7.2 can be used to find robust Bayes acts, (ii) give the corresponding maximum entropy distributions and (iii) exhibit the associated generalized exponential family and specific entropy function.

7.6.1. *Brier score.* Consider the Brier score for $\mathcal{X} = \{x_1, \ldots, x_N\}$. By (17) we may write this score as

$$S(x, Q) = 1 - 2q(x) + \sum_j q(j)^2.$$

To try to apply Theorem 7.1 we search for a linear distribution $P_\tau \in \Gamma_\tau$. That is, we must find $(\beta_j)$ such that, for all $x \in \mathcal{X}$,

(92) $$1 - 2p_\tau(x) + \sum_y p_\tau(y)^2 = \beta_0 + \sum_{j=1}^{k} \beta_j t_j(x).$$

Equivalently, we must find $(\alpha_j)$ such that, for all $x$,

(93) $$p_\tau(x) \equiv \alpha_0 + \sum_{j=1}^{k} \alpha_j t_j(x).$$

The mean-value constraints

$$\sum_x t_j(x) p_\tau(x) = \tau_j, \qquad j = 1, \ldots, k,$$

together with the normalization constraint

$$\sum_x p_\tau(x) = 1,$$

will typically determine a unique solution for the $k+1$ coefficients $(\alpha_j)$ in (93). As long as this procedure leads to a nonnegative value for each $p_\tau(x)$, by Theorem 7.1 and the fact that the Brier score is proper we shall then have obtained a saddle-point $(P_\tau, P_\tau)$.

However, as we shall see below, for certain values of $\tau$ this putative "solution" for $P_\tau$ might have some $p_\tau(x)$ negative—showing that it is simply not possible to satisfy (92). By Theorem 5.2 we know that, even in this case a saddle-point $(P_\tau, P_\tau)$ exists. We can find it by applying Theorem 7.2: we first restrict the sample space to some $\mathcal{X}^* \subseteq \mathcal{X}$ and try to find a probability distribution $P_\tau$ satisfying the mean-value and normalization constraints, such that $p_\tau(x) = 0$ for $x \notin \mathcal{X}^*$ and for which, for some $(\beta_j)$ (92) holds for all $x \in \mathcal{X}^*$ [or, equivalently, for some $(\alpha_j)$ (93) holds for all $x \in \mathcal{X}^*$]. Among all such restrictions $\mathcal{X}^*$ that lead to an everywhere nonnegative solution for $(p_\tau(x))$, we choose that yielding the largest value of $H$. Then the resulting distribution $P_\tau$ will supply a saddle-point and so, simultaneously, (i) will have $H(P_\tau) = h(\tau)$, the maximum possible generalized entropy $1 - \sum_x p(x)^2$ subject to the mean-value constraints, and (ii) (which we regard as more important) will be robust Bayes for the Brier score against all distributions satisfying that constraint.

A more intuitive and more efficient geometric variant of the above procedure will be given in Section 8.

EXAMPLE 7.1. Suppose $\mathcal{X} = \{-1, 0, 1\}$ and $T \equiv X$. Consider the constraint $\mathrm{E}(X) = \tau$, for $\tau \in [-1, 1]$. We first look for linear acts satisfying (93). The mean-value constraint $\sum_x x\, p_\tau(x) = \tau$ and normalization constraint $\sum_x p_\tau(x) = 1$ provide two independent linear equations for the coefficients $(\alpha_0, \alpha_1)$ in (93), so uniquely determining $(\alpha_0, \alpha_1)$, and hence $p_\tau$. We easily find $\alpha_0 = \frac{1}{3}$, $\alpha_1 = \frac{1}{2}\tau$ and thus $p_\tau(x) = \frac{1}{3} + \frac{1}{2}\tau x$ $(x = -1, 0, 1)$ (whence $\beta_1 = -\tau$, $\beta_0 = \frac{2}{3} + \frac{1}{3}\tau^2$). We thus obtain a nonnegative solution for $(p_\tau(-1), p_\tau(0), p_\tau(1))$ only so long as $\tau \in [-2/3, 2/3]$: in this and only this case the act $p_\tau$ is linear. When $\tau$ falls outside this interval we can proceed by trying the restricted sample spaces $\{-1\}$, $\{0\}$, $\{1\}$, $\{0, 1\}$, $\{-1, 0\}$, $\{-1, 1\}$, as indicated above. All in all, we find that the optimal distribution $P_\tau$ has probabilities, entropy and $\beta$ satisfying Definition 7.3, as given in Table 1.

The family $\{P_\tau : -1 \leq \tau \leq 1\}$ constitutes the regular generalized exponential family over $\mathcal{X}$ generated by the Brier score and the statistic $T \equiv X$. The location of this family in the probability simplex is depicted in Figure 2.

We note that $h(\tau) = \beta_0 + \beta_1 \tau$ and $\beta_1 = h'(\tau)$ $(-1 < \tau < 1)$. The function $h(\tau)$ is plotted in Figure 3; Figure 4 shows the correspondence between $\beta_1$ and $\tau$.

By Theorem 7.1(i), the robust Bayes act $P_\tau$ will be an equalizer rule when $\tau$ is linear, that is, for $\tau \in [-\frac{2}{3}, \frac{2}{3}]$, and also (trivially) when $\tau = \pm 1$.

The above example demonstrates the need for condition (i) or (ii) in Theorem 7.3 and Corollary 7.3: typically both these conditions fail here for $\tau \notin [-\frac{2}{3}, \frac{2}{3}]$. Thus let $P^*$ have probabilities $(p^*(-1), p^*(0), p^*(1)) = (0.9, 0, 0.1)$, so that $\tau^* = \mathrm{E}_{P^*}(X) = -0.8$ and $\zeta_{\tau^*} = (0.8, 0.2, 0)$. From (18) we find $s^*(\zeta_{\tau^*}) = -0.24$. However, $\zeta_{\tau^*} = \zeta_{-0.8}$ is not the maximum support act in $\mathcal{F}^r$ in this case: it can be checked that this is given by $\zeta_{-0.95} = (0.95, 0.05, 0)$, having support $s^*(\zeta_\tau) = -0.195$.

7.6.2. *Log loss.* We now specialize the analysis of Section 7.3 to the case $\mathcal{X} = \{-1, 0, 1\}$, $T \equiv X$, with $\mu$ counting measure.

For $\tau \in (-1, 1)$, the maximum entropy distribution $P_\tau$ will have (robust Bayes) probability mass function of the form $p_\tau(x) = \exp -(\beta_0 + \beta_1 x)$. That is, the probability vector $p_\tau = (p_\tau(-1), p_\tau(0), p_\tau(1))$ will be of the form $(pe^{\beta_1}, p, pe^{-\beta_1})$, subject to the normalization and mean-value constraints

$$p(1 + e^{\beta_1} + e^{-\beta_1}) = 1, \tag{94}$$

$$p(e^{-\beta_1} - e^{\beta_1}) = \tau, \tag{95}$$

which uniquely determine $p \in (0, 1)$, $\beta_1 \in \mathcal{R}$. Then $h(\tau) = \beta_0 + \beta_1 \tau$, where $\beta_0 = -\log p$.

We thus have

$$p = (1 + e^{\beta_1} + e^{-\beta_1})^{-1}, \tag{96}$$

$$\tau = p(e^{-\beta_1} - e^{\beta_1}), \tag{97}$$

$$h = -\log p + \beta_1 \tau. \tag{98}$$

Table 1
*Brier score: maximum entropy distributions*

| | $p_\tau(-1)$ | $p_\tau(0)$ | $p_\tau(1)$ | $h(\tau)$ | $\beta_0$ | $\beta_1$ |
|---|---|---|---|---|---|---|
| $\tau = -1$ | 1 | 0 | 0 | 0 | $= \beta_1$ | $\beta_1 \geq 2$ |
| $-1 < \tau \leq -\frac{2}{3}$ | $-\tau$ | $1 + \tau$ | 0 | $-2\tau(1+\tau)$ | $2\tau^2$ | $-2 - 4\tau$ |
| $-\frac{2}{3} < \tau < \frac{2}{3}$ | $\frac{1}{3} - \frac{1}{2}\tau$ | $\frac{1}{3}$ | $\frac{1}{2}\tau + \frac{1}{3}$ | $\frac{2}{3} - \frac{1}{2}\tau^2$ | $\frac{2}{3} + \frac{1}{2}\tau^2$ | $-\tau$ |
| $\frac{2}{3} \leq \tau < 1$ | 0 | $1 - \tau$ | $\tau$ | $2\tau(1-\tau)$ | $2\tau^2$ | $2 - 4\tau$ |
| $\tau = 1$ | 0 | 0 | 1 | 0 | $= -\beta_1$ | $\beta_1 \leq -2$ |

40 P. D. GRÜNWALD AND A. P. DAWID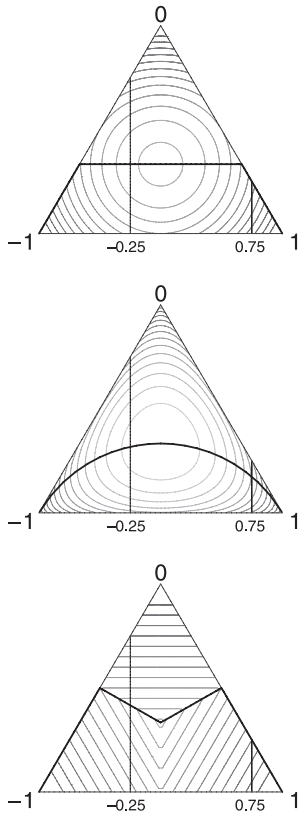

FIG. 2. *Brier score, logarithmic score and zero–one loss: the probability simplex for $\mathcal{X} = \{-1, 0, 1\}$, with entropy contours and generalized exponential family (maximum entropy distributions for the constraint $\mathrm{E}(X) = \tau$, $\tau \in [-1, 1]$). The set of distributions satisfying $\mathrm{E}(X) = \tau$ corresponds to a vertical line intersecting the base at $\tau$; this is displayed for $\tau = -0.25$ and $\tau = 0.75$. The intersection of the bold curve and the vertical line corresponding to $\tau$ represents the maximum entropy distribution for constraint $\mathrm{E}(X) = \tau$.*

On varying $\beta_1$ in $(-\infty, \infty)$, we obtain the parametric curve $(\tau, h)$ displayed in Figure 3; Figure 4 displays the correspondence between $\beta_1$ and $\tau$. It is readily verified that $dh/d\tau = (dh/d\beta_1)/(d\tau/d\beta_1) = \beta_1$, in accordance with (57).

In the terminology of Section 7.4, the above family $\{P_\tau : \tau \in (0,1)\}$ constitutes the natural exponential family associated with the logarithmic score and the statistic $T$. It is also the usual exponential family for this problem. However, the full exponential family further includes $\tau = \pm 1$. The family $\Gamma_1$ consists of the single distribution $P_1$ putting all its mass on the point 1. Then trivially $P_1$ is maximum entropy [with specific entropy $h(1) = 0$], and $p_1 = (0, 0, 1)$, with loss vector $L(\cdot, p_1) = (\infty, \infty, 0)$, is unique Bayes against $P_1$ and robust Bayes against $\Gamma_1$. Clearly (59) fails in this case, but even though $\tau = 1$ is not regular the property of Lemma 7.2 does hold there (albeit trivially). Similar properties apply at $\tau = -1$.

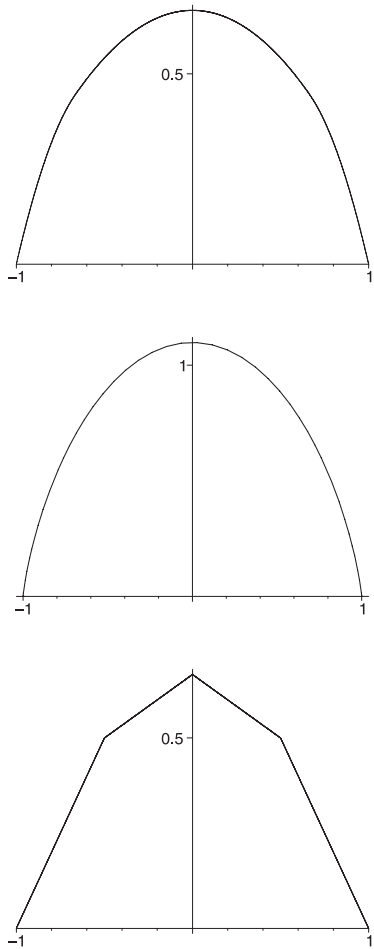

Fig. 3. *Specific entropy function $h(\tau)$ for Brier score, logarithmic score and zero–one loss.*

7.6.3. *Zero–one loss.* We now consider the zero–one loss (22) and seek robust Bayes acts against mean-value constraints $\Gamma_\tau$ of form (76). Once again we can try to apply Theorem 7.1 by looking for an act $\zeta_\tau \in \mathcal{Z}$ that is Bayes against some $P_\tau \in \Gamma_\tau$, and such that

(99) $$L(x, \zeta_\tau) \equiv 1 - \zeta_\tau(x) = \beta_0 + \sum_{j=1}^{k} \beta_j \, t_j(x)$$

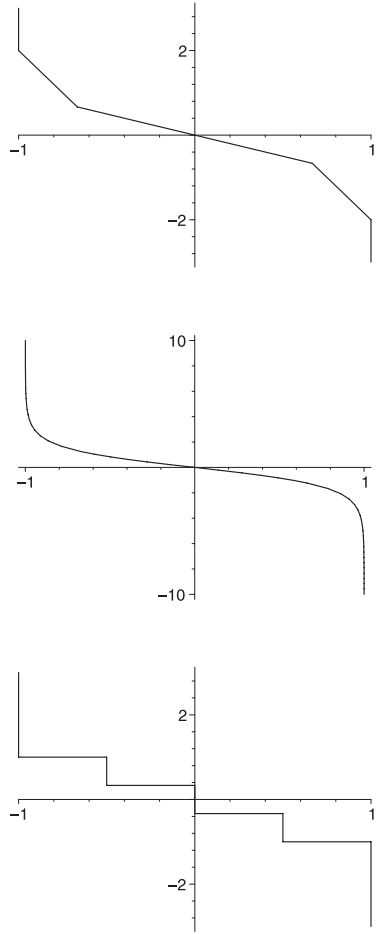

Fig. 4. *Correspondence between mean-value parameter $\tau$ (x-axis) and natural parameter $\beta_1$ (y-axis) of generalized exponential family, for Brier score, logarithmic score and zero–one loss.*

for all $x \in \mathcal{X}$. When this proves impossible, we can again proceed by restricting the sample space and using Theorem 7.2. The distribution $P_\tau$ will again maximize the generalized entropy. However, in this problem, in contrast to the log and Brier score cases, because of nonsemistrictness the Bayes act against $P_\tau$ may be nonunique—and, if we want to ensure that (99) (or its restricted version) holds, it may matter which of the Bayes acts (including randomized acts) we pick. Thus the familiar routine "maximize the generalized entropy, and then use a Bayes act against this distribution" is not, by itself, fully adequate to derive the robust Bayes act: additional care must be taken to select the *right* Bayes act.

EXAMPLE 7.2. Again take $\mathcal{X} = \{-1, 0, 1\}$ and $T \equiv X$. Consider the constraint $\mathrm{E}(X) = \tau$, where $\tau \in [-1, 1]$. We find that for each $\tau$ a unique maximum entropy $P_\tau$ exists. By some algebra we can then find the probabilities $(p_\tau(-1), p_\tau(0), p_\tau(1))$; they are given in Table 2, together with the corresponding specific entropy $h(\tau)$ (also plotted in Figure 3).

The family of distributions $\{P_\tau : \tau \in [-1, 1]\}$ thus constitutes the full generalized exponential family over $\mathcal{X}$ generated by the zero–one loss and the statistic $T \equiv X$. The location of this family in the probability simplex is depicted in Figure 2.

How can we determine the robust Bayes acts $\zeta_\tau$? We know that any such $\zeta_\tau$ is Bayes against $P_\tau$ and thus puts all its mass on the modes of $P_\tau$. As can be seen, for $-0.5 \leq \tau \leq 0.5$ the set $\mathcal{A}_{P_\tau}$ of these modes has more than one element. We additionally use (99), restricted to $x$ in the support of $P_\tau$, to find out which $\zeta_\tau \in \mathcal{A}_{P_\tau}$ are robust Bayes. For $\tau \in [-\frac{1}{2}, \frac{1}{2}]$ this requires

$$
\begin{align}
-\beta_1 + \beta_0 &= 1 - \zeta_\tau(-1), \\
\beta_0 &= 1 - \zeta_\tau(0), \\
\beta_1 + \beta_0 &= 1 - \zeta_\tau(1),
\end{align}
\tag{100}
$$

from which we readily deduce $\beta_0 = \frac{2}{3}$. The condition that $\zeta_\tau$ put all its mass on the modes of $P_\tau$ then uniquely determines $\zeta_\tau$ for $-0.5 \leq \tau < 0$ and for $0 < \tau < 0.5$. If $\tau = 0$, all acts $\zeta$ are Bayes for some $P \in \Gamma_\tau$ (take $P$ uniform), and hence by Theorem 7.1 all solutions to (100) [i.e., such that $\zeta_\tau(0) = \frac{1}{3}$] are robust Bayes acts. Finally, for $\tau = 0.5$ (the case $\tau = -0.5$ is similar) we must have $\zeta_\tau(-1) = 0$, and we can use the "supporting hyperplane" property (56) to deduce that $\zeta_\tau(0) \leq \frac{1}{3}$.

Table 3 gives the robust Bayes acts $\zeta_\tau$ for each $\tau \in [-1, 1]$, together with the corresponding values of $\beta_0, \beta_1$. Thus $\zeta_\tau$ is a linear act for $-0.5 \leq \tau \leq 0.5$ (where we must choose $a = \frac{1}{3}$ at the endpoints). Again we see that $h(\tau) = \beta_0 + \beta_1 \tau$, and that $\beta_1 = h'(\tau)$ where this exists.

TABLE 2
*Zero–one loss: maximum entropy distributions*

|  | $p_\tau(-1)$ | $p_\tau(0)$ | $p_\tau(1)$ | $h(\tau)$ |
|---|---|---|---|---|
| $\tau = -1$ | 1 | 0 | 0 | 0 |
| $-1 < \tau < -\frac{1}{2}$ | $-\tau$ | $1 + \tau$ | 0 | $1 + \tau$ |
| $\tau = -\frac{1}{2}$ | $\frac{1}{2}$ | $\frac{1}{2}$ | 0 | $\frac{1}{2}$ |
| $-\frac{1}{2} < \tau < 0$ | $\frac{1-\tau}{3}$ | $\frac{1-\tau}{3}$ | $\frac{1+2\tau}{3}$ | $\frac{2+\tau}{3}$ |
| $\tau = 0$ | $\frac{1}{3}$ | $\frac{1}{3}$ | $\frac{1}{3}$ | $\frac{2}{3}$ |
| $0 < \tau < \frac{1}{2}$ | $\frac{1-2\tau}{3}$ | $\frac{1+\tau}{3}$ | $\frac{1+\tau}{3}$ | $\frac{2-\tau}{3}$ |
| $\tau = \frac{1}{2}$ | 0 | $\frac{1}{2}$ | $\frac{1}{2}$ | $\frac{1}{2}$ |
| $\frac{1}{2} < \tau < 1$ | 0 | $1 - \tau$ | $\tau$ | $1 - \tau$ |
| $\tau = 1$ | 0 | 0 | 1 | 0 |

Figure 4 shows the relationship between $\beta_1$ and $\tau$. In this case the uniqueness part of Condition 7.3 is not satisfied, with the consequence that neither $\beta_1$ nor $\tau$ uniquely determines the other. However, the full exponential family $\{P_\tau : -1 \leq \tau \leq 1\}$ is clearly specified by the one-one map $\tau \mapsto P_\tau$, and most of the properties of such families remain valid.

**8. Relative entropy, discrepancy, divergence.** Analogous to our generalized definition of entropy, we now introduce *generalized relative entropy* with respect to a decision problem, and we show how the negative relative entropy has a natural interpretation as a measure of discrepancy. This allows us to extend our minimax results to a more general setting and leads to a generalization of the Pythagorean property of the relative Shannon entropy [Csiszár (1975)].

We first introduce the concept of the *discrepancy* between a distribution $P$ and a (possibly randomized) act $\zeta$, induced by a decision problem.

8.1. *Discrepancy.* Suppose first that $H(P)$ is finite. We define, for any $\zeta \in \mathcal{Z}$, the *discrepancy* $D(P, \zeta)$ between the distribution $P$ and the act $\zeta$ by

$$(101) \qquad D(P, \zeta) := L(P, \zeta) - H(P).$$

In the general terminology of decision theory, $D(P, \zeta)$ measures DM's *regret* [Berger (1985), Section 5.5.5] associated with taking action $\zeta$ when Nature generates $X$ from $P$. Also, since $-D(P, \zeta)$ differs from $-L(P, \zeta)$ by a term only involving $P$, we can use it in place of the support function $s_P(\zeta)$: thus maximizing support is equivalent to minimizing discrepancy.

We note that, if a Bayes act $\zeta_P$ against $P$ exists, then

$$(102) \qquad D(P, \zeta) = E_P\{L(X, \zeta) - L(X, \zeta_P)\}.$$

Table 3
*Zero–one loss: robust Bayes acts*

|  | $\zeta_\tau(-1)$ | $\zeta_\tau(0)$ | $\zeta_\tau(1)$ | $\beta_0$ | $\beta_1$ |
|---|---|---|---|---|---|
| $\tau = -1$ | 1 | 0 | 0 | $= \beta_1$ | $\beta_1 \geq 1$ |
| $-1 < \tau < -\frac{1}{2}$ | 1 | 0 | 0 | 1 | 1 |
| $\tau = -\frac{1}{2}$ | $1-a$ | $a \leq \frac{1}{3}$ | 0 | $1-a$ | $1-2a$ |
| $-\frac{1}{2} < \tau < 0$ | $\frac{2}{3}$ | $\frac{1}{3}$ | 0 | $\frac{2}{3}$ | $\frac{1}{3}$ |
| $\tau = 0$ | $a \leq \frac{2}{3}$ | $\frac{1}{3}$ | $\frac{2}{3} - a$ | $\frac{2}{3}$ | $a - \frac{1}{3}$ |
| $0 < \tau < \frac{1}{2}$ | 0 | $\frac{1}{3}$ | $\frac{2}{3}$ | $\frac{2}{3}$ | $-\frac{1}{3}$ |
| $\tau = \frac{1}{2}$ | 0 | $a \leq \frac{1}{3}$ | $1-a$ | $1-a$ | $2a - 1$ |
| $\frac{1}{2} < \tau < 1$ | 0 | 0 | 1 | 1 | $-1$ |
| $\tau = 1$ | 0 | 0 | 1 | $= -\beta_1$ | $\beta_1 \leq -1$ |

We shall also use (102) as the *definition* of $D(P, \zeta)$ when $P \notin \mathcal{P}$, or $H(P)$ is not finite, but $P$ has a Bayes act (in which case it will not matter which such Bayes act we choose). This definition can itself be generalized further to take account of some cases where no Bayes act exists; we omit the details.

The function $D$ has the following properties:

(i) $D(P, \zeta) \in [0, \infty]$.
(ii) $D(P, \zeta) = 0$ if and only if $\zeta$ is Bayes against $P$.
(iii) For any $a, a' \in \mathcal{A}$, $D(P, a) - D(P, a')$ is linear in $P$ (in the sense of Lemma 3.2).
(iv) $D$ is a convex function of $P$.

Conversely, under regularity conditions any function $D$ satisfying (i)–(iii) above can be generated from a suitable decision problem by means of (101) or (102) [Dawid (1998)].

8.1.1. *Discrepancy and divergence.* When our loss function is a $\mathcal{Q}$-proper scoring rule $S$, we shall typically denote the corresponding discrepancy function by $d$. Thus for $P, Q \in \mathcal{Q}$ with $H(P)$ finite,

$$(103) \qquad d(P, Q) = S(P, Q) - H(P).$$

We now have $d(P, Q) \geq 0$, with equality when $Q = P$; if $S$ is $\mathcal{Q}$-strict, then $d(P, Q) > 0$ for $Q \neq P$. Conversely, if for any scoring rule $S$, $S(P, Q) - S(P, P)$ is nonnegative for all $P, Q \in \mathcal{Q}$, then the scoring rule $S$ is $\mathcal{Q}$-proper. We refer to $d(P, Q)$ as the *divergence* between the distributions $P$ and $Q$. As we shall see in Section 10, divergence can be regarded as analogous to a measure of squared Euclidean distance.

The following lemma, generalizing Lemmas 4 and 7 of Topsøe (1979), follows easily from (103) and the linearity of $S(P, Q)$ in $P$.

LEMMA 8.1. *Let $S$ be a proper scoring rule, with associated entropy function $H$ and divergence function $d$. Let $P_1, \ldots, P_n$ have finite entropies, and let $(p_1, \ldots, p_n)$ be a probability vector. Then*

$$(104) \qquad H(\overline{P}) = \sum p_i H(P_i) + \sum p_i d(P_i, \overline{P}),$$

$$(105) \qquad d(\overline{P}, Q) = \sum p_i d(P_i, Q) - \sum p_i d(P_i, \overline{P}),$$

*where $\overline{P} := \sum p_i P_i$.*

We can also associate a divergence with a more general decision problem, with loss function $L$ such that $\mathcal{Z}_Q$ is nonempty for all $Q \in \mathcal{Q}$, by

$$(106) \qquad d(P, Q) := D(P, \zeta_Q) = \mathrm{E}_P\{L(X, \zeta_Q) - L(X, \zeta_P)\},$$

where again for each $Q \in \mathcal{Q}$ we suppose we have selected some specific Bayes act $\zeta_Q$. This will then be identical with the divergence associated directly [using, e.g., (103)] with the corresponding scoring rule given by (15), and (104) and (105) will continue to hold with this more general definition.

8.2. *Relative loss.* Given a game $\mathcal{G} = (\mathcal{X}, \mathcal{A}, L)$, choose, once and for all, a *reference act* $\zeta_0 \in \mathcal{Z}$. We can construct a new game $\mathcal{G}_0 = (\mathcal{X}, \mathcal{A}, L_0)$, where the new loss function $L_0$ is given by

$$(107) \qquad L_0(x, a) := L(x, a) - L(x, \zeta_0).$$

This extends naturally to randomized acts: $L_0(x, \zeta) := L(x, \zeta) - L(x, \zeta_0)$. We call $L_0$ the *relative loss function* and $\mathcal{G}_0$ the *relative game* with respect to the reference act $\zeta_0$. In order that $L_0 > -\infty$ we shall require $L(x, \zeta_0) < \infty$ for all $x \in \mathcal{X}$. We further restrict attention to distributions in $\mathcal{P}' := \{P : L_0(P, a) \text{ is defined for all } a \in \mathcal{A}\}$ and randomized acts in $\mathcal{Z}' := \{\zeta : L_0(P, \zeta) \text{ is defined for all } P \in \mathcal{P}'\}$. In general, $\mathcal{P}'$ and $\mathcal{Z}'$ may not be identical with $\mathcal{P}$ and $\mathcal{Z}$.

The expected relative loss $L_0(P, \zeta)$ satisfies

$$(108) \qquad L_0(P, \zeta) = L(P, \zeta) - L(P, \zeta_0)$$

whenever $L(P, \zeta_0)$ is finite. Whether or not this is so, it is easily seen that the Bayes acts against any $P$ are the same in both games.

DEFINITION 8.1. An act $\zeta_0 \in \mathcal{Z}$ is called *neutral* if the loss function $L(x, \zeta_0)$ is a finite constant, $k$ say, on $\mathcal{X}$.

If a neutral act exists, and we use it as our reference act, then $L_0(P, \zeta) = L(P, \zeta) - k$, all $P \in \mathcal{P}$. The relative game $\mathcal{G}_0$ is then effectively the same as the original game $\mathcal{G}$, and maximum entropy distributions, saddle-points, and other properties of the two games, or of their restricted subgames, will coincide. However, these equivalences are typically not valid for more general relative games.

8.3. *Relative entropy.* When a Bayes act $\zeta_P$ against $P$ exists, the *generalized relative entropy* $H_0(P) := \inf_{a \in \mathcal{A}} L_0(P, a)$ associated with the relative loss $L_0$ is seen to be

$$(109) \qquad H_0(P) = \mathrm{E}_P\{L(X, \zeta_P) - L(X, \zeta_0)\}.$$

[In particular, we must have $-\infty \leq H_0(P) \leq 0$.] When $L(P, \zeta_0)$ is finite, this becomes

$$(110) \qquad H_0(P) = H(P) - L(P, \zeta_0).$$

Comparing (109) with (102), we observe the simple but fundamental relation

$$(111) \qquad H_0(P) = -D(P, \zeta_0).$$

The *maximum generalized relative entropy criterion* thus becomes identical to the *minimum discrepancy criterion*:

*Choose $P \in \Gamma$ to minimize, over $P \in \Gamma$, its discrepancy $D(P, \zeta_0)$ from the reference act $\zeta_0$.*

Note that, even though Bayes acts are unaffected by changing from $L$ to the relative loss $L_0$, the corresponding entropy function (110) is *not* unaffected. Thus in general the maximum entropy criterion (for the same constraints) will deliver different solutions in the two problems. Related to this, we can also expect to obtain different robust Bayes acts in the two problems.

Suppose we construct the relative loss taking as our reference act $\zeta_0$ a Bayes act against a fixed *reference distribution* $P_0$. Alternatively, start with a proper scoring rule $S$, and construct directly the relative score with reference to the act $P_0$. The minimum discrepancy criterion then becomes the *minimum divergence criterion*: choose $P \in \Gamma$ to minimize the divergence $d(P, P_0)$ from the reference distribution $P_0$.

This reinterpretation can often assist in finding a maximum relative entropy distribution. If moreover we can choose $P_0$ to be neutral, this minimum divergence criterion becomes equivalent to maximizing entropy in the original game.

### 8.4. *Relative loss and generalized exponential families.*

#### 8.4.1. *Invariance relative to linear acts.*

Suppose the reference act $\zeta_0$ is linear with respect to $L$ and $T$, so that we can write

$$(112) \qquad L(x, \zeta_0) = \delta_0 + \delta^{\mathrm{T}} t(x).$$

Then if $\mathrm{E}_P(T)$ exists,

$$(113) \qquad L_0(P, \zeta) = L(P, \zeta) - \delta_0 - \delta^{\mathrm{T}} \mathrm{E}_P(T),$$

$$(114) \qquad H_0(P) = H(P) - \delta_0 - \delta^{\mathrm{T}} \mathrm{E}_P(T).$$

In particular, for all $P \in \Gamma_\tau$,

$$(115) \qquad L_0(P, \zeta) = L(P, \zeta) - \delta_0 - \delta^{\mathrm{T}} \tau,$$

$$(116) \qquad H_0(P) = H(P) - \delta_0 - \delta^{\mathrm{T}} \tau.$$

We see immediately from the definitions that the full, the natural, the regular and the linear generalized exponential families generated by $L_0$ and $T$ are identical with those generated by $L$ and $T$. The correspondence $\tau \mapsto P_\tau$ is unaffected; for the natural case, if $Q_\beta$ arises from $L$ and $Q_{0,\beta}$ from $L_0$, we have $Q_{0,\beta} = Q_{\beta+\delta}$.

Suppose in particular that we take any $P_\sigma \in \mathcal{E}^l$. In this case we can take $\zeta_0$ having property (112) to be the corresponding Bayes act $\zeta_\sigma$. We thus see that a generalized exponential family is unchanged when the loss function is redefined by taking it relative to some linear member of the family. This property is well known for the case of a standard exponential family, where every regular member is linear (with respect to the logarithmic score). In that case, the relative loss can also be interpreted as the logarithmic score when the base measure $\mu$ is changed to $P_\sigma$; the exponential family is unchanged by such a choice.

8.4.2. *Lafferty additive models.* Lafferty (1999) defines *the additive model relative to a Bregman divergence $d$, reference measure $P_0$ and constraint random variable $T:\mathcal{X} \to \mathcal{R}$* as the family of probability measures $\{Q_\beta : \beta \in \mathcal{R}\}$ where

$$(117) \qquad Q_\beta := \arg\min_{P \in \mathcal{P}} \beta \mathrm{E}_P\{T(X)\} + d(P, P_0).$$

We note that $P_0 = Q_0$ is in this family.

Let $S$ be the Bregman score (29) associated with $d$ and let $S_0$ be the associated relative score $S_0(x, Q) \equiv S(x, Q) - S(x, P_0)$. Note that by (111) $d(P, P_0) = -H_0(P)$, where $H_0(P)$ is the entropy associated with $S_0$. Lafferty's additive models are thus special cases of our natural generalized exponential families as defined in Section 7.4, being generated by the specific loss function $S_0$ and statistic $T$. As shown in Section 8.4.1, when $P_0$ is linear (with respect to $S$ and $T$) the previous sentence remains true on replacing $S_0$ by $S$.

These considerations do not rely on any special Bregman properties, and so extend directly to any loss-based divergence function $d$ of the form given by (103) or (106).

8.5. *Examples.*

8.5.1. *Brier score.* In the case of the Brier score, the divergence between $P$ and $Q$ is given by the squared Euclidean distance between their probability vectors:

$$(118) \qquad d(P, Q) = \|p - q\|^2 = \sum_j \{p(j) - q(j)\}^2.$$

Using a reference distribution $P_0$, the relative entropy thus becomes

$$(119) \qquad H_0(P) = -\sum_j \{p(j) - p_0(j)\}^2.$$

The uniform distribution over $\mathcal{X}$ is neutral. Therefore the distribution within a set $\Gamma$ that maximizes the Brier entropy is just that minimizing the discrepancy from the uniform reference distribution $P_0$.

To see the consequences of this for the construction of generalized Brier exponential families, let $\mathcal{X} = \{-1, 0, 1\}$ and consider the Brier score picture in Figure 2. The bold line depicts the maximum entropy distributions for constraints $\mathrm{E}(T) = \tau$, $\tau \in [-1, 1]$. By the preceding discussion, these coincide with the minimum $P_0$-discrepancy distributions. For each fixed value of $\tau$, the set $\Gamma_\tau = \{P : \mathrm{E}_P(X) = \tau\}$ is represented by the vertical line through the simplex intersecting the base line at the coordinate $\tau$. In Figure 2 the cases $\tau = -0.25$ and $\tau = 0.75$ are shown explicitly. The minimum discrepancy distribution within $\Gamma_\tau$ will be given by the point on that line within the simplex that is nearest to the center of the simplex. This gives us a simple geometric means to find the minimum relative discrepancy distributions for $\tau \in [-1, 1]$, involving less work than the procedure detailed in Section 7.6.1. We easily see that for $\tau \in [-2/3, 2/3]$ the minimizing point $p_\tau$ is in the interior of the line segment, while for $\tau$ outside this interval the minimizing point is at one end of the segment.

8.5.2. *Logarithmic score.* For $P \in \mathcal{M}$ (i.e., $P \ll \mu$) any version $p$ of the density $dP/d\mu$ is Bayes against $P$. Then, with $q$ any version of $dQ/d\mu$, $D(P, q) = \mathrm{E}_P[\log\{p(X)/q(X)\}]$ is the Kullback–Leibler divergence $\mathrm{KL}(P, Q)$ and does not depend on the choice of the versions of either $p$ or $q$. Again, for $P, Q \in \mathcal{M}$ we can treat $S$ as a proper scoring rule $S(x, Q)$, with $d(P, Q) \equiv \mathrm{KL}(P, Q)$ as its associated divergence. [For $P \notin \mathcal{M}$ there is no Bayes act (see Section 3.5.2), and so, according to our definition (102), the discrepancy $D(P, q)$ is not defined: we might define it as $+\infty$ in this case.] Maximizing the relative entropy is thus equivalent to minimizing the Kullback–Leibler divergence in this case.

There is a simple relationship between the choice of base measure $\mu$, which is a necessary input to our specification of the decision problem, and the use of a reference distribution for defining relative loss. If we had constructed our logarithmic loss using densities starting with a different choice $\mu_0$ of base measure, where $\mu_0$ is mutually absolutely continuous with $\mu$, we should have obtained instead the loss function $S_0(x, Q) = -\log q_0(x)$, with $q_0(x) = (dQ/d\mu_0)(x) = (dQ/d\mu)(x) \times (d\mu/d\mu_0)(x)$. Thus $S_0(x, Q) = S(x, Q) + k(x)$, with $k(x) \equiv -\log d(x)$, where $d$ is some version of $d\mu/d\mu_0$. In particular, when $\mu_0$ is a probability measure, this is exactly the relative loss function (107) with respect to the reference distribution $\mu_0$, when we start from the problem constructed in terms of $\mu$ (in particular, it turns out that this relative game will not depend on the starting measure $\mu$). As already determined, the corresponding relative entropy function is $H_0(P) = -\mathrm{KL}(P, \mu_0)$.

8.5.3. *Zero–one loss.* In this case, the discrepancy between $P$ and an act $\zeta \in \mathcal{Z}$ is given by

$$D(P, \zeta) = p_{\max} - \sum_{j \in \mathcal{X}} p(j)\zeta(j). \tag{120}$$

When $\mathcal{X}$ has finite cardinality $N$, and $\zeta_0$ is the randomized act that chooses uniformly from $\mathcal{X}$, we have $S(x, \zeta_0) \equiv 1 - 1/N$, so that this choice of $\zeta_0$ is neutral.

Take $\mathcal{X} = \{-1, 0, 1\}$ and $T \equiv X$, let $\zeta_0$ be uniform on $\mathcal{X}$ and consider the minimum zero–one $\zeta_0$-discrepancy distributions shown in Figure 2. Determining this family of distributions geometrically is easy once one has determined the contours of constant generalized entropy, since these are also the contours of constant discrepancy from $\zeta_0$.

8.5.4. *Bregman divergence.* In a finite sample space, the Bregman score (29) generates the Bregman divergence (30). Thus minimizing the Bregman divergence is equivalent to maximizing the associated relative entropy, which is in turn equivalent to finding a distribution that is robust Bayes against the associated relative loss function. Minimizing a Bregman divergence has become a popular tool in the construction and analysis of on-line learning algorithms [Lafferty (1999) and Azoury and Warmuth (2001)], on account of numerous pleasant properties it enjoys. As shown by properties (i)–(iv) of Section 8.1 and as will further be seen in Section 10, many of these properties generalize to an arbitrary decision-based divergence function as defined by (103) or (106).

In more general sample spaces, the separable Bregman score (34) generates the separable Bregman divergence $d_\psi$ given by (37). When the measure $\mu$ appearing in these formulae is itself a probability distribution, $\mu$ will be neutral (uniquely so if $\psi$ is strictly convex); then minimizing over $P$ the separable Bregman divergence $d_\psi(P, \mu)$ of $P$ from $\mu$ becomes equivalent to maximizing the separable Bregman entropy $H(P)$ as given by (38).

**9. Statistical problems: discrepancy as loss.** In this section we apply the general ideas presented so far to more specifically statistical problems.

9.1. *Parametric prediction problems.* In a statistical decision problem, we have a family $\{P_\omega : \omega \in \Omega\}$ of distributions for an observable $X$ over $\mathcal{X}$, labelled by the values $\omega$ of a parameter $\Omega$ ranging over $\Omega$; the consequence of taking an action $a$ depends on the value of $\Omega$. We shall show how one can construct a suitable loss function for this purpose, starting from a general decision problem $\mathcal{G}$ with loss depending on the value of $X$, and relate the minimax properties of the derived statistical game $\widehat{\mathcal{G}}$ to those of the underlying basic game $\mathcal{G}$.

In our context $X$ is best thought of as a future outcome to be predicted, perhaps after conducting a statistical experiment to learn about $\Omega$. The distributions of $X$ given $\Omega = \omega$ would often be taken to be the same as those governing the data in the experiment, but this is not essential. Our emphasis is thus on statistical models for prediction, rather than for observed data: the latter will not enter directly. For applications of this predictive approach to problems of experimental design, see Dawid (1998) and Dawid and Sebastiani (1999).

9.2. *Technical framework.* Let $(\mathcal{X}, \mathcal{B})$ be a separable metric space with its Borel $\sigma$-field, and let $\mathcal{P}_0$ be the family of all probability distributions over $(\mathcal{X}, \mathcal{B})$. We shall henceforth want to consider $\mathcal{P}_0$ itself (and subsets thereof) as an abstract "parameter space." When we wish to emphasize this point of view we shall denote $\mathcal{P}_0$ by $\Theta_0$, and its typical member by $\theta$; when $\theta$ is considered in its original incarnation as a probability distribution on $(\mathcal{X}, \mathcal{B})$, we may also denote it by $P_\theta$.

$\Theta_0$ becomes a metric space under the Prohorov metric in $\mathcal{P}_0$, and the associated topology is then identical with the weak topology on $\mathcal{P}_0$ [Billingsley (1999), page 72]. We denote the set of all probability distributions, or *laws*, on the Borel $\sigma$-field $\mathcal{C}$ in $\Theta_0$ by $\mathcal{L}_0$. Such a law can be regarded as a "prior distribution" for a parameter random variable $\Theta$ taking values in $\Theta_0$. For such a law $\Pi \in \mathcal{L}_0$, we denote by $P_\Pi \in \mathcal{P}_0$ its mean, given by $P_\Pi(A) = \mathrm{E}_\Pi\{P_\Theta(A)\}$ ($A \in \mathcal{B}$): this is just the marginal "predictive" (mixture) distribution for $X$ over $\mathcal{X}$, obtained by first generating a value $\theta$ for $\Theta$ from $\Pi$, and then generating $X$ from $P_\theta$.

9.3. *The derived game.* Starting from a basic game $\mathcal{G} = (\mathcal{X}, \mathcal{A}, L)$, we construct a new *derived game*, $\widehat{\mathcal{G}} := (\Theta, \mathcal{A}, \widehat{L})$. The new loss function $\widehat{L}$ on $\Theta \times \mathcal{A}$ is just the discrepancy function for the original game $\widehat{\mathcal{G}}$,

$$\widehat{L}(\theta, a) := D(P_\theta, a), \tag{121}$$

and the original sample space $\mathcal{X}$ is replaced by $\Theta := \{\theta \in \Theta_0 : D(P_\theta, a) \text{ is defined for all } a \in \mathcal{A}\}$.

We have

$$\widehat{L}(\theta, a) = L(P_\theta, a) - H(P_\theta) \tag{122}$$

when $H(P_\theta)$ is finite. Properties (121) and (122) then extend directly to randomized acts $\zeta \in \mathcal{Z}$ for DM. A randomized act for Nature in $\widehat{\mathcal{G}}$ is a law putting all its mass on $\Theta \subseteq \Theta_0$. We shall denote the set of such laws by $\mathcal{L} \subseteq \mathcal{L}_0$.

Note that $\widehat{L}(\theta, a)$ is just the regret associated with taking action $a$ when $X \sim P_\theta$. It is nonnegative, and it vanishes if and only if $a$ is Bayes against $P_\theta$. Such a regret function will often be a natural loss function to use in a statistical decision problem.

Since $\widehat{L} \geq 0$, the expected loss $\widehat{L}(\Pi, \zeta)$ is defined in $[0, \infty]$ for all $\Pi \in \mathcal{L}$, $\zeta \in \mathcal{Z}$. From (122) we obtain

$$\widehat{L}(\Pi, \zeta) = L(P_\Pi, \zeta) - \int H(P_\theta) \, d\Pi(\theta) \tag{123}$$

when the integral exists. An act $\zeta_0$ will thus be Bayes against $\Pi$ in $\widehat{\mathcal{G}}$ if and only if it is Bayes against $P_\Pi$ in $\mathcal{G}$. More generally, this equivalence follows from the property $\mathrm{E}_\Pi\{\widehat{L}(\Theta,\zeta) - \widehat{L}(\Theta,\zeta_0)\} = \mathrm{E}_{P_\Pi}\{L(X,\zeta) - L(X,\zeta_0)\}$. In particular, if $L$ is a $\mathcal{Q}$-proper scoring rule in the basic game $\mathcal{G}$, and the mixture distribution $P_\Pi \in \mathcal{Q}$, then $P_\Pi$ will be Bayes against $\Pi$ in $\widehat{\mathcal{G}}$.

The *derived entropy function* is

$$\widehat{H}(\Pi) = H(P_\Pi) - \int H(P_\theta)\, d\Pi(\theta) \tag{124}$$

(when the integral exists) and is nonnegative. This measures the expected reduction in uncertainty about $X$ obtainable by learning the value of $\Theta$, when initially $\Theta \sim \Pi$: it is the *expected value of information* [DeGroot (1962)] in $\Theta$ about $X$.

The derived discrepancy is just

$$\widehat{D}(\Pi,\zeta) = D(P_\Pi,\zeta). \tag{125}$$

9.4. *A statistical model.* Let $\Omega \subseteq \Theta_0$: for example, $\Omega$ might be a parametric family of distributions for $X$. We can think of $\Omega$ as the statistical model for the generation of $X$. We will typically write $\omega$ or $P_\omega$ for a member of $\Omega$ and use $\Omega$ to denote the parameter $\Theta$ when it is restricted to taking values in $\Omega$. We denote by $\Delta \subseteq \mathcal{L}_0$ the class of laws on $\Theta_0$ that give all their mass to $\Omega$ and can thus serve as priors for the parameter $\Omega$ of the model; we denote by $\Gamma \subseteq \mathcal{P}_0$ the family $\{P_\Pi : \Pi \in \Delta\}$ of all distributions for $X$ obtainable as mixtures over the model $\Omega$. Clearly both $\Delta$ and $\Gamma$ are convex.

LEMMA 9.1. *Suppose that the family $\Omega$ of distributions on $(\mathcal{X},\mathcal{B})$ is tight. Then so too are $\Gamma$ and $\Delta$ [the latter as a family of laws on $(\Theta_0,\mathcal{C})$].*

PROOF. The tightness of $\Gamma$ follows easily from the definition.

Let $\overline{\Omega}$ denote the closure of $\Omega$ in $\Theta_0$. Since $\Omega$ is tight, so is $\overline{\Omega}$ [use, e.g., Theorem 3.1.5(iii) of Stroock (1993)], and then Prohorov's theorem [Billingsley (1999), Theorem 5.1] implies that $\overline{\Omega}$ is compact in the weak topology. Any collection (in particular, $\Delta$) of distributions on $(\Theta_0,\mathcal{C})$ supported on $\overline{\Omega}$ is then necessarily tight. □

9.5. *Minimax properties.* Now consider a statistical model with $\Omega \subseteq \Theta$ (so that $\Delta \subseteq \mathcal{L}$). We can tailor the derived game $\widehat{\mathcal{G}}$ to this model by simply restricting the domain of $\widehat{L}$ to $\Omega \times \mathcal{A}$. We would thus be measuring the loss (regret) of taking act $\zeta \in \mathcal{Z}$, when the true parameter value is $\omega \in \Omega$, by $\widehat{L}(\omega,\zeta) = D(P_\omega,\zeta)$. Alternatively, and equivalently, we can focus attention on the *restricted game* $\widehat{\mathcal{G}}^\Delta$ as defined in Section 4.2, with $\Delta$ the family of laws supported on the model $\Omega$. In the present context we shall denote this by $\widehat{\mathcal{G}}^\Omega$.

We will often be interested in the existence and characterization of a value, saddle-point, maximum entropy (maximin) prior $\widehat{\Pi}^*$ or robust Bayes (minimax) act $\widehat{\zeta}^*$, in the game $\widehat{\mathcal{G}}^\Omega$. Note in particular that, when we do have a saddle-point $(\widehat{\Pi}^*, \widehat{\zeta}^*)$ in $\widehat{\mathcal{G}}^\Omega$, with value $\widehat{H}^*$, we can use Lemma 4.2 to deduce that $\widehat{\Pi}^*$ must put all its mass on $\Upsilon := \{\omega \in \Omega : D(P_\omega, \widehat{\zeta}^*) = \widehat{H}^*\}$: in particular, with $\widehat{\Pi}^*$-prior probability 1 the discrepancy from the minimax act is constant. When, as will typically hold, $\Upsilon$ is a proper subset of $\Omega$, we further deduce from Corollary 4.4 that $\widehat{\zeta}^*$ is not an equalizer rule in $\widehat{\mathcal{G}}^\Omega$.

To investigate further the minimax and related properties of the game $\widehat{\mathcal{G}}^\Omega$, we could try to verify directly for this game the requirements of the general theorems already proved in Sections 5–7. However, under suitable conditions these required properties will themselves follow from properties of the basic game $\mathcal{G}$. We now detail this relationship for the particular case of Theorem 6.4.

We shall impose the following condition:

CONDITION 9.1. There exists $K \in \mathcal{R}$ such that $H(P_\omega) \geq K$ for all $\omega \in \Omega$.

By concavity of $H$, Condition 9.1 is equivalent to $H(Q) \geq K$ for all $Q \in \Gamma$. The following lemma is proved in the Appendix.

LEMMA 9.2. *Suppose Condition 9.1 holds. Then if Conditions 6.1 and 6.3 hold for $L$ and $\Gamma$ (in $\mathcal{G}$), they likewise hold for $\widehat{L}$ and $\Delta$ (in $\widehat{\mathcal{G}}$).*

The next theorem now follows directly from Lemmas 9.1 and 9.2 and Theorem 6.4.

THEOREM 9.1. *Suppose Conditions 6.1, 6.3 and 9.1 all hold for $L$ and $\Gamma$ in $\mathcal{G}$ and, in addition, the statistical model $\Omega$ is tight. Then $\widehat{H}^* := \sup_{\Pi \in \Delta} \widehat{H}(\Pi)$ is finite, the game $\widehat{\mathcal{G}}^\Omega$ has value $\widehat{H}^*$ and there exists a minimax (robust Bayes) act $\widehat{\zeta}^*$ in $\widehat{\mathcal{G}}^\Omega$ such that*

$$\sup_{\omega \in \Omega} \widehat{L}(\omega, \widehat{\zeta}^*) = \inf_{\zeta \in \mathcal{Z}} \sup_{\omega \in \Omega} \widehat{L}(\omega, \zeta) = \sup_{\Pi \in \Delta} \inf_{a \in \mathcal{A}} \widehat{L}(\Pi, a) = \widehat{H}^*. \tag{126}$$

We remark that the convexity requirement on $\Gamma$ in Condition 6.3 will be satisfied automatically, while the finite entropy requirement is likewise guaranteed by Condition 9.1 and the assumed finiteness of $H^*$.

The proof of Theorem A.2 shows that we can take $\widehat{\zeta}^*$ to be Bayes in $\widehat{\mathcal{G}}$ against some law $\widehat{\Pi}^*$ in the weak closure $\overline{\Delta}$ of $\Delta$ (or, equivalently, Bayes in $\mathcal{G}$ against $\widehat{P}^* := P_{\widehat{\Pi}^*}$ in the weak closure $\overline{\Gamma}$ of $\Gamma$). However, in general, if $\Delta$ is not weakly closed, $\widehat{\zeta}^*$ need not be a Bayes act in $\widehat{\mathcal{G}}$ against any prior distribution $\Pi \in \Delta$ (equivalently, not Bayes in $\mathcal{G}$ against any mixture distribution $P_\Pi \in \Gamma$).

On noting that for any reference act $\zeta_0$ the games $\mathcal{G}^\Gamma$ and $\mathcal{G}_0^\Gamma$ induce the same derived game, and using (111), we have the following.

COROLLARY 9.1. *Suppose that there exists $\zeta_0 \in \mathcal{Z}$ such that Conditions 6.1 and 6.3 hold for $L_0$ and $\Gamma$ in the relative game $\mathcal{G}_0^\Gamma$, and, in addition, that $\mathcal{L}$ is tight. Suppose further that $D(P_\omega, \zeta_0)$ is bounded above for $\omega \in \Omega$. Then there exists a minimax (robust Bayes) act $\hat{\zeta}^*$ in the game $\widehat{\mathcal{G}}^\Omega$.*

If the boundedness condition in Corollary 9.1 fails, we shall have

(127) $$\sup_{\omega \in \Omega} \widehat{L}(\omega, \zeta_0) = \sup_{\omega \in \Omega} D(P_\omega, \zeta_0) = \infty.$$

It can thus fail for all $\zeta_0 \in \mathcal{Z}$ only when $\inf_{\zeta \in \mathcal{Z}} \sup_{\omega \in \Omega} \widehat{L}(\omega, \zeta) = \infty$; that is, the upper value of the game $\widehat{\mathcal{G}}^\Omega$ is $\infty$. In this case the game has no value, and any $\zeta \in \mathcal{Z}$ will trivially be minimax in $\widehat{\mathcal{G}}^\Omega$. In the contrary case, we would normally expect to be able to find a suitable $\zeta_0 \in \mathcal{Z}$ to satisfy all the conditions of Corollary 9.1 and thus demonstrate the existence of a robust Bayes act $\widehat{\zeta}^*$ in $\widehat{\mathcal{G}}^\Omega$.

9.6. *Kullback–Leibler loss: the redundancy-capacity theorem.* An important special case arises when the model $\Omega$ is dominated by a $\sigma$-finite measure $\mu$, and the loss function $L$ in $\mathcal{G}$ is given by the logarithmic score (20) with respect to $\mu$. In this case, for any possible choice of $\mu$, the derived loss is just the Kullback–Leibler divergence, $\widehat{L}(\omega, P) \equiv \mathrm{KL}(P_\omega, P)$. We call such a game a *Kullback–Leibler game*. The corresponding derived entropy $\widehat{H}(\Pi)$, as given by (124), becomes the *mutual information*, $I_\Pi(X, \Omega)$, between $X$ and $\Omega$, in their joint distribution generated by the prior distribution $\Pi$ for $\Omega$ [Lindley (1956)]. There has been much research, especially for asymptotic problems, into the existence and properties of a maximin "reference" prior distribution $\Pi$ over $\Omega$ maximizing this mutual information, or of a minimax act (which can be regarded as a distribution $\widehat{P}^* \in \mathcal{M}$ over $\mathcal{X}$) for DM [Bernardo (1979), Berger and Bernardo (1992), Clarke and Barron (1990, 1994), Haussler (1997) and Xie and Barron (2000)].

The following result follows immediately from Corollary 9.1 and Proposition A.1.

THEOREM 9.2. *Suppose that loss on $\Omega \times \mathcal{A}$ is measured by $\widehat{L}(\omega, P) = \mathrm{KL}(P_\omega, P)$, and that the model $\Omega$ is tight. Then there exists a minimax act $\widehat{P}^* \in \mathcal{M}$ for $\widehat{\mathcal{G}}^\Omega$, achieving $\inf_{P \in \mathcal{M}} \sup_{\omega \in \Omega} \mathrm{KL}(P_\omega, P)$. When this quantity is finite it is the value of the game and equals the maximum attainable mutual information, $I^* := \sup_{\Pi \in \Delta} I_\Pi(X, \Omega)$.*

Theorem 9.2, a version of the "redundancy-capacity theorem" of information theory [Gallager (1976), Ryabko (1979), Davisson and Leon-Garcia (1980) and Krob and Scholl (1997)], constitutes the principal result (Lemma 3) of Haussler (1997). Our proof techniques are different, however.

If $I^*$ is achieved for some $\widehat{\Pi}^* \in \Delta$, then $(\widehat{\Pi}^*, \widehat{P}^*)$ is a saddle-point in $\widehat{\mathcal{G}}^\Omega$, whence, since $\widehat{P}^*$ is then Bayes in $\widehat{\mathcal{G}}$ against $\widehat{\Pi}^*$, $\widehat{P}^*$ is the mixture distribution $P_{\widehat{\Pi}^*} = \int P_\omega \, d\widehat{\Pi}^*(\omega)$. Furthermore, since Lemma 4.2 applies in this case, we find that $\widehat{\Pi}^*$ must be supported on the subspace $\Upsilon := \{\omega \in \Omega : \mathrm{KL}(P_\omega, \widehat{P}^*) = I^*\}$. As argued in Section 4.3, for the case of a continuous parameter-space $\widehat{\Pi}^*$ will typically be a discrete distribution. Notwithstanding this, it is known that, for suitably regular problems, as sample size increases this discrete maximin prior converges weakly to the absolutely continuous Jeffreys invariant prior distribution [Bernardo (1979), Clarke and Barron (1994) and Scholl (1998)].

**10. The Pythagorean inequality.** The Kullback–Leibler divergence satisfies a property reminiscent of squared Euclidean distance. This property was called the *Pythagorean property* by Csiszár (1975). The Pythagorean property leads to an interpretation of minimum relative entropy inference as an *information projection* operation. This view has been emphasized by Csiszár and others in various papers [Csiszár (1975, 1991) and Lafferty (1999)]. Here we investigate the Pythagorean property in our more general framework and show how it is intrinsically related to the minimax theorem: essentially, a *Pythagorean inequality* holds for a discrepancy function $D$ if and only if the loss function $L$ on which $D$ is based admits a saddle-point in a suitable restricted game. Below we formally state and prove this; in Section 10.2 we shall give several examples.

Let $\Gamma \subseteq \mathcal{P}$ be a family of distributions over $\mathcal{X}$, and let $\zeta_0$ be a reference act, such that $L(P, \zeta_0)$ is finite for all $P \in \Gamma$ [so that $L_0(P, \zeta)$ is defined for all $P \in \Gamma$, $\zeta \in \mathcal{Z}$]. We impose no further restrictions on $\Gamma$ (in particular, convexity is not required). Consider the relative restricted game $\mathcal{G}_0^\Gamma$, with loss function $L_0(P, a)$, for $P \in \Gamma$, $a \in \mathcal{A}$. We allow randomization over $\mathcal{A}$ but not over $\Gamma$. The entropy function for this game is $H_0(P) = -D(P, \zeta_0)$ and is always nonpositive.

THEOREM 10.1. *Suppose $(P^*, \zeta^*)$ is a saddle-point in $\mathcal{G}_0^\Gamma$. Then for all $P \in \Gamma$,*

(128) $$D(P, \zeta^*) + D(P^*, \zeta_0) \leq D(P, \zeta_0).$$

*Conversely, if* (128) *holds with its right-hand side finite for all $P \in \Gamma$, then $(P^*, \zeta^*)$ is a saddle-point in $\mathcal{G}_0^\Gamma$.*

Proof. Let $H_0^* := H_0(P^*) = -D(P^*, \zeta_0)$. If $(P^*, \zeta^*)$ is a saddle-point in $\mathcal{G}_0^\Gamma$, then $H_0^* = L_0(P^*, \zeta^*)$ and is finite. Also, for all $P \in \Gamma$,

$$(129) \quad L_0(P, \zeta^*) \leq H_0^*.$$

If $H_0(P) = -\infty$, then $D(P, \zeta_0) = \infty$, so that (128) holds trivially. Otherwise, (129) is equivalent to

$$(130) \quad \{L_0(P, \zeta^*) - H_0(P)\} + \{-H_0^*\} \leq \{-H_0(P)\},$$

which is just (128).

Conversely, in the case that $D(P, \zeta_0)$ is finite for all $P \in \Gamma$, (128) implies (129). Also, putting $P = P^*$ in (128) gives $D(P^*, \zeta^*) = 0$, which is equivalent to $\zeta^*$ being Bayes against $P^*$. Moreover, $H(P^*) = D(P^*, \zeta_0)$ is finite. By (44), $(P^*, \zeta^*)$ is a saddle-point in $\mathcal{G}_0^\Gamma$. □

COROLLARY 10.1. *If $S$ is a $\mathcal{Q}$-proper scoring rule and $\Gamma \subseteq \mathcal{Q}$, then in the restricted relative game $\mathcal{G}_0^\Gamma$ having loss $S_0(P, Q)$ (for fixed reference distribution $P_0 \in \mathcal{Q}$), if $(P^*, P^*)$ is a saddle-point (in which case $P^*$ is both maximum entropy and robust Bayes), then for all $P \in \Gamma$,*

$$(131) \quad d(P, P^*) + d(P^*, P_0) \leq d(P, P_0).$$

*Conversely, if (131) holds and $d(P, P_0) < \infty$ for all $P \in \Gamma$, then $(P^*, P^*)$ is a saddle-point in $\mathcal{G}_0^\Gamma$.*

We shall term (128), or its special case (131), the *Pythagorean inequality*.

We deduce from (128), together with $D(P, \zeta_0) = -H_0(P)$, that for all $P \in \Gamma$,

$$(132) \quad H_0(P^*) - H_0(P) \geq D(P, \zeta^*),$$

providing a quantitative strengthening of the maximum relative entropy property, $H_0(P^*) - H_0(P) \geq 0$, of $P^*$. Similarly, (131) yields

$$(133) \quad H_0(P^*) - H_0(P) \geq d(P, P^*).$$

Often we are interested not in the relative game $\mathcal{G}_0^\Gamma$ but in the original game $\mathcal{G}^\Gamma$. The following corollary relates the Pythagorean inequality to this original game:

COROLLARY 10.2. *Suppose that in the restricted game $\mathcal{G}^\Gamma$ there exists an act $\zeta_0 \in \mathcal{Z}$ such that $L(P, \zeta_0) = k \in \mathcal{R}$, for all $P \in \Gamma$ (in particular, this will hold if $\zeta_0$ is neutral). Then, if $(P^*, \zeta^*)$ is a saddle-point in $\mathcal{G}^\Gamma$, (128) holds for all $P \in \Gamma$; the converse holds if $H(P)$ is finite for all $P \in \Gamma$.*

10.1. *Pythagorean equality.* Related work to date has largely confined itself to the case of equality in (128). This has long been known to hold for the Kullback–Leibler divergence of Section 8.5.2 [Csiszár (1975)]. More recently [Jones and Byrne (1990), Csiszár (1991) and Della Pietra, Della Pietra and Lafferty (2002)], it has been shown to hold for a general Bregman divergence under certain additional conditions. This result extends beyond our framework in that it allows for divergences not defined on probability spaces. On the other hand, when we try to apply it to probability spaces as in Section 3.5.4, its conditions are seen to be highly restrictive, requiring not only differentiability but also, for example, that the tangent space $\nabla H(q)$ of $H$ at $q$ should become infinitely steep as $q$ approaches the boundary of the probability simplex. This is not satisfied even for such simple cases as the Brier score: see Section 10.2.1, where we obtain strict inequality in (128).

The following result follows easily on noting that we have equality in (128) if and only if we have it in (129):

THEOREM 10.2. *Suppose $(P^*, \zeta^*)$ is a saddle-point in $\mathcal{G}_0^\Gamma$. If $\zeta^*$ is an equalizer rule in $\mathcal{G}_0^\Gamma$ [i.e., $L_0(P, \zeta^*) = H_0(P^*)$ for all $P \in \Gamma$], then (128) holds with equality for all $P \in \Gamma$. Conversely, if (128) holds with equality, then $L_0(P, \zeta^*) = H_0(P^*)$ for all $P \in \Gamma$ such that $D(P, \zeta_0) < \infty$; in particular, if $D(P, \zeta_0) < \infty$ for all $P \in \Gamma$, $\zeta^*$ is an equalizer rule in $\mathcal{G}_0^\Gamma$.*

Combining Theorem 10.2 with Theorem 7.1(i) or Corollary 7.2 now gives the following:

COROLLARY 10.3. *Let $\Gamma = \Gamma_\tau = \{P \in \mathcal{P} : \mathrm{E}_P\{t(X)\} = \tau\}$. Suppose $(P^*, \zeta^*) := (P_\tau, \zeta_\tau)$ is a saddle-point in $\mathcal{G}_0^\tau$. If either $(P_\tau, \zeta_\tau)$ is a linear pair or $P \ll P_\tau$, then (128) holds with equality.*

10.2. *Examples.* We now illustrate the Pythagorean theorem and its consequences for our running examples.

10.2.1. *Brier score.* Let $\mathcal{X}$ be finite. As remarked in Section 8.5.1, the Brier divergence $d(P, Q)$ between two distributions $P$ and $Q$ is just $\|p - q\|^2$. Let $\Gamma \subseteq \mathcal{P}$ be closed and convex. By Theorem 5.2, we know that there then exists a $P^* \in \Gamma$ such that $(P^*, P^*)$ is a saddle-point in the relative game $\mathcal{G}_0^\Gamma$. Therefore, by Corollary 10.1 we have, for all $P \in \Gamma$,

$$\|p - p^*\|^2 + \|p^* - p_0\|^2 \leq \|p - p_0\|^2, \tag{134}$$

or equivalently,

$$(p - p^*)^{\mathrm{T}}(p^* - p_0) \leq 0. \tag{135}$$

The distribution $P^*$ within $\Gamma$ that maximizes the Brier entropy relative to $P_0$, or equivalently that minimizes the Brier discrepancy to $P_0$, is given by the point closest to $P_0$ in $\Gamma$, that is, the Euclidean projection of $P_0$ onto $\Gamma$. That this distribution is also a saddle-point is reflected in the fact that the angle $\angle(p, p^*, p_0) \geq 90°$ for all $P \in \Gamma$.

Consider again the case $\mathcal{X} = \{-1, 0, 1\}$ and constraint $\mathrm{E}_P(X) = \tau$. For $\tau \in [-2/3, 2/3]$, where (except for the extreme cases) the minimizing point $p_\tau$ is in the interior of the line segment, (135), and so (134), holds with equality for all $P \in \Gamma_\tau$; while for $\tau$ outside this interval, where the minimizing point is at one end of the segment, (135) and (134) hold with strict inequality for all $P \in \Gamma_\tau \setminus \{P_\tau\}$. Note further that in the former case $p_\tau$ is linear; for $\tau \in (-2/3, 2/3)$ $p_\tau$ is in the interior of the simplex, so that $P_\tau$ has full support. Hence, by Theorem 7.1(i) or Corollary 7.2, $p_\tau$ is an equalizer rule. In the latter case $P_\tau$ does not have full support, and indeed the strict inequality in (134) implies by Theorem 10.2 that it cannot be an equalizer rule.

We can also use (135) to investigate the existence of a saddle-point for certain nonconvex $\Gamma$. Thus suppose, for example, that $\Gamma$ is represented in the simplex by a spherical surface. Then the necessary and sufficient condition (135) for a saddle-point will hold for a reference point $p^0$ outside the sphere, but fail for $p^0$ inside. In the latter case Corollary 4.1 does not apply, and the maximum Brier entropy distribution in $\Gamma$ (the point in $\Gamma$ closest to the center of the simplex) will *not* be robust Bayes against $\Gamma$.

10.2.2. *Logarithmic score.* In this case $d(P, Q)$ becomes the Kullback–Leibler divergence $\mathrm{KL}(P, Q)$ ($P, Q \in \mathcal{M}$). This has been intensively studied for the case of mean-value constraints $\Gamma_\tau^{\mathcal{M}} = \{P \in \mathcal{M} : \mathrm{E}_P(T) = \tau\}$ ($\tau \in \mathcal{T}^0$), when the Pythagorean property (131) holds with equality [Csiszár (1975)]. By Theorem 10.2 this is essentially equivalent to the equalizer property of the maximum relative entropy density $p_\tau$, as already demonstrated (in a way that even extends to distributions $P \in \Gamma_\tau \setminus \mathcal{M}$) in Section 7.3. (Recall from Section 8.5.2 that in this case the relative entropy, with respect to a reference distribution $P_0$, is simply the ordinary entropy under base measure $P_0$.)

In the simple discrete example studied in Section 7.6.2, the above equalizer property also extended (trivially) to the boundary points $\tau = \pm 1$. Such an extension also holds for more general discrete sample spaces, since the condition of Corollary 7.2 can be shown to apply when $\tau$ is on the boundary of $\mathcal{T}$. So in all such cases the Pythagorean inequality (131) is in fact an equality.

10.2.3. *Zero–one loss.* For the case $\mathcal{X} = \{-1, 0, 1\}$ and constraint $\mathrm{E}_P(X) = \tau$, with $\zeta_0$ uniform on $\mathcal{X}$, we have $H_0(P) = H(P) - 1 + 1/N$, and then (132)

(equivalent to both the Pythagorean and the saddle-point property) asserts: for all $P \in \Gamma_\tau$,

$$H(P_\tau) - H(P) \geq D(P, \zeta_\tau). \tag{136}$$

Using (25) and (120), (136) becomes

$$p_{\tau,\max} \leq \sum p(x)\, \zeta_\tau(x). \tag{137}$$

This can be confirmed for the specifications of $P_\tau$ and $\zeta_\tau$ given in Tables 2 and 3. Specifically, for $0 \leq \tau < \frac{1}{2}$, both sides of (137) are $(1+\tau)/3$ (the equality confirming that in this case we have an equalizer rule), while, for $\frac{1}{2} < \tau \leq 1$, (137) becomes $\tau \leq p(1)$, which holds since $\tau = p(1) - p(-1)$ (in particular we have strict inequality, and hence do not have an equalizer rule, unless $\tau = 1$). For $\tau = \frac{1}{2}$, we calculate $\sum p(x)\zeta_\tau(x) - p_{\tau,\max} = (1 - 3a)p(-1)$, which is nonnegative since $a \leq 1/3$, so verifying the Pythagorean inequality, and hence the robust Bayes property of $\zeta_{1/2} = (0, a, 1-a)$ for $a \leq \frac{1}{3}$—although this will be an equalizer rule only for $a = \frac{1}{3}$. Similar results hold when $-1 \leq \tau < 0$.

## 11. Conclusions and further work.

11.1. *What has been achieved.* In this paper we started by interpreting the Shannon entropy of a distribution $P$ as the smallest expected logarithmic loss a DM can achieve when the data are distributed according to $P$. We showed how this interpretation (a) allows for a reformulation of the maximum entropy procedure as a robust Bayes procedure and (b) can be generalized to supply a natural extension of the concept of entropy to arbitrary loss functions. Both these ideas were already known. Our principal novel contribution lies in the combination of the two: the generalized entropies typically still possess a minimax property, and therefore maximum generalized entropy can again be justified as a robust Bayes procedure. For some simple decision problems, as in Section 5, this result is based on an existing minimax theorem due to Ferguson (1967); see the Appendix, Section A.1. For others, as in Section 6, we need more general results, such as Lemma A.1, which uses a (so far as we know) novel proof technique.

We have also considered in detail in Section 7 the special minimax results available when the constraints have the form of known expectations for certain quantities. Arising out of this is our second novel contribution: an extension of the usual definition of "exponential family" to a more general decision framework, as described in Section 7.4. We believe that this extension holds out the promise of important new general statistical theory, such as variations on the concept of sufficiency.

Our third major contribution lies in relating the above theory to the problem of minimizing a discrepancy between distributions. This in turn leads to two further results: in Section 9.5 we generalize Haussler's minimax theorem for the Kullback–Leibler divergence to apply to arbitrary discrepancies; in Section 10 we demonstrate the equivalence between the existence of a saddle-point and a "Pythagorean inequality."

11.2. *Possible developments.*  We end by discussing some possible extensions of our work.

11.2.1. *Moment inequalities.*  As an extension to the moment equalities discussed in Section 7, one may consider robust Bayes problems for moment inequalities, of the form $\Gamma = \{P : E_P(T) \in A\}$, where $A$ is a general (closed, convex) subset of $\mathcal{R}^k$. A direct approach to (39) is complicated by the combination of inner maximization and outer minimization [Noubiap and Seidel (2001)]. Replacement of this problem by a single maximization of entropy over $\Gamma$ could well simplify analysis.

11.2.2. *Nonparametric robust Bayes.*  Much of robust Bayes analysis involves "nonparametric" families $\Gamma$: for example, we might have a reference distribution $P_0$, but, not being sure of its accurate specification, wish to guard against any $P$ in the "$\varepsilon$-neighborhood" of $P_0$, that is, $\{P_0 + c(P - P_0) : |c| \leq \varepsilon, P \text{ arbitrary}\}$. Such a set being closed and convex, a saddle-point will typically exist, and then we can again, in principle, find the robust Bayes act by maximizing the generalized entropy. However, in general it may not be easy to determine or describe the solution to this problem.

11.2.3. *Other generalizations of entropy and entropy optimization problems.*  It would be interesting to make connections between the generalized entropies and discrepancies defined in this text and the several other generalizations of entropy and relative entropy which exist in the literature. Two examples are the Rényi entropies [Rényi (1961)] and the family of entropies based on expected Fisher information considered by Borwein, Lewis and Noll (1996).

Finally, very recently, Harremoës and Topsøe [Topsøe (2002) and Harremoës and Topsøe (2002)] have proposed a generalization of Topsøe's original minimax characterization of entropy [Topsøe (1979)]. They show that a whole range of entropy-related optimization problems can be interpreted from a minimax perspective. While Harremoës and Topsøe's results are clearly related to ours, the exact relation remains a topic of further investigation.

## APPENDIX: PROOFS OF MINIMAX THEOREMS

We first prove Theorem 6.1, which can be used for loss functions that are bounded from above, and Theorem 6.2, which relates saddle-points to differentiability of the entropy. We then prove a general lemma, Lemma A.1, which can be used for unbounded loss functions but imposes other restrictions. This lemma is used to prove Theorem 6.3. Next we demonstrate a general result, Theorem A.2, which implies Theorem 6.4. Finally we prove Lemma 9.2.

**A.1. Theorem 6.1:** $L$ **upper-bounded,** $\Gamma$ **closed and tight.** The following result follows directly from Theorem 2 of Ferguson [(1967), page 85].

THEOREM A.1. *Consider a game* $(\mathcal{X}, \mathcal{A}, L)$. *Suppose that $L$ is bounded below and that there is a topology on $\mathcal{Z}$, the space of randomized acts, such that the following hold:*

(i) $\mathcal{Z}$ *is compact.*
(ii) $L: \mathcal{X} \times \zeta \to \mathcal{R}$ *is lower semicontinuous in $\zeta$ for all $x \in \mathcal{X}$.*

*Then the game has a value, that is,* $\sup_{P \in \mathcal{P}} \inf_{a \in \mathcal{A}} L(P, a) = \inf_{\zeta \in \mathcal{Z}} \sup_{x \in \mathcal{X}} L(x, \zeta)$. *Moreover, a minimax $\zeta$, attaining $\inf_{\zeta \in \mathcal{Z}} \sup_{x \in \mathcal{X}} L(x, \zeta)$, exists.*

Note that $\mathcal{Z}$ could be any convex set. By symmetry considerations, we thus have the following.

COROLLARY A.1. *Consider a game* $(\Gamma, \mathcal{A}, L)$. *Suppose that $L$ is bounded above and there is a topology on $\Gamma$ such that the following hold:*

(i) $\Gamma$ *is convex and compact.*
(ii) $L: \Gamma \times \mathcal{A} \to \mathcal{R}$ *is upper semicontinuous in $P$ for all $a \in \mathcal{A}$.*

*Then the game has a value, that is,* $\inf_{\zeta \in \mathcal{Z}} \sup_{x \in \mathcal{X}} L(x, \zeta) = \sup_{P \in \Gamma} \inf_{a \in \mathcal{A}} L(P, a)$. *Moreover, a maximin $P$, attaining $\sup_{P \in \Gamma} \inf_{a \in \mathcal{A}} L(P, a)$, exists.*

PROOF OF THEOREM 6.1. Since $\Gamma$ is tight and weakly closed, by Prohorov's theorem [Billingsley (1999), Theorem 5.1] it is weakly compact. Also, under the conditions imposed $L(P, a)$ is, for each $a \in \mathcal{A}$, upper semicontinuous in $P$ in the weak topology [Stroock (1993), Theorem 3.1.5(v)]. Theorem 6.1 now follows from Corollary A.1. □

**A.2. Theorems 6.2 and 6.3:** $L$ **unbounded,** $\sup H(P)$ **achieved.** Throughout this section, we assume that $\Gamma$ is convex and that $H^* := \sup_{P \in \Gamma} H(P)$ is finite and is achieved for some $P^* \in \Gamma$ admitting a not necessarily unique Bayes act $\zeta^*$.

To prove that $(P^*, \zeta^*)$ is a saddle-point, it is sufficient to show that $L(P, \zeta^*) \leq L(P^*, \zeta^*) = H^*$ for all $P \in \Gamma$.

PROOF OF THEOREM 6.2. By Lemma 3.2, $L(P, \zeta^*)$ and $L(P_0, \zeta^*)$ are finite, and $f(\lambda) := L(Q_\lambda, \zeta^*)$ is linear in $\lambda \in [0, 1]$. Also, $f(\lambda) \geq H(Q_\lambda)$ for all $\lambda$ and $f(\lambda^*) = H(Q_{\lambda^*}) = H^*$. Thus $f(\lambda)$ must coincide with the tangent to the curve $H(Q_\lambda)$ at $\lambda = \lambda^*$. It follows that

$$(138) \quad L(P, \zeta^*) = f(1) = H^* + (1 - \lambda)\left\{\left(\frac{d}{d\lambda}\right)H(Q_\lambda)\right\}_{\lambda = \lambda^*}.$$

However,

$$\left\{\left(\frac{d}{d\lambda}\right)H(Q_\lambda)\right\}_{\lambda = \lambda^*} = \lim_{\lambda \downarrow \lambda^*}\frac{H(Q_\lambda) - H^*}{\lambda - \lambda^*} \leq 0,$$

since $H(Q_\lambda) \leq H^*$ for $\lambda > \lambda^*$. We deduce $L(P, \zeta^*) \leq H^*$. □

NOTE. If $P_0$ in the statement of Theorem 6.2 can be chosen to be in $\Gamma$, then we further have $H(Q_\lambda) \leq H^*$ for $\lambda < \lambda^*$, which implies $\{(d/d\lambda)H(Q_\lambda)\}_{\lambda = \lambda^*} = 0$, and hence $L(P, \zeta^*) = H^*$. In particular, if this can be done for all $P \in \Gamma$ (i.e., $P^*$ is an "algebraically interior" point of $\Gamma$), then $\zeta^*$ will be an equalizer rule.

From this point on, for any $P \in \Gamma$, $\lambda \in [0, 1]$ we write $P_\lambda := \lambda P + (1 - \lambda)P^*$. Then, since we are assuming $\Gamma$ convex, $P_\lambda \in \Gamma$.

LEMMA A.1. *Suppose Conditions 6.3 and 6.4 hold. Let $\zeta_\lambda$ be Bayes against $P_\lambda$ (in particular, $\zeta^* := \zeta_0$ is Bayes against $P^*$, and $\zeta_1$ is Bayes against $P$). Then*

$$(139) \quad L(P, \zeta_\lambda) - L(P^*, \zeta_\lambda) = \frac{H(P_\lambda) - L(P^*, \zeta_\lambda)}{\lambda}$$

$$(140) \quad \leq 0$$

*($0 < \lambda < 1$). Moreover, $\lim_{\lambda \downarrow 0} L(P^*, \zeta_\lambda)$ and $\lim_{\lambda \downarrow 0} L(P, \zeta_\lambda)$ both exist as finite numbers, and*

$$(141) \quad \lim_{\lambda \downarrow 0} L(P^*, \zeta_\lambda) = H^*.$$

PROOF. First note that, since $H(P_\lambda) = L(P_\lambda, \zeta_\lambda)$ is finite, by Lemma 3.2 both $L(P, \zeta_\lambda)$ and $L(P^*, \zeta_\lambda)$ are finite for $0 < \lambda < 1$. Also by Lemma 3.2, for all $\zeta \in \mathcal{Z}$, $L(P_\lambda, \zeta)$ is, when finite, a linear function of $\lambda \in [0, 1]$. Then

$$\lambda L(P, \zeta) + (1 - \lambda)L(P^*, \zeta) = L(P_\lambda, \zeta)$$

$$(142) \quad \geq H(P_\lambda) = L(P_\lambda, \zeta_\lambda)$$

$$(143) \quad = \lambda L(P, \zeta_\lambda) + (1 - \lambda)L(P^*, \zeta_\lambda).$$

On putting $\zeta = \zeta_\lambda$ we have equality in (142); then rearranging yields (139), and (140) follows from $L(P^*, \zeta_\lambda) \geq H^*$ and $H(P_\lambda) \leq H^*$.

For general $\zeta \in \mathcal{Z}$ we obtain (when all terms are finite)

$$(144) \qquad \lambda\{L(P,\zeta_\lambda) - L(P,\zeta)\} \leq (1-\lambda)\{L(P^*,\zeta) - L(P^*,\zeta_\lambda)\}.$$

Put $\zeta = \zeta_1$, so that $L(P,\zeta_1) = H(P)$ is finite, and first suppose that $L(P^*,\zeta_1)$ is finite. Then the left-hand side of (144) is nonnegative, and so $L(P^*,\zeta_1) \geq L(P^*,\zeta_\lambda)$ ($0 \leq \lambda \leq 1$)—which inequality clearly also holds if $L(P^*,\zeta_1) = \infty$. An identical argument can be applied on first replacing $\zeta_1$ by $\zeta_{\lambda'}$ ($0 < \lambda' < 1$), and we deduce that $L(P^*,\zeta_{\lambda'}) \geq L(P^*,\zeta_\lambda)$ ($0 \leq \lambda \leq \lambda' \leq 1$). That is to say, $L(P^*,\zeta_\lambda)$ is a nondecreasing function of $\lambda$ on $[0,1]$. It follows that

$$(145) \qquad \lim_{\lambda \downarrow 0} L(P^*,\zeta_\lambda) \geq L(P^*,\zeta_0) = H^*.$$

A parallel argument, interchanging the roles of $P^*$ and $P$, shows that $L(P,\zeta_\lambda)$ is nonincreasing in $\lambda \in [0,1]$. Since, by (140), for all $\lambda \in (0, 0.5]$, $L(P,\zeta_\lambda) \leq L(P^*,\zeta_\lambda) \leq L(P^*,\zeta_{0.5}) < \infty$, it follows that $\lim_{\lambda \downarrow 0} L(P,\zeta_\lambda)$ exists and is finite.

Since $P^*$ maximizes entropy over $\Gamma$,

$$(146) \qquad \begin{aligned} H(P^*) - L(P^*,\zeta_\lambda) &\geq H(P_\lambda) - L(P^*,\zeta_\lambda) \\ &= \lambda\{L(P,\zeta_\lambda) - L(P^*,\zeta_\lambda)\}, \end{aligned}$$

by (143). On noting $L(P^*,\zeta_\lambda) \leq L(P^*,\zeta_1)$ since $L(P^*,\zeta_\lambda)$ is nondecreasing, and using $L(P,\zeta_\lambda) \geq H(P)$, (146) implies $H^* - L(P^*,\zeta_\lambda) \geq \lambda\{H(P) - L(P^*,\zeta_1)\}$. If $L(P^*,\zeta_1) < \infty$, then letting $\lambda \downarrow 0$ we obtain $H^* \geq \lim_{\lambda \downarrow 0} L(P^*,\zeta_\lambda)$, which, together with (145), establishes (141). Otherwise, noting that $L(P^*,\zeta_{0.5}) < \infty$, we can repeat the argument with $P$ replaced by $P_{0.5}$. $\square$

COROLLARY A.2.

$$(147) \qquad \lim_{\lambda \downarrow 0} L(P,\zeta_\lambda) - H^* = \lim_{\lambda \downarrow 0} \frac{H(P_\lambda) - L(P^*,\zeta_\lambda)}{\lambda}.$$

COROLLARY A.3 (Condition for existence of a saddle-point). $L(P,\zeta^*) \leq H(P^*)$ *if and only if*

$$(148) \qquad \lim_{\lambda \downarrow 0} \frac{H(P_\lambda) - L(P^*,\zeta_\lambda)}{\lambda} \leq \lim_{\lambda \downarrow 0} L(P,\zeta_\lambda) - L(P,\zeta^*).$$

PROOF OF THEOREM 6.3. The conditions of Lemma A.1 are satisfied. By Corollary A.3 and (140), we see that it is sufficient to prove that, for all $P \in \Gamma$,

$$(149) \qquad 0 \leq \lim_{\lambda \downarrow 0} L(P,\zeta_\lambda) - L(P,\zeta^*).$$

However, (149) is implied by Condition 6.1. $\square$

**A.3. If $\sup_{P \in \Gamma} H(P)$ is not achieved.** In some cases $\sup_{P \in \Gamma} H(P)$ may not be achieved in $\Gamma$ [Topsøe (1979)]. We might then think of enlarging $\Gamma$ to, say, its weak closure $\overline{\Gamma}$. However, this can be much bigger than $\Gamma$. For example, for uncountable $\mathcal{X}$, the weak closure of a set, all of whose members are absolutely continuous with respect to $\mu$, typically contains distributions that are not. Then Theorem 6.3 may not be applicable.

EXAMPLE A.1. Consider the logarithmic score, as in Section 3.5.2, with $\mathcal{X} = \mathcal{R}$ and $\mu$ Lebesgue measure, and let $\Gamma = \{P : P \ll \mu, \mathrm{E}(X) = 0, \mathrm{E}(X^2) = 1\}$. Then $\overline{\Gamma}$ contains the distribution $P$ with $P(X=1) = P(X=-1) = 1/2$, for which $H(P) = -\infty$. There is no Bayes act against this $P$.

This example illustrates that, in case $\sup_{P \in \Gamma} H(P)$ is not achieved [for an instance of this, see Cover and Thomas (1991), Chapter 9], we cannot simply take its closure and then apply Theorem 6.3, since Condition 6.3 could still be violated.

The following theorem, which in turn implies Theorem 6.4 of Section 6, shows that the game $(\Gamma, \mathcal{A}, L)$ will often have a value even when $\Gamma$ is not weakly closed. We need to impose an additional condition:

CONDITION A.1. Every sequence $(Q_n)$ of distributions in $\Gamma$ such that $H(Q_n)$ converges to $H^*$ has a weak limit point in $\mathcal{P}_0$.

THEOREM A.2. *Suppose Conditions 6.1, 6.3 and A.1 hold. Then there exists $\zeta^* \in \mathcal{Z}$ such that*

$$(150) \quad \sup_{P \in \Gamma} L(P, \zeta^*) = \inf_{\zeta \in \mathcal{Z}} \sup_{P \in \Gamma} L(P, \zeta) = \sup_{P \in \Gamma} \inf_{a \in \mathcal{A}} L(P, a) = H^*.$$

*In particular, the game $\mathcal{G}^\Gamma$ has value $H^*$, and $\zeta^*$ is robust Bayes against $\Gamma$.*

PROOF. Let $(Q_n)$ be a sequence in $\Gamma$ such that $H(Q_n)$ converges to $H^*$. In particular, $(H(Q_n))$ is bounded below. On choosing a subsequence if necessary, we can suppose by Condition A.1 that $(Q_n)$ has a weak limit $P^*$, and further that for all $n$ $H^* - H(Q_n) < 1/n$. By Condition 6.1, $P^*$ has a Bayes act $\zeta^*$.

Now pick any $P \in \Gamma$. We will show that $L(P, \zeta^*) \leq H^*$. First fix $n$ and define $R^n_\lambda := \lambda P + (1-\lambda) Q_n$, $H^n_\lambda := H(R^n_\lambda)$ $(0 \leq \lambda \leq 1)$. In particular, $R^n_0 = Q_n$, $R^n_1 = P$. Then $R^n_\lambda \in \Gamma$, with Bayes act $\zeta^n_\lambda$, say. We have $H^n_\lambda = L(R^n_\lambda, \zeta^n_\lambda) = \lambda L(P, \zeta^n_\lambda) + (1 - \lambda) L(R^n_0, \zeta^n_\lambda)$, while $H^n_0 \leq L(R^n_0, \zeta^n_\lambda)$. It follows that

$$(151) \quad L(P, \zeta^n_\lambda) \leq H^n_0 + (H^n_\lambda - H^n_0)/\lambda.$$

Since $H^n_0 = H(Q_n) > H^* - 1/n$ and $H^n_0, H^n_\lambda \leq H^*$, we obtain

$$(152) \quad L(P, \zeta^n_{1/\sqrt{n}}) \leq H^* + 1/n + 1/\sqrt{n}.$$

Now with $Q'_n := R^n_{1/\sqrt{n}}$, $(Q'_n)$ converges weakly to $P^*$. Moreover, $H(Q'_n) \geq (1/\sqrt{n})H(P) + (1 - 1/\sqrt{n})H(Q_n)$ is bounded below. On applying Condition 6.1 to $Q'_n$, and using (152), we deduce

(153) $$L(P, \zeta^*) \leq H^*.$$

It now follows that

(154) $$\inf_{\zeta \in \mathcal{Z}} \sup_{P \in \Gamma} L(P, \zeta) \leq \sup_{P \in \Gamma} L(P, \zeta^*) \leq H^*.$$

However,

(155) $$H^* = \sup_{P \in \Gamma} \inf_{a \in \mathcal{A}} L(P, a) = \sup_{P \in \Gamma} \inf_{\zeta \in \mathcal{Z}} L(P, \zeta) \leq \inf_{\zeta \in \mathcal{Z}} \sup_{P \in \Gamma} L(P, \zeta),$$

where the the second equality follows from Proposition 3.1 and the third inequality is standard. Together, (154) and (155) imply the theorem. □

PROOF OF THEOREM 6.4. If $\Gamma$ is tight, then by Prohorov's theorem any sequence $(Q_n)$ in $\Gamma$ must have a weak limit point, so that, in particular, Condition A.1 holds. □

It should be noted that, for $P^*$ appearing in the above proof, we may have $H(P^*) \neq H^*$. In the case of Shannon entropy, we have $H(P^*) \leq H^*$; a detailed study of the case of strict inequality has been carried out by Harremoës and Topsøe (2001).

We now show, following Csiszár (1975) and Topsøe (1979), that the conditions of Theorem A.2 are satisfied by the logarithmic score. We take $L = S$, the logarithmic score (20) defined with respect to a measure $\mu$. This is $\mathcal{M}$-strictly proper, where $\mathcal{M}$ is the set of all probability distributions absolutely continuous with respect to $\mu$.

PROPOSITION A.1. *Conditions A.1 and 6.2 are satisfied for the logarithmic score $S$ relative to a measure $\mu$ if* either *of the following holds:*

(i) $\mu$ *is a probability measure and* $\mathcal{Q} = \mathcal{M}$;
(ii) $\mathcal{X}$ *is countable,* $\mu$ *is counting measure and* $\mathcal{Q} = \{P \in \mathcal{P} : H(P) < \infty\}$.

PROOF. To show Condition A.1, under either (i) or (ii), let $(Q_n)$ be a sequence of distributions in $\Gamma$ such that $H(Q_n)$ converges to $H^*$. Given $\varepsilon > 0$, choose $N$ such that, for $n \geq N$, $H^* - H(Q_n) < \varepsilon$. Then for $n, m \geq N$, on applying (104) we have

(156) $$\begin{aligned} H^* &\geq H\{\tfrac{1}{2}(Q_n + Q_m)\} \\ &= \tfrac{1}{2}[H(Q_n) + H(Q_m) + \mathrm{KL}\{Q_n, \tfrac{1}{2}(Q_n + Q_m)\} \\ &\quad + \mathrm{KL}\{Q_m, \tfrac{1}{2}(Q_n + Q_m)\}] \\ &\geq H^* - \varepsilon + \tfrac{1}{16}\|Q_n - Q_m\|^2, \end{aligned}$$

where $\|\cdot\|$ denotes total variation and the last inequality is an application of Pinsker's inequality $\mathrm{KL}(P_1, P_2) \geq (1/4)\|P_1 - P_2\|^2$ [Pinsker (1964)]. That is, $n, m \geq N \Rightarrow \|Q_n - Q_m\|^2 \leq 16\varepsilon$, so that $(Q_n)$ is a Cauchy sequence under $\|\cdot\|$. Since the total variation metric is complete, $(Q_n)$ has a limit $Q$ in the uniform topology, which is then also a weak limit [Stroock (1993), Section 3.1]. This shows Condition A.1.

To demonstrate Condition 6.2, suppose $Q_n \in \mathcal{Q}$, $H(Q_n) \geq K > -\infty$ for all $n$, and $(Q_n)$ converges weakly to some distribution $Q_0 \in \mathcal{P}_0$. By Posner (1975), Theorem 1, $\mathrm{KL}(P, Q)$ is jointly weakly lower semicontinuous in both arguments. In case (i), the entropy $H(P) \equiv -\mathrm{KL}(P, \mu)$ is thus upper semicontinuous in $P \in \mathcal{P}$, and it follows that $0 \geq H(Q_0) \geq K > -\infty$, which implies $Q_0 \in \mathcal{M} = \mathcal{Q}$. In case (ii), the entropy function is lower semicontinuous [Topsøe (2001)], whence $0 \leq H(Q_0) < \infty$, and again $Q_0 \in \mathcal{Q}$. In either case, the lower semicontinuity of $\mathrm{KL}(P, Q)$ in $Q$ then implies that, for $P \in \mathcal{Q}$, $S(P, Q_0) = \mathrm{KL}(P, Q_0) + H(P) \leq \liminf_{n \to \infty} \{\mathrm{KL}(P, Q_n) + H(P)\} = \liminf_{n \to \infty} S(P, Q_n)$.

□

Theorem A.2 essentially extends the principal arguments and results of Topsøe (1979) to nonlogarithmic loss functions. In such cases it might sometimes be possible to establish the required conditions by methods similar to Proposition A.1, but in general this could require new techniques.

**A.4. Proof of Lemma 9.2.** Suppose Condition 9.1 holds, and Conditions 6.1 and 6.3 hold for $L$ and $\Gamma$ in $\mathcal{G}$. We note that $H(P_\omega)$ is then bounded below by $K$ and above by $H^*$ for $\omega \in \Omega$; for $\Pi \in \Delta$, the integral in (123) and (124) is then bounded by the same quantities.

To show Condition 6.1 holds for $\widehat{L}$ and $\Delta$ in $\widehat{\mathcal{G}}$, let $\Pi_n \in \Delta$, with Bayes act $\zeta_n \in \mathcal{Z}$ in $\widehat{\mathcal{G}}$, be such that $(\widehat{H}(\Pi_n))$ is bounded below and $(\Pi_n)$ converges weakly to $\Pi_0 \in \overline{\Delta}$. Defining $Q_n := P_{\Pi_n}, Q_0 := P_{\Pi_0}$, we then have $Q_n \in \Gamma$, with Bayes act $\zeta_n \in \mathcal{Z}$ in $\mathcal{G}$. Now let $f: \mathcal{X} \to \mathcal{R}$ be bounded and continuous, and define $g: \Theta_0 \to \mathcal{R}$ by $g(\theta) = \mathrm{E}_{P_\theta}\{f(X)\}$. By the definition of weak convergence, the function $g$ is continuous. It follows that $\mathrm{E}_{Q_n}\{f(X)\} = \mathrm{E}_{\Pi_n}\{g(\Theta)\} \to \mathrm{E}_{\Pi_0}\{g(\Theta)\} = \mathrm{E}_{Q_0}\{f(X)\}$. This shows that $(Q_n)$ converges weakly to $Q_0$. Also, by (124) and Condition 9.1, the sequence $(H(Q_n))$ is bounded below. It now follows from Condition 6.1 in $\mathcal{G}^\Gamma$ that $Q_0$ has a Bayes act $\zeta_0$ in $\mathcal{G}$—any such act likewise being Bayes against $\Pi_0$ in $\widehat{\mathcal{G}}$. Also, for an appropriate choice of the Bayes acts $(\zeta_n)$ and $\zeta_0$, $L(P, \zeta_0) \leq \liminf_{n \to \infty} L(P, \zeta_n)$, for all $P \in \Gamma$. By finiteness of the integral in (123) we then obtain $\widehat{L}(\Pi, \zeta_0) \leq \liminf_{n \to \infty} \widehat{L}(\Pi, \zeta_n)$, for all $\Pi \in \Delta$.

We now show that Condition 6.3 holds for $\widehat{L}$ and $\Delta$ in $\widehat{\mathcal{G}}$. First it is clear that $\Delta$ is convex. Since $\Pi \in \Delta$ and $P_\Pi \in \Gamma$ have the same Bayes acts (in their respective games), if $P_\Pi \in \Gamma$ has a Bayes act, then so does $\Pi$. Also, the integral in (123) is bounded as a function of $\Pi$, whence $\widehat{H}(\Pi)$ is finite if $H(P_\Pi)$ is, and $\sup_{\Pi \in \Delta} \widehat{H}(\Pi)$ is finite if $\sup_{P \in \Gamma} H(P)$ is.

CWI Amsterdam  
P.O. Box 94079  
NL-1090 GB Amsterdam  
The Netherlands  
e-mail: pdg@cwi.nl

Department of Statistical Science  
University College London  
Gower Street  
London WC1E 6BT  
United Kingdom  
e-mail: dawid@stats.ucl.ac.uk